\documentclass{article}
\usepackage{graphicx,amsthm}
\usepackage{verbatim,amsmath,amscd,amssymb,color}
\usepackage[mathscr]{eucal}
\usepackage{epic,eepic,multicol}
\numberwithin{equation}{section}
\input xy
\xyoption{all}
\def\m@th{\mathsurround=0pt}

\def\fsquare(#1,#2){
\hbox{\vrule$\hskip-0.4pt\vcenter to #1{\normalbaselines\m@th
\hrule\vfil\hbox to #1{\hfill$\scriptstyle #2$\hfill}\vfil\hrule}$\hskip-0.4pt
\vrule}}

\def\addsquare(#1,#2){\hbox{$
	\dimen1=#1 \advance\dimen1 by -0.8pt
	\vcenter to #1{\hrule height0.4pt depth0.0pt%
	\hbox to #1{%
	\vbox to \dimen1{\vss%
	\hbox to \dimen1{\hss$\scriptstyle~#2~$\hss}%
	\vss}%
	\vrule width0.4pt}%
	\hrule height0.4pt depth0.0pt}$}}

\def\Fsquare(#1,#2){
\hbox{\vrule$\hskip-0.4pt\vcenter to #1{\normalbaselines\m@th
\hrule\vfil\hbox to #1{\hfill$#2$\hfill}\vfil\hrule}$\hskip-0.4pt
\vrule}}

\def\Addsquare(#1,#2){\hbox{$
	\dimen1=#1 \advance\dimen1 by -0.8pt
	\vcenter to #1{\hrule height0.4pt depth0.0pt%
	\hbox to #1{%
	\vbox to \dimen1{\vss%
	\hbox to \dimen1{\hss$~#2~$\hss}%
	\vss}%
	\vrule width0.4pt}%
	\hrule height0.4pt depth0.0pt}$}}

\def\hfourbox(#1,#2,#3,#4){%
	\fsquare(0.3cm,#1)\addsquare(0.3cm,#2)\addsquare(0.3cm,#3)\addsquare(0.3cm,#4)}

\def\Hfourbox(#1,#2,#3,#4){%
	\Fsquare(0.4cm,#1)\Addsquare(0.4cm,#2)\Addsquare(0.4cm,#3)\Addsquare(0.4cm,#4)}

\def\HHfourbox(#1,#2,#3,#4){%
	\Fsquare(0.8cm,#1)\Addsquare(0.8cm,#2)\Addsquare(0.8cm,#3)\Addsquare(0.8cm,#4)}

\def\hthreebox(#1,#2,#3){%
	\fsquare(0.3cm,#1)\addsquare(0.3cm,#2)\addsquare(0.3cm,#3)}

\def\htwobox(#1,#2){%
	\fsquare(0.3cm,#1)\addsquare(0.3cm,#2)}

\def\vfourbox(#1,#2,#3,#4){%
	\normalbaselines\m@th\offinterlineskip
	\vcenter{\hbox{\fsquare(0.3cm,#1)}
	      \vskip-0.4pt
	      \hbox{\fsquare(0.3cm,#2)}	
	      \vskip-0.4pt
	      \hbox{\fsquare(0.3cm,#3)}	
	      \vskip-0.4pt
	      \hbox{\fsquare(0.3cm,#4)}}}

\def\VVfourbox(#1,#2,#3,#4){%
	\normalbaselines\m@th\offinterlineskip
	\vcenter{\hbox{\Fsquare(0.8cm,#1)}
	      \vskip-0.4pt
	      \hbox{\Fsquare(0.8cm,#2)}	
	      \vskip-0.4pt
	      \hbox{\Fsquare(0.8cm,#3)}	
	      \vskip-0.4pt
	      \hbox{\Fsquare(0.8cm,#4)}}}

\def\Vfourbox(#1,#2,#3,#4){%
	\normalbaselines\m@th\offinterlineskip
	\vcenter{\hbox{\Fsquare(0.4cm,#1)}
	      \vskip-0.4pt
	      \hbox{\Fsquare(0.4cm,#2)}	
	      \vskip-0.4pt
	      \hbox{\Fsquare(0.4cm,#3)}	
	      \vskip-0.4pt
	      \hbox{\Fsquare(0.4cm,#4)}}}

\def\vthreebox(#1,#2,#3){%
	\normalbaselines\m@th\offinterlineskip
	\vcenter{\hbox{\fsquare(0.3cm,#1)}
	      \vskip-0.4pt
	      \hbox{\fsquare(0.3cm,#2)}	
	      \vskip-0.4pt
	      \hbox{\fsquare(0.3cm,#3)}}}

\def\vtwobox(#1,#2){%
	\normalbaselines\m@th\offinterlineskip
	\vcenter{\hbox{\fsquare(0.4cm,#1)}
	      \vskip-0.4pt
	      \hbox{\fsquare(0.4cm,#2)}}}

\def\Hthreebox(#1,#2,#3){%
	\Fsquare(0.4cm,#1)\Addsquare(0.4cm,#2)\Addsquare(0.4cm,#3)}

\def\HHthreebox(#1,#2,#3){%
	\Fsquare(0.8cm,#1)\Addsquare(0.8cm,#2)\Addsquare(0.8cm,#3)}

\def\Htwobox(#1,#2){%
	\Fsquare(0.4cm,#1)\Addsquare(0.4cm,#2)}

\def\H6twobox(#1,#2){%
	\Fsquare(0.6cm,#1)\Addsquare(0.6cm,#2)}

\def\HHtwobox(#1,#2){%
	\Fsquare(0.8cm,#1)\Addsquare(0.8cm,#2)}

\def\Vthreebox(#1,#2,#3){%
	\normalbaselines\m@th\offinterlineskip
	\vcenter{\hbox{\Fsquare(0.4cm,#1)}
	      \vskip-0.4pt
	      \hbox{\Fsquare(0.4cm,#2)}	
	      \vskip-0.4pt
	      \hbox{\Fsquare(0.4cm,#3)}}}

\def\Vtwobox(#1,#2){%
	\normalbaselines\m@th\offinterlineskip
	\vcenter{\hbox{\Fsquare(0.6cm,#1)}
	      \vskip-0.4pt
	      \hbox{\Fsquare(0.6cm,#2)}}}

\def\twoone(#1,#2,#3){%
	\normalbaselines\m@th\offinterlineskip
	\vcenter{\hbox{\htwobox({#1},{#2})}
	      \vskip-0.4pt
	      \hbox{\fsquare(0.3cm,#3)}}}

\def\Twoone(#1,#2,#3){%
	\normalbaselines\m@th\offinterlineskip
	\vcenter{\hbox{\Htwobox({#1},{#2})}
	      \vskip-0.4pt
	      \hbox{\Fsquare(0.4cm,#3)}}}

\def\TTwoone(#1,#2,#3){%
	\normalbaselines\m@th\offinterlineskip
	\vcenter{\hbox{\H6twobox({#1},{#2})}
	      \vskip-0.4pt
	      \hbox{\Fsquare(0.6cm,#3)}}}

\def\threeone(#1,#2,#3,#4){%
	\normalbaselines\m@th\offinterlineskip
	\vcenter{\hbox{\hthreebox({#1},{#2},{#3})}
	      \vskip-0.4pt
	      \hbox{\fsquare(0.3cm,#4)}}}

\def\Threeone(#1,#2,#3,#4){%
	\normalbaselines\m@th\offinterlineskip
	\vcenter{\hbox{\Hthreebox({#1},{#2},{#3})}
	      \vskip-0.4pt
	      \hbox{\Fsquare(0.4cm,#4)}}}

\def\Threetwo(#1,#2,#3,#4,#5){%
	\normalbaselines\m@th\offinterlineskip
	\vcenter{\hbox{\Hthreebox({#1},{#2},{#3})}
	      \vskip-0.4pt
	      \hbox{\Htwobox({#4},{#5})}}}

\def\threetwo(#1,#2,#3,#4,#5){%
	\normalbaselines\m@th\offinterlineskip
	\vcenter{\hbox{\hthreebox({#1},{#2},{#3})}
	      \vskip-0.4pt
	      \hbox{\htwobox({#4},{#5})}}}

\def\threethree(#1,#2,#3,#4,#5,#6){%
	\normalbaselines\m@th\offinterlineskip
	\vcenter{\hbox{\hthreebox({#1},{#2},{#3})}
	      \vskip-0.4pt
	      \hbox{\hthreebox({#4},{#5},{#6})}}}

\def\twotwo(#1,#2,#3,#4){%
	\normalbaselines\m@th\offinterlineskip
	\vcenter{\hbox{\htwobox({#1},{#2})}
	      \vskip-0.4pt
	      \hbox{\htwobox({#3},{#4})}}}

\def\Twotwo(#1,#2,#3,#4){%
	\normalbaselines\m@th\offinterlineskip
	\vcenter{\hbox{\Htwobox({#1},{#2})}
	      \vskip-0.4pt
	      \hbox{\Htwobox({#3},{#4})}}}

\def\TTwotwo(#1,#2,#3,#4){%
	\normalbaselines\m@th\offinterlineskip
	\vcenter{\hbox{\H6twobox({#1},{#2})}
	      \vskip-0.4pt
	      \hbox{\H6twobox({#3},{#4})}}}

\def\twooneone(#1,#2,#3,#4){%
	\normalbaselines\m@th\offinterlineskip
	\vcenter{\hbox{\htwobox({#1},{#2})}
	      \vskip-0.4pt
	      \hbox{\fsquare(0.3cm,#3)}
	      \vskip-0.4pt
	      \hbox{\fsquare(0.3cm,#4)}}}

\def\Twooneone(#1,#2,#3,#4){%
	\normalbaselines\m@th\offinterlineskip
	\vcenter{\hbox{\Htwobox({#1},{#2})}
	      \vskip-0.4pt
	      \hbox{\Fsquare(0.4cm,#3)}
	      \vskip-0.4pt
	      \hbox{\Fsquare(0.4cm,#4)}}}

\def\Twotwoone(#1,#2,#3,#4,#5){%
	\normalbaselines\m@th\offinterlineskip
	\vcenter{\hbox{\Htwobox({#1},{#2})}
	      \vskip-0.4pt
	      \hbox{\Htwobox({#3},{#4})}
              \vskip-0.4pt
	      \hbox{\Fsquare(0.4cm,#5)}}}

\def\twotwoone(#1,#2,#3,#4,#5){%
	\normalbaselines\m@th\offinterlineskip
	\vcenter{\hbox{\htwobox({#1},{#2})}
	      \vskip-0.4pt
	      \hbox{\htwobox({#3},{#4})}
              \vskip-0.4pt
	      \hbox{\fsquare(0.3cm,#5)}}}

\def\twotwotwo(#1,#2,#3,#4,#5,#6){%
	\normalbaselines\m@th\offinterlineskip
	\vcenter{\hbox{\htwobox({#1},{#2})}
	      \vskip-0.4pt
	      \hbox{\htwobox({#3},{#4})}
              \vskip-0.4pt
	      \hbox{\htwobox({#5},{#6})}}}

\def\threetwoone(#1,#2,#3,#4,#5,#6){%
	\normalbaselines\m@th\offinterlineskip
	\vcenter{\hbox{\hthreebox({#1},{#2},{#3})}
	      \vskip-0.4pt
	      \hbox{\htwobox({#4},{#5})}
              \vskip-0.4pt
	      \hbox{\fsquare(0.3cm,#6)}}}

\def\TThreetwoone(#1,#2,#3,#4,#5,#6){%
	\normalbaselines\m@th\offinterlineskip
	\vcenter{\hbox{\HHthreebox({#1},{#2},{#3})}
	      \vskip-0.4pt
	      \hbox{\HHtwobox({#4},{#5})}
              \vskip-0.4pt
	      \hbox{\fsquare(0.8cm,#6)}}}

\def\Threeoneone(#1,#2,#3,#4,#5){%
	\normalbaselines\m@th\offinterlineskip
	\vcenter{\hbox{\Hthreebox({#1},{#2},{#3})}
	      \vskip-0.4pt
	      \hbox{\Fsquare(0.4cm,#4)}
              \vskip-0.4pt
	      \hbox{\Fsquare(0.4cm,#5)}}}

\def\threeoneone(#1,#2,#3,#4,#5){%
	\normalbaselines\m@th\offinterlineskip
	\vcenter{\hbox{\hthreebox({#1},{#2},{#3})}
	      \vskip-0.4pt
	      \hbox{\fsquare(0.3cm,#4)}
              \vskip-0.4pt
	      \hbox{\fsquare(0.3cm,#5)}}}

\def\FFourtwoone(#1,#2,#3,#4,#5,#6,#7){%
	\normalbaselines\m@th\offinterlineskip
	\vcenter{\hbox{\HHfourbox({#1},{#2},{#3},{#4})}
	      \vskip-0.4pt
	      \hbox{\HHtwobox({#5},{#6})}
              \vskip-0.4pt
	      \hbox{\fsquare(0.8cm,#7)}}}

\def\a{\fsquare(0.3cm){1}\addsquare(0.3cm)(2)\addsquare(0.3cm)(3)}

\def\b{\hbox{%
	\normalbaselines\m@th\offinterlineskip
	\vcenter{\hbox{\fsquare(0.3cm){2}}\vskip-0.4pt\hbox{\fsquare(0.3cm){2}}}}}

\def\c{\hbox{\normalbaselines\m@th\offinterlineskip%
	\vcenter{\hbox{\a}\vskip-0.4pt\hbox{\b}}}}

\def\ffsquare#1{%
	\fsquare(0.4cm,\hbox{#1})}

\def\naga{%
	\hbox{$\vcenter to 0.4cm{\normalbaselines\m@th
	\hrule\vfil\hbox to 1.2cm{\hfill$\cdots$\hfill}\vfil\hrule}$}}

\def\vnaga{\normalbaselines\m@th\baselineskip0pt\offinterlineskip%
	\vrule\vbox to 1.2cm{\vskip7pt\hbox to \dimen1{$\hfil\vdots\hfil$}\vfil}\vrule}

\def\dvbox{\hbox{\normalbaselines\m@th\baselineskip0pt\offinterlineskip\vbox{%
	  \hbox{$\ffsquare 1$}\vskip-0.4pt\hbox{$\vnaga$}\vskip-0.4pt\hbox{$\ffsquare N$}}}}

\def\sq(#1){\fsquare(0.4cm,#1)}
\def\Sq(#1){\fsquare(0.5cm,#1)}
\def\SSq(#1){\fsquare(0.9cm,#1)}

\def\mapright#1{\smash{\mathop{\longrightarrow}\limits^{#1}}}

\font\germ=eufm10

\newcommand{\cD}{{\mathcal D}}

\newcommand{\cV}{{\mathcal V}}

\newcommand{\frg}{\mathfrak g}

\newcommand{\frt}{\mathfrak t}

\newcommand{\bbA}{\mathbb A}

\newcommand{\bbC}{\mathbb C}

\newcommand{\bbP}{\mathbb P}
\newcommand{\bbN}{\mathbb N}

\newcommand{\bbZ}{\mathbb Z}

\newcommand{\wt}{\mathrm{wt}}

\newcommand{\lbr}{\begin{bmatrix}}
\newcommand{\rbr}{\end{bmatrix}}
\newcommand{\forb}{\bigcirc\kern-2.8ex \because}
\newcommand{\forbb}{\bigcirc\kern-3.0ex \because}
\newcommand{\forbbb}{\bigcirc\kern-3.1ex \because}
\newcommand{\cd}{commutative diagram }

\newcommand\C{\mathbb C}

\newcommand\Z{{\mathbb Z}}

\def\aa{{\bf a}}
\def\bb{{\bf b}}
\def\ge{\frg}

\def\al{\alpha}

\def\beneme{\begin{enumerate}}
\def\beq{\begin{equation}}
\def\beqn{\begin{eqnarray}}
\def\beqnn{\begin{eqnarray*}}

\def\bbra#1,#2,#3{\left\{\begin{array}{c}\hspace{-5pt}
#1;#2\\ \hspace{-5pt}#3\end{array}\hspace{-5pt}\right\}}
\def\cd{\cdots}

\def\cc{{\bf c}}
\def\dd{{\bf d}}

\def\del{\delta}
\def\Del{\Delta}

\def\Delre{\Delta^{\rm re}}

\def\eit{\tilde{e}_i}
\def\eneme{\end{enumerate}}
\def\ep{\epsilon}
\def\eeq{\end{equation}}
\def\eeqn{\end{eqnarray}}
\def\eeqnn{\end{eqnarray*}}
\def\fit{\tilde{f}_i}

\def\gau#1,#2{\left[\begin{array}{c}\hspace{-5pt}#1\\
\hspace{-5pt}#2\end{array}\hspace{-5pt}\right]}

\def\ify{\infty}

\def\kify{B^{k,\infty}}

\def\lan{\langle}

\def\lm{\lambda}
\def\Lm{\Lambda}
\def\mapright#1{\smash{\mathop{\longrightarrow}\limits^{#1}}}

\def\nd{\noindent}
\def\nn{\nonumber}

\def\ot{\otimes}

\def\osigma{\ovl\sigma}
\def\ovl{\overline}

\def\qq{\qquad}
\def\q{\quad}
\def\qed{\hfill\framebox[2mm]{}}

\def\ra{\rightarrow}
\def\ran{\rangle}

\def\syl{\scriptstyle}

\def\til{\tilde}
\def\tm{\times}
\def\tt{\frt}
\def\TY(#1,#2,#3){{#1^{(#2)}_{#3}}}
\def\TYO(#1,#2,#3){\ovl{#1}^{(#2)}_{#3}}
\def\TYOS(#1,#2,#3){\ovl{#1}^{*(#2)}_{#3}}
\def\aTY(#1,#2,#3){{{}^1 #1^{(#2)}_{#3}}}
\def\bTY(#1,#2,#3){{{}^2 #1^{(#2)}_{#3}}}
\def\jTY(#1,#2,#3){{{}^j #1^{(#2)}_{#3}}}
\def\aTYO(#1,#2,#3){{{}^1 \ovl{#1}^{(#2)}_{#3}}}
\def\bTYO(#1,#2,#3){{{}^2 \ovl{#1}^{(#2)}_{#3}}}
\def\TYS(#1,#2,#3){#1^{*(#2)}_{#3}}
\def\UD{{\mathcal UD}}

\def\uqp{U'_q(\ge)}

\def\vep{\varepsilon}
\def\vp{\varphi}

\def\W1{W(\varpi_1)}

\def\wt{{\rm wt}}
\def\wtil{\widetilde}
\def\what{\widehat}

\def\ZZ{\mathbb Z}

\def\m@th{\mathsurround=0pt}
\def\fsquare(#1,#2){
\hbox{\vrule$\hskip-0.4pt\vcenter to #1{\normalbaselines\m@th
\hrule\vfil\hbox to #1{\hfill$\scriptstyle #2$\hfill}\vfil\hrule}$\hskip-0.4pt
\vrule}}

\theoremstyle{definition}
\newtheorem{df}{Definition}[section]
\newtheorem{thm}[df]{Theorem}
\newtheorem{pro}[df]{Proposition}
\newtheorem{lem}[df]{Lemma}

\newtheorem{cor}[df]{Corollary}

\newcommand{\seteq}{\mathbin{:=}}
\newcommand{\cl}{\colon}
\newcommand{\be}{\begin{enumerate}}
\newcommand{\ee}{\end{enumerate}}

\newcommand{\ba}{\begin{array}}
\newcommand{\ea}{\end{array}}

\newcommand{\eq}{\begin{eqnarray}}
\newcommand{\eneq}{\end{eqnarray}}

\title{Affine Geometric Crystal of $\TY(A,1,n)$
and Limit of Kirillov-Reshetikhin Perfect Crystals}
\author{\begin{tabular}{c}
Kailash C. M\textsc{isra}
\thanks{supported in part by Simons Foundation Grant $\sharp  307555$.
E-mail:misra@ncsu.edu
}
\end{tabular}
\begin{tabular}{c}
Toshiki N\textsc{akashima}\thanks{supported in part by 
JSPS Grants in Aid for 
Scientific Research $\sharp 15K04794$.
E-mail:toshiki@ sophia.ac.jp}
{}\quad
\end{tabular}}

\date{}

\begin{document}
\maketitle
\begin{abstract}
Let $\ge$ be an affine Lie algebra with index set 
$I = \{0, 1, 2, \cdots , n\}$ and $\ge^L$ be its Langlands dual. 
It is conjectured in \cite{KNO} that  for each 
$k \in I \setminus \{0\}$ the affine Lie algebra $\ge$ has 
a positive geometric crystal whose ultra-discretization is 
isomorphic to the limit of certain coherent family of perfect crystals 
for $\ge^L$. Motivated by this conjecture we construct a positive geometric crystal for the affine Lie algebra $\ge =  \TY(A,1,n)$ for each Dynkin index $k\in I\setminus\{0\}$ and show that its ultra-discretization is isomorphic to the limit of a coherent family of perfect crystals for $ \TY(A,1,n)$ given in \cite{OSS}. In the process we develop and use some lattice-path combinatorics. 

\end{abstract}

\renewcommand{\thesection}{\arabic{section}}
\section{Introduction}
\setcounter{equation}{0}
\renewcommand{\theequation}{\thesection.\arabic{equation}}
Let $\ge = \ge(A)$ be an affine Lie algebra with Cartan matrix  $A= (a_{ij})_{i,j \in I}, I = \{0, 1, \cdots , n\}$ and Cartan datum $(A, \{\al_i\}_{i \in I},\al^\vee_i\}_{\i\in I})$.  Let $U_q(\ge)$ denote the quantum affine algebra associated with $\ge$. We denote $\ge_i$ to be the subalgebra of $\ge$ with index set 
$I \setminus \{i\}$. Note that both $\ge_0$ and $\ge_n$ are isomorphic to the simple Lie algebra 
$A_n$.
We denote the affine weight lattice 
and the dual affine weight lattice of $\ge$ by $P= \Z \Lambda_0 \oplus \Z \Lambda_1\oplus 
\cdots \oplus \Z \Lambda_n \oplus \Z\delta$ and  
$P^\vee = \Z \al^\vee_0 \oplus \Z \al^\vee_1 \oplus \cdots \oplus 
\Z \al^\vee_n \oplus \Z d $  respectively, where $\delta$ is the null root. 
For a dominant weight  $\lambda \in P^+ = \{\mu \in P \mid \mu (h_i) 
\geq 0 \quad  {\rm for \ \  all} \quad i \in I \}$ of level 
$l = \lambda (c)$ ($c =$ canonical central element of $\ge$), 
Kashiwara defined the crystal base $(L(\lambda), B(\lambda))$
\cite{Kas1} 
for the integrable highest weight $U_q(\ge)$-module $V(\lambda)$. 
The crystal $B(\lambda)$ is the $q= 0$ limit of the canonical basis 
\cite{Lu} or the global crystal basis \cite{Kas2}. There are many known explicit 
realizations of the affine crystal $B(\lambda)$. One such realization is the path
realization \cite{KMN1} using perfect crystals. A perfect crystal is a crystal for 
certain finite dimensional module called Kirillov-Reshetikhin module (KR-module for short) 
of the quantum affine algebra $U_q(\ge)$ (\cite{KR}, \cite{HKOTY,HKOTT}). 
The KR-modules are parametrized by two integers $(k, l)$, where 
$k \in I \setminus \{0\}$ and $l$ any positive integer. 
Let $\{\varpi_k\}_{k\in I\setminus\{0\}}$ be the set of level $0$ 
fundamental weights \cite{K0} . Hatayama et al  
(\cite{HKOTY,HKOTT}) conjectured that any KR-module
$W(l\varpi_k)$ admit a crystal base $B^{k,l}$ in the sense of Kashiwara 
and furthermore $B^{k,l}$ is perfect if $l$ is a multiple of 
$c_k^\vee\seteq \mathrm{max }(1,\frac{2}{(\al_k,\al_k)})$. 
This conjecture has been proved for quantum affine algebras 
$U_q(\ge)$ of classical types (\cite{OS,FOS1,FOS2}). 
When $\{B^{k,l}\}_{l\geq 1}$ is a coherent family of perfect crystals \cite{KKM}
we denote its limit by $B^\infty (\varpi_k)$ (or just $B^{k,\infty}$).

In 2000, Berenstein and Kazhdan \cite{BK} introduced the notion of geometric
crystal for reductive algebraic groups which was extended to Kac-Moody groups
in \cite{N}.
For a given Cartan datum $(A, \{\alpha_i\}_{i \in I},
\{\al^\vee_i\}_{\i\in I} )$, 
the geometric crystal is defined as a quadruple 
$\cV(\ge)=(X, \{e_i\}_{i \in I}, \{\gamma_i\}_{i \in I},$ 
$\{\vep_i\}_{i\in I})$, 
where $X$ is an algebraic variety,  $e_i:\bbC^\times\times
X\longrightarrow X$ 
are rational $\bbC^\times$-actions and  
$\gamma_i,\vep_i:X\longrightarrow 
\bbC$ $(i\in I)$ are rational functions satisfying certain conditions  
( see Definition \ref{def-gc}). 
A geometric 
crystal is said to be a positive geometric crystal 
if it admits a positive structure (see Definition 2.5).
A remarkable relation between positive geometric crystals 
and algebraic crystals is the ultra-discretization functor $\mathcal
{UD}$ between them (see Section 2.4). Applying this functor, positive rational 
functions are transfered to piecewise-linear 
functions by the simple correspondence:
$$
x \times y \longmapsto x+y, \qquad \frac{x}{y} \longmapsto x - y, 
\qquad x + y \longmapsto {\rm max}\{x, y\}.
$$

It was conjectured in \cite{KNO} that for each affine Lie algebra $\ge$ and 
each Dynkin index $i \in I \setminus {0}$, there exists a positive geometric crystal
$\cV(\ge)=(X, \{e_i\}_{i \in I}, \{\gamma_i\}_{i \in I}, 
\{\vep_i\}_{i\in I})$ whose ultra-discretization $\mathcal{UD}(\cV)$ is isomorphic 
to the limit $B^{\infty}$ of a coherent family of perfect crystals for the Langlands dual $\ge^L$.
This conjecture has been shown to be true for $k=1$ and $\ge = A_n^{(1)}, 
B_n^{(1)}, C_n^{(1)}, D_n^{(1)}, A_{2n-1}^{(2)}, A_{2n}^{(2)},$
$D_{n+1}^{(2)}$ \cite{KNO}, $D_4^{(3)}$ \cite{IMN}, $G_2^{(1)}$ \cite{N4}.
For $\ge = A_n^{(1)}, k=2$ we showed the conjecture to be true in \cite{MN}.

In this paper we have constructed positive geometric crystals associated
with each Dynkin index $k \in I \setminus {0}$ for the affine Lie algebra $A_n^{(1)}$ and have
proved that its ultra-discretization is isomorphic to the limit $B^{k,\infty}$
of the coherent family of perfect crystals $\{B^{k,l}\}_{l \geq 1}$ given in (\cite{KMN2,OSS}). To construct the $A_n^{(1)}$ positive geometric crystal we proceed as follows. We construct $\ge_0$ and $\ge_n$ positive geometric crystals $\cV_1
=\{\cV_1(x),\{e_i\},\{\gamma_i\},\{\vep_i\} \mid 1 \leq i \leq n\}$ and $\cV_2
=\{\cV_2(y),\{\ovl e_i\},\{\ovl \gamma_i\},\{\ovl\vep_i\}\mid 0 \leq i \leq n-1\}$ respectively, in the fundamental representation $W(\varpi_k)$. Then we define a birational
isomorphism $\osigma$ between $\cV_1$ and $\cV_2$, and using this isomorphism
we patch them together to obtain an affine geometric crystal 
$\cV(\TY(A,1,n))$. In this case we conjecture that $y = \osigma(x)$ is the unique solution of $\cV_2(y) = a(x) \cV_1(x)$ which would prove Conjecture 1.2 in \cite{KNO} completely.

This paper is organized as follows. In Section 2, 
we recall necessary definitions and facts about geometric crystals and ultra-discretization. In Section 3 we give $\ge_0$ (resp. $\ge_n$) positive geometric crystal $\cV_1$ (resp. $\cV_2$) explicitly. In Section 4 we develop some lattice-path combinatorics which is used in the sequel.  In Section 5 we give the birational bi-positive map $\osigma$ between $\cV_1$ and $\cV_2$. In Section 6 we construct the positive affine geometric crystal $\cV(\TY(A,1,n))$. In Section 7, we show that the ultra-discretization of the affine geometric crystal $\cV(\TY(A,1,n))$ is isomorphic to the crystal $B^{k,\infty}$ given in \cite{OSS}. Finally in Section 8 as an application we give the explicit actions of the affine Weyl group on the geometric crystal $\cV(\TY(A,1,n))$ and its ultra-discretization $B^{k,\infty}$.

\newpage
\renewcommand{\thesection}{\arabic{section}}
\section{Geometric crystals}
\setcounter{equation}{0}
\renewcommand{\theequation}{\thesection.\arabic{equation}}

We review Kac-Moody groups and geometric crystals
following  \cite{BK,Ku2,N,PK}.
\subsection{Kac-Moody algebras and Kac-Moody groups}
\label{KM}
Let  $A=(a_{ij})_{i,j\in I}$ be a symmetrizable generalized Cartan matrix
 with a finite index set $I$.
and  $(\tt,\{\al_i\}_{i\in I}, \{\al^\vee_i\}_{i\in I})$ 
the associated
root data, where ${\tt}$ is a vector space 
over $\bbC$ and
$\{\al_i\}_{i\in I}\subset\tt^*$ and 
$\{\al^\vee_i\}_{i\in I}\subset\tt$
are linearly independent 
satisfying $\al_j(\al^\vee_i)=a_{ij}$.

The Kac-Moody Lie algebra $\ge=\ge(A)$ associated with $A$
is the Lie algebra over $\bbC$ generated by $\tt$, the 
Chevalley generators $e_i$ and $f_i$ $(i\in I)$
with the usual defining relations (\cite{KP,PK}).
There is the root space decomposition 
$\ge=\bigoplus_{\al\in \tt^*}\ge_{\al}$.
Denote the set of roots by 
$\Delta:=\{\al\in \tt^*|\al\ne0,\,\,\ge_{\al}\ne(0)\}$.
Set $Q=\sum_i\bbZ \al_i$, $Q_+=\sum_i\bbZ_{\geq0} \al_i$,
$Q^\vee:=\sum_i\bbZ \al^\vee_i$
and $\Delta_+:=\Delta\cap Q_+$.
An element of $\Delta_+$ is called 
a {\it positive root}.
Let $P\subset \tt^*$ be a weight lattice such that 
$\bbC\ot P=\tt^*$, whose element is called a
weight.

Define simple reflections $s_i\in{\rm Aut}(\tt)$ $(i\in I)$ by
$s_i(h):=h-\al_i(h)\al^\vee_i$, which generate the Weyl group $W$.
It induces the action of $W$ on $\tt^*$ by
$s_i(\lm):=\lm-\lm(\al^\vee_i)\al_i$.
Set $\Delre:=\{w(\al_i)|w\in W,\,\,i\in I\}$, whose element 
is called a real root.

Let $\ge'$ be the derived Lie algebra 
of $\ge$ and let 
$G$ be the Kac-Moody group associated 
with $\ge'$(\cite{PK}).
Let $U_{\al}:=\exp\ge_{\al}$ $(\al\in \Delre)$
be the one-parameter subgroup of $G$.
The group $G$ is generated by $U_{\al}$ $(\al\in \Delre)$.
Let $U^{\pm}$ be the subgroup generated by $U_{\pm\al}$
($\al\in \Delre_+=\Delre\cap Q_+$), {\it i.e.,}
$U^{\pm}:=\lan U_{\pm\al}|\al\in\Del^{\rm re}_+\ran$.

For any $i\in I$, there exists a unique homomorphism;
$\phi_i:SL_2(\bbC)\rightarrow G$ such that
\[
\hspace{-2pt}\phi_i\left(
\left(
\begin{array}{cc}
c&0\\
0&c^{-1}
\end{array}
\right)\right)=c^{\al^\vee_i},\,
\phi_i\left(
\left(
\begin{array}{cc}
1&t\\
0&1
\end{array}
\right)\right)=\exp(t e_i),\,
 \phi_i\left(
\left(
\begin{array}{cc}
1&0\\
t&1
\end{array}
\right)\right)=\exp(t f_i).
\]
where $c\in\bbC^\times$ and $t\in\bbC$.
Set $\al^\vee_i(c):=c^{\al^\vee_i}$,
$x_i(t):=\exp{(t e_i)}$, $y_i(t):=\exp{(t f_i)}$, 
$G_i:=\phi_i(SL_2(\bbC))$,
$T_i:=\phi_i(\{{\rm diag}(c,c^{-1})\vert 
c\in\bbC^{\vee}\})$ 
and 
$N_i:=N_{G_i}(T_i)$. Let
$T$ (resp. $N$) be the subgroup of $G$ 
with the Lie algebra $\tt$
(resp. generated by the $N_i$'s), 
which is called a {\it maximal torus} in $G$, and let
$B^{\pm}=U^{\pm}T$ be the Borel subgroup of $G$.
We have the isomorphism
$\phi:W\mapright{\sim}N/T$ defined by $\phi(s_i)=N_iT/T$.
An element $\ovl s_i:=x_i(-1)y_i(1)x_i(-1)
=\phi_i\left(
\left(
\begin{array}{cc}
0&\pm1\\
\mp1&0
\end{array}
\right)\right)$ is in 
$N_G(T)$, which is a representative of 
$s_i\in W=N_G(T)/T$. 

\subsection{Geometric crystals}
Let $X$ be an ind-variety , 
{$\gamma_i:X\rightarrow \bbC$} and 
$\vep_i:X\longrightarrow \bbC$ ($i\in I$) 
rational functions on $X$, and
{$e_i:\bbC^\times \times X\longrightarrow X$}
$((c,x)\mapsto e^c_i(x))$ a
rational $\bbC^\times$-action.

\begin{df}
\label{def-gc}
A quadruple $(X,\{e_i\}_{i\in I},\{\gamma_i,\}_{i\in I},
\{\vep_i\}_{i\in I})$ is a 
$G$ (or $\ge$)-{\it geometric} {\it crystal} 
if
\begin{enumerate}
\item
$\{1\}\times X\subset dom(e_i)$ 
for any $i\in I$.
\item
$\gamma_j(e^c_i(x))=c^{a_{ij}}\gamma_j(x)$.
\item $\{e_i\}_{i\in I}$ satisfy the following relations.
\[
 \begin{array}{lll}
&\hspace{-20pt}e^{c_1}_{i}e^{c_2}_{j}
=e^{c_2}_{j}e^{c_1}_{i}&
{\rm if }\,\,a_{ij}=a_{ji}=0,\\
&\hspace{-20pt} e^{c_1}_{i}e^{c_1c_2}_{j}e^{c_2}_{i}
=e^{c_2}_{j}e^{c_1c_2}_{i}e^{c_1}_{j}&
{\rm if }\,\,a_{ij}=a_{ji}=-1,\\
&\hspace{-20pt}
e^{c_1}_{i}e^{c^2_1c_2}_{j}e^{c_1c_2}_{i}e^{c_2}_{j}
=e^{c_2}_{j}e^{c_1c_2}_{i}e^{c^2_1c_2}_{j}e^{c_1}_{i}&
{\rm if }\,\,a_{ij}=-2,\,
a_{ji}=-1,\\
&\hspace{-20pt}
e^{c_1}_{i}e^{c^3_1c_2}_{j}e^{c^2_1c_2}_{i}
e^{c^3_1c^2_2}_{j}e^{c_1c_2}_{i}e^{c_2}_{j}
=e^{c_2}_{j}e^{c_1c_2}_{i}e^{c^3_1c^2_2}_{j}e^{c^2_1c_2}_{i}
e^{c^3_1c_2}_je^{c_1}_i&
{\rm if }\,\,a_{ij}=-3,\,
a_{ji}=-1,
\end{array}
\]
\item
$\vep_i(e_i^c(x))=c^{-1}\vep_i(x)$ and $\vep_i(e_j^c(x))=\vep_i(x)$ if 
$a_{i,j}=a_{j,i}=0$.
\end{enumerate}
\end{df}

The condition (iv) is slightly modified from the one in 
\cite{IN,N3,N4}.

Let $W$ be the  Weyl group associated with $\ge$. 
For $w \in W$ define $R(w)$  by
\[
 R(w):=\{(i_1,i_2,\cd,i_l)\in I^l|w=s_{i_1}s_{i_2}\cd s_{i_l}\},
\]
where $l$ is the length of $w$.
Then $R(w)$ is the set of reduced words of $w$.
For a word ${\bf i}=(i_1,\cd,i_l)\in R(w)$ 
$(w\in W)$, set 
$\al^{(j)}:=s_{i_l}\cd s_{i_{j+1}}(\al_{i_j})$ 
$(1\leq j\leq l)$ and 
\begin{eqnarray*}
e_{\bf i}:&T\times X\rightarrow &X\\
&(t,x)\mapsto &e_{\bf i}^t(x):=e_{i_1}^{\al^{(1)}(t)}
e_{i_2}^{\al^{(2)}(t)}\cd e_{i_l}^{\al^{(l)}(t)}(x).
\end{eqnarray*}
Note that the condition (iii) above is 
equivalent to the following:
{$e_{\bf i}=e_{\bf i'}$}
for any 
$w\in W$, ${\bf i}$.
${\bf i'}\in R(w)$.

\subsection{Geometric crystal on Schubert cell}
\label{schubert}

Let $w\in W$ be a Weyl group element and take a 
reduced expression $w=s_{i_1}\cd s_{i_l}$. 
Let $X:=G/B$ be the flag
variety, which is an ind-variety 
and $X_w\subset X$ the
Schubert cell associated with $w$, which satisfies $X=\sqcup_{w\in W}X_w$ and 
has a natural geometric crystal structure
(\cite{BK,N}).
For a reduced word ${\bf i}:=(i_1,\cd,i_k)$ of $w\in W$, set 
\begin{equation}
B_{\bf i}^-
:=\{Y_{\bf i}(c_1,\cd,c_k)
:=Y_{i_1}(c_1)\cd Y_{i_l}(c_k)
\,\vert\, c_1\cd,c_k\in\bbC^\times\}\subset B^-,
\label{bw1}
\end{equation}
where $Y_i(c):=y_i(\frac{1}{c})\al^\vee_i(c)$.
If $I=\{i_1,\cd,i_k\}$, this has a geometric crystal structure(\cite{N})
isomorphic to $X_w$. 
The explicit forms of the action $e^c_i$, the rational 
function $\vep_i$  and $\gamma_i$ on 
$B_{\bf i}^-$ are given by
\begin{eqnarray}
&& e_i^c(Y_{\bf i}(c_1,\cd,c_k))
=Y_{\bf i}({\mathcal C}_1,\cd,{\mathcal C}_k)),\nn \\
&&\text{where}\nn\\
&&{\mathcal C}_j:=
c_j\cdot \frac{\displaystyle \sum_{1\leq m\leq j,i_m=i}
 \frac{c}
{c_1^{a_{i_1,i}}\cd c_{m-1}^{a_{i_{m-1},i}}c_m}
+\sum_{j< m\leq k,i_m=i} \frac{1}
{c_1^{a_{i_1,i}}\cd c_{m-1}^{a_{i_{m-1},i}}c_m}}
{\displaystyle\sum_{1\leq m<j,i_m=i} 
 \frac{c}
{c_1^{a_{i_1,i}}\cd c_{m-1}^{a_{i_{m-1},i}}c_m}+
\mathop\sum_{j\leq m\leq k,i_m=i}  \frac{1}
{c_1^{a_{i_1,i}}\cd c_{m-1}^{a_{i_{m-1},i}}c_m}},
\label{eici}\\
&& \vep_i(Y_{\bf i}(c_1,\cd,c_k))=
\sum_{1\leq m\leq k,i_m=i} \frac{1}
{c_1^{a_{i_1,i}}\cd c_{m-1}^{a_{i_{m-1},i}}c_m},
\label{vep-i}\\
&&\gamma_i(Y_{\bf i}(c_1,\cd,c_k))
=c_1^{a_{i_1,i}}\cd c_k^{a_{i_k,i}}.
\label{gamma-i}
\end{eqnarray}
{\sl Remark.}
As in \cite{N}, the above setting requires the condition
$I=\{i_1,\cd,i_k\}$.
Otherwise, set $J:=\{i_1,\cd,i_k\}\subsetneq I$ and let $\ge_J\subsetneq \ge$ 
be the corresponding subalgebra.
Then, by arguing similarly to \cite[4.3]{N}, we can define the $\ge_J$-geometric crystal 
structure on $B^-_{\bf i}$.
\subsection{Positive structure and 
Ultra-discretizations}
\label{positive-str}

Let us recall the notions of 
positive structure and  ultra-discretization following \cite{BK,NKNO,}.

Let $T=(\bbC^\times)^l$ be an algebraic torus over $\bbC$ and 
$X^*(T):={\rm Hom}(T,\bbC^\times)\cong \ZZ^l$ 
(resp. $X_*(T):={\rm Hom}(\bbC^\times,T)\cong \ZZ^l$) 
be the lattice of characters
(resp. co-characters)
of $T$.  Define
$$
\begin{array}{cccc}
v:&\bbC(c)\setminus\{0\}&\longrightarrow &\ZZ\\
&f(c)&\mapsto
&{\rm deg}(f(c)),
\end{array}
$$
where $\rm deg$ is the degree of poles at $c=\ify$. 
Here note that for $f_1,f_2\in \bbC(c)\setminus\{0\}$, we have
\begin{equation}
v(f_1 f_2)=v(f_1)+v(f_2),\q
v\left(\frac{f_1}{f_2}\right)=v(f_1)-v(f_2)
\label{ff=f+f}
\end{equation}
A non-zero rational function on
an algebraic torus $T$ is called {\em positive} if
it can be written as $g/h$ where
$g$ and $h$ are positive linear combinations of
characters of $T$.
\begin{df}
Let 
$f\cl T\rightarrow T'$ be 
a rational morphism between
two algebraic tori $T$ and 
$T'$.
We say that $f$ is {\em positive},
if $\eta\circ f$ is positive
for any character $\eta\cl T'\to \C$.
\end{df}
Denote by ${\rm Mor}^+(T,T')$ the set of 
positive rational morphisms from $T$ to $T'$.

\begin{lem}[\cite{BK}]
\label{TTT}
For any $f\in {\rm Mor}^+(T_1,T_2)$             
and $g\in {\rm Mor}^+(T_2,T_3)$, 
the composition $g\circ f$
is well-defined and belongs to ${\rm Mor}^+(T_1,T_3)$.
\end{lem}

By this lemma, we can define a category ${\mathcal T}_+$
whose objects are algebraic tori over $\bbC$ and arrows
are positive rational morphisms.

Let $f\cl T\rightarrow T'$ be a 
positive rational morphism
of algebraic tori $T$ and 
$T'$.
We define a map $\what f\cl X_*(T)\rightarrow X_*(T')$ by 
\[
\langle\eta,\what f(\xi)\rangle
=v(\eta\circ f\circ \xi),
\]
where $\eta\in X^*(T')$ and $\xi\in X_*(T)$.
\begin{lem}[\cite{BK}]
For any algebraic tori $T_1$, $T_2$, $T_3$, 
and positive rational morphisms 
$f\in {\rm Mor}^+(T_1,T_2)$, 
$g\in {\rm Mor}^+(T_2,T_3)$, we have
$\what{g\circ f}=\what g\circ\what f.$
\end{lem}
Let ${\hbox{\germ Set}}$ denote the category of sets with the morphisms being set maps.
By the above lemma, we obtain a functor: 
\[
\begin{array}{cccc}
{\mathcal{UD}}:&{\mathcal T}_+&\longrightarrow &{{\hbox{\germ Set}}}\\
&T&\mapsto& X_*(T)\\
&(f:T\rightarrow T')&\mapsto& 
(\what f:X_*(T)\rightarrow X_*(T')))
\end{array}
\]

\begin{df}[\cite{BK}]
Let $\chi=(X,\{e_i\}_{i\in I},\{{\rm wt}_i\}_{i\in I},
\{\vep_i\}_{i\in I})$ be a 
finite dimensional geometric crystal, $T'$ an algebraic torus
and $\theta:T'\rightarrow X$ 
a birational isomorphism.
The isomorphism $\theta$ is called 
{\it positive structure} on
$\chi$ if it satisfies
\begin{enumerate}
\item for any $i\in I$ the rational functions
$\gamma_i\circ \theta:T'\rightarrow \bbC$ and 
$\vep_i\circ \theta:T'\rightarrow \bbC$ 
are positive.
\item
For any $i\in I$, the rational morphism 
$e_{i,\theta}:\bbC^\tm \tm T'\rightarrow T'$ defined by
$e_{i,\theta}(c,t)
:=\theta^{-1}\circ e_i^c\circ \theta(t)$
is positive.
\end{enumerate}
\end{df}
Let $\theta:T\rightarrow X$ be a positive structure on 
a geometric crystal $\chi=(X,\{e_i\}_{i\in I},$
$\{{\rm wt}_i\}_{i\in I},
\{\vep_i\}_{i\in I})$.
Applying the functor ${\mathcal{UD}}$ 
to positive rational morphisms
$e_{i,\theta}:\bbC^\tm \tm T\rightarrow T$ and
$\gamma_i\circ \theta,\vep_i\circ\theta:T\ra \bbC$
(the notations are
as above), we obtain
\begin{eqnarray*}
\til e_i&:=&{\mathcal{UD}}(e_{i,\theta}):
\ZZ\tm X_*(T) \rightarrow X_*(T)\\
{\rm wt}_i&:=&{\mathcal{UD}}(\gamma_i\circ\theta):
X_*(T)\rightarrow \bbZ,\\
\vep_i&:=&{\mathcal{UD}}(\vep_i\circ\theta):
X_*(T)\rightarrow \bbZ.
\end{eqnarray*}
Now, for given positive structure $\theta:T'\rightarrow X$
on a geometric crystal 
$\chi=(X,\{e_i\}_{i\in I},$
$\{{\rm wt}_i\}_{i\in I},
\{\vep_i\}_{i\in I})$, we associate 
the quadruple $(X_*(T'),\{\til e_i\}_{i\in I},
\{{\rm wt}_i\}_{i\in I},\{\vep_i\}_{i\in I})$
with a free pre-crystal structure (see \cite[Sect.7]{BK}) 
and denote it by ${\mathcal{UD}}_{\theta,T'}(\chi)$.
We have the following theorem:

\begin{thm}[\cite{BK,N}]
For any geometric crystal 
$\chi=(X,\{e_i\}_{i\in I},\{\gamma_i\}_{i\in I},$
$\{\vep_i\}_{i\in I})$ and positive structure
$\theta:T'\rightarrow X$, the associated pre-crystal
${\mathcal{UD}}_{\theta,T'}(\chi)=
(X_*(T'),\{\til e_i\}_{i\in I},\{{\rm wt}_i\}_{i\in I},
\{\vep_i\}_{i\in I})$ 
is a crystal {\rm (see \cite[Sect.7]{BK})}
\end{thm}

Now, let ${\mathcal{GC}}^+$ be a category whose 
object is a triplet
$(\chi,T',\theta)$ where 
$\chi=(X,\{e_i\},\{\gamma_i\},\{\vep_i\})$ 
is a geometric crystal and $\theta:T'\rightarrow X$ 
is a positive structure on $\chi$, and morphism
$f:(\chi_1,T'_1,\theta_1)\longrightarrow 
(\chi_2,T'_2,\theta_2)$ is given by a rational map
$\vp:X_1\longrightarrow X_2$  
($\chi_i=(X_i,\cd)$) such that 
\begin{eqnarray*}
&&\vp\circ e^{X_1}_{i}=e^{X_2}_{i}\circ\vp,\q \gamma^{X_2}_{i}\circ\vp=\gamma^{X_1}_{i},\q
\vep^{X_2}_{i}\circ\vp=\vep^{X_1}_{i},\\
&&\text{ and }f:=\theta_2^{-1}\circ\vp\circ\theta_1:T'_1\longrightarrow T'_2,
\end{eqnarray*}
is a positive rational morphism. Let ${\mathcal{CR}}$
be the category of crystals. 
Then by the theorem above, we have
\begin{cor}
\label{cor-posi}
The map $ \mathcal{UD} = \mathcal{UD}_{\theta,T'}$ defined above is a functor
\begin{eqnarray*}
 {\mathcal{UD}}&:&{\mathcal{GC}}^+\longrightarrow {\mathcal{CR}},\\
&&(\chi,T',\theta)\mapsto X_*(T'),\\
&&(f:(\chi_1,T'_1,\theta_1)\rightarrow 
(\chi_2,T'_2,\theta_2))\mapsto
(\what f:X_*(T'_1)\rightarrow X_*(T'_2)).
\end{eqnarray*}

\end{cor}
We call the functor $\mathcal{UD}$
{\it ``ultra-discretization''} as in (\cite{N,N2})
instead of ``tropicalization'' as in \cite{BK}.
And 
for a crystal $B$, if there
exists a geometric crystal $\chi$ and a positive 
structure $\theta:T'\rightarrow X$ on $\chi$ such that 
${\mathcal{UD}}(\chi,T',\theta)\cong B$ as crystals, 
we call an object $(\chi,T',\theta)$ in ${\mathcal{GC}}^+$
a {\it tropicalization} of $B$, which is not standard but
we use such a terminology as before.


\renewcommand{\thesection}{\arabic{section}}
\section{$A_n$-
Geometric Crystals $\cV_1$ and $\cV_2$ in $W(\varpi_k)$ }
\setcounter{equation}{0}
\renewcommand{\theequation}{\thesection.\arabic{equation}}
\subsection{Affine Lie Algebra $\TY(A,1,n)$}

In the sequel, we assume $\ge$ to be the affine Lie algebra 
$\TY(A,1,n)$ $(n \geq 2)$. 
The Cartan matrix $A=(a_{ij})_{i,j \in I}$, $I = \{0, 1, \cdots , n\}$
is given by:
\begin{eqnarray*}
a_{ij}= \begin{cases}
2 \qquad &{\rm if} \quad i=j ,\\
-1 \qquad &{\rm if} \quad  i \equiv (j \pm 1) \ \ {\rm mod} (n+1),\\
0 \qquad &{\rm otherwise,}
\end{cases}
\end{eqnarray*}
and its Dynkin diagram is as follows. 

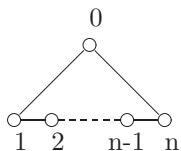
\begin{figure}[h!]
\begin{center}

\setlength{\unitlength}{.5cm}
\begin{picture}(5,2)
\linethickness{0.05mm}
\put(.175,0){\line(1,0){.65}}
\put(3.175,0){\line(1,0){.65}}
\put(.15,.15){\line(1,1){1.7}}
\put(2.15,1.85){\line(1,-1){1.7}}
\multiput(1.175,0)(.3,0){6}{\line(1,0){.15}}
\put(0,0){\circle{.35}}
\put(0,-.8){1}
\put(1,0){\circle{.35}}
\put(1,-.8){2}
\put(3,0){\circle{.35}}
\put(2.5,-.8){n-1}
\put(4,0){\circle{.35}}
\put(4,-.8){n}
\put(2,2){\circle{.35}}
\put(2,2.5){0}
\end{picture}

\end{center}
\label{Dynkin}
\end{figure}

For the affine Lie algebra $A_n^{(1)}$, let 
$\{\alpha_0, \alpha_1, \cdots \alpha_n\}$, 
$\{\al^\vee_0, \al^\vee_1, \cdots  \al^\vee_n\}$ and 
$\{\Lm_0, \Lm_1, \cdots \Lm_n\}$ be the set of 
simple roots, simple coroots and fundamental weights, respectively.

The standard null root $\delta$ 
and the canonical central element $c$ are 
given by
\[
\delta=\alpha_0+\alpha_1+ \cdots +\alpha_n
\quad\text{and}\quad c=\al^\vee_0+\al^\vee_1+ \cdots +\al^\vee_n, 
\]
where 
$\al_0=2\Lm_0-\Lm_1-\Lm_n+\del,\q
\al_i=-\Lm_{i-1}+2\Lm_i-\Lm_{i+1}, 1 \le i \le n-1 \q
\al_n=-\Lm_0 - \Lm_{n-1}+2\Lm_n.$

Let $\sigma$ denote the Dynkin diagram 
automorphism such that 
$\sigma (\al_i) = \al_{\overline{i+1}}$, 
where $\overline{i+1} = (i+1) \ {\rm mod} (n+1)$.
Consider the level $0$ fundamental weight $\varpi_k:=\Lm_k-\Lm_0$
$(1\leq  k\leq  n)$.
Let $I_0 = I \setminus \{0\}, \ \ I_n = I \setminus \{n\}$, 
and $\ge _ i$ denote the subalgebra of $\ge$ associated with the index sets 
$I_i (i= 0, n)$. Then $\ge_0$ as well as $\ge_n$ is 
isomorphic to the simple Lie algebra of type $A_n$.

Let $W(\varpi_k)$ be the fundamental representation 
of $\uqp$ associated with $\varpi_k$ (\cite{K0}).
By \cite[Theorem 5.17]{K0}, $W(\varpi_k)$ is a
finite-dimensional irreducible integrable 
$\uqp$-module and has a global basis
with a simple crystal. Thus, we can consider 
the specialization $q=1$ and obtain the 
finite-dimensional $\TY(A,1,n)$-module $W(\varpi_k)$, 
which we call a fundamental representation
of $\TY(A,1,n)$ and use the same notation as above.
We shall present the explicit form of 
$W(\varpi_k)$ below.

\subsection{Fundamental representation 
$W(\varpi_k)$ for $\TY(A,1,n)$}
\label{fundamental}
The $\TY(A,1,n)$-module $W(\varpi_k)$ is an 
${}_{n+1}C_k=\frac{(n+1)!}{k!(n+1-k)!}$-dimensional 
module with the basis:
\[
\{(i_1,i_2,\cd, i_k) \mid 1 \le i_1 <i_2<\cd <i_k \le n+1\}, 
\]
where $(i_1,i_2,\cd, i_k)$ denotes the one-column Young tableaux 
with entries $i_1,i_2,\cd,i_k$.
The explicit actions of $e_i$, $f_i$ and $\al_i^\vee(c)=c^{h_i}$ 
($c\in \bbC^\times$) on these basis vectors
are given as follows:

\nd 
For $1 \le i \le n$, we have
\begin{eqnarray*}
f_i(i_1,\cd,i_j,i_{j+1},\cd, i_k) &=& \begin{cases}
(i_1,\cd,i+1,i_{j+1},\cd, i_k), & \exists j;i_j=i<i_{j+1}-1, \\
0, &\text{otherwise}.
\end{cases} \\
e_i(i_1,\cd,i_{j-1},i_j,\cd, i_k) &=& \begin{cases}
(i_1,\cd,i_{j-1},i,\cd, i_k), & \exists j;i_{j-1}+1<i+1=i_{j+1}, \\
0, &\text{otherwise}.
\end{cases}\\
\al_i^\vee(c)(i_1,\cd, i_k)&=&
\begin{cases}
c(i_1,\cd, i_k),&\exists j;\,i_j=i<i_{j+1}-1,\\
c^{-1}(i_1,\cd, i_k),&\exists j;\,i_j+1<i_{j+1}=i+1,\\
(i_1,\cd, i_k),&{\rm otherwise},
\end{cases}
\end{eqnarray*}
\begin{eqnarray*}
f_0(i_1,\cd, i_k) &=& \begin{cases}
(1,i_1,\cd, i_{k-1}), & i_1\ne1,\,\,i_k=n+1, \\
0, &\text{otherwise}.
\end{cases} \\
e_0(i_1,\cd, i_k) &=& \begin{cases}
(i_2,\cd, i_{k},n+1), & i_1=1,\,\,i_k\ne n+1, \\
0, & \text{otherwise}.
\end{cases}\\
\al_0^\vee(c)(i_1,\cd, i_k)&=&
\begin{cases}
c(i_1,\cd, i_k),&i_1\ne1,\,\,i_k=n+1,\\
c^{-1}(i_1,\cd, i_k),&i_1=1,\,\,i_k\ne n+1,\\
(i_1,\cd, i_k),&{\rm otheriwse}.
\end{cases}
\end{eqnarray*}
Note that in $W(\varpi_k)$ the vector 
$(1,2,\cd,k-1,k)=:u_1$ (resp. $(1,2,\cd,k-1, n+1)=:u_2$) 
is the $\ge_0$ (resp. $\ge_n$) highest weight vector
with weight $\varpi_k = \Lm_k - \Lm_0$ 
(resp. $\sigma^{-1}\varpi_k = \Lm_{k-1} - \Lm_n$).
We also find that $\eit^2=\fit^2=0$ on $W(\varpi_k)$ and then 
$x_i(c)=1+ce_i$ and $y_i(c)=1+cf_i$ on $W(\varpi_k)$. 
\subsection{Geometric Crystals $\cV_1$ and $\cV_2$ in $W(\varpi_k)$ }
\label{v1v2}

To construct the affine geometric crystal $\cV(\TY(A,1,n))$ 
in $W(\varpi_k)$ $(1\leq  k\leq  n)$ explicitly, we shall 
introduce two $A_n$-geometric crystals 
$\cV_1, \cV_2\subset W(\varpi_k)$.

For $\xi\in (\frt^*_{\rm cl})_0$, let $t(\xi)$ be the 
translation as in \cite[Sect 4]{K0} and 
$\wtil\varpi_k:=\max(1,\frac{2}{(\al_k,\al_k)})\varpi_k$ 
as in \cite{K1}. Indeed, 
$\wtil\varpi_k= \varpi_k$ in our case.
Then we have 
\begin{eqnarray*}
&& t(\varpi_k)=
\sigma^k(s_{k'}s_{k'-1}\cd s_1)(s_{k'+1}s_{k'}\cd s_2)\cd(s_ns_{n-1}\cd s_k),
\\
&&t(\sigma^{-1}\varpi_k)=\sigma^k
(s_{k'-1}\cd s_1s_0)(s_{k'}s_{k'-1}\cd s_1)
\cd(s_{n-1}\cd s_{k-1}),
\end{eqnarray*}
where $k'=n+1-k$ and $\sigma$ is the  Dynkin diagram automorphism
as above.
Set $w_1:=\sigma^{-k}t(\varpi_k)\in W$ and 
$w_2:=\sigma^{-k}t(\sigma^{-1}\varpi_k)\in W$.
Associated with these Weyl group elements $w_1, w_2 \in W$,
we shall define algebraic varieties $\cV_1, \ \cV_2\subset W(\varpi_k)$.
Take two variables $x,y$ of $(\bbC^\times)^{kk'}$ as follows:
\[
\hspace{-20pt}
x=\begin{pmatrix}
\TY(x,k,1)&\TY(x,k,2)&\cd&\TY(x,k,k')\\
\cd&\cd&\cd\\
\TY(x,2,k-1)&\TY(x,2,k)&\cd&\TY(x,2,n-1)\\
\TY(x,1,k)&\TY(x,1,k+1)&\cd&\TY(x,1,n)
\end{pmatrix},
\,\,
y=\begin{pmatrix}
\TY(y,k,0)&\TY(y,k,1)&\cd&\TY(y,k,k'-1)\\
\cd&\cd&\cd\\
\TY(y,2,k-2)&\TY(y,2,k-1)&\cd&\TY(y,2,n-2)\\
\TY(y,1,k-1)&\TY(y,1,k)&\cd&\TY(y,1,n-1)
\end{pmatrix}
\in (\bbC^\times)^{kk'},
\]and set 
\begin{eqnarray*}
&&Y_{w_1}(x)=Y_{k'}(\TY(x,k,k'))\cd Y_1(\TY(x,k,1))\cd
Y_{n}(\TY(x,1,n))\cd Y_k(\TY(x,1,k)),\\
&&Y_{w_2}(y)=Y_{k'-1}(\TY(y,k,k'-1))\cd Y_0(\TY(y,k,0))\cd
Y_{n-1}(\TY(y,1,n-1))\cd Y_{k-1}(\TY(y,1,k-1)).
\end{eqnarray*}
Now, we define $\cV_1$ and $\cV_2$ as 
\begin{eqnarray}
&&\cV_1
:=\{v_1(x):=Y_{w_1}(x)u_1
|x\in (\bbC^\times)^{kk'}\}\subset W(\varpi_k),\label{v1}\\
&&\cV_2:=\{v_2(y):=Y_{w_2}(y)u_2
|y\in (\bbC^\times)^{kk'}\}\subset W(\varpi_k).\label{v2}
\end{eqnarray}
Note that the dimensions of $\cV_1$ and $\cV_2$ are 
$kk'$.

Let us see the explicit $\ge_0$(resp. $\ge_n$)
-geometric crystal structure on $\cV_1$ (resp. $\cV_2$)
according to \ref{schubert}.
Indeed, the geometric crystal structure on $\cV_j$ $(j=1,2)$ coincides 
with those on $B^-_{w_j}$ $(j=1,2)$, that is,
\begin{eqnarray*}
&&\gamma_i(Y_{w_j}(x)u_j)=\gamma_i(Y_{w_j}(x)),\q
\vep_i(Y_{w_j}(x)u_j)=\vep_i(Y_{w_j}(x)),\\
&&
e_i^c((Y_{w_j}(x))=Y_{w_j}(x')\Longrightarrow
e_i^c((Y_{w_j}(x)u_j)=Y_{w_j}(x')u_j,\q(j=1,2).
\end{eqnarray*}

Set $\aa:=\max(k-i+1,1)$ and $\bb=\min(k,n-i+1)$. 
Note that $\bb\geq \aa$.

\begin{pro}\label{V1}
The $\ge_0$-geometric crystal structure of $\cV_1
=(\cV_1,\{e_i\},\{\gamma_i\},\{\vep_i\})_{i=1,2,\cd,n}$
is described as follows: for $v_1(x)=Y_{w_1}(\TY(x,i,j))u_1$ and 
$i=1,2,\cd, n$, set
\begin{equation}
\TY(D,l,i)(x):=
\frac{\TY(x,l,i)(\TY(x,l+1,i)\cd \TY(x,\bb-1,i)\TY(x,\bb,i))^2}
{\TY(x,l+1,i-1)\TY(x,l+2,i-1)\cd \TY(x,\bb,i-1)\TY(x,\bb+1,i-1)\cdot
\TY(x,l,i+1)\TY(x,l+1,i+1)\cd \TY(x,\bb-1,i+1)\TY(x,\bb,i+1)}.
\label{Dli}
\end{equation}
\begin{eqnarray}
&&\gamma_i(v_1(x)):=
\TY(D,\aa,i)\frac{\TY(x,\aa,i)}
{\TY(x,\aa,i-1)\TY(x,\aa-1,i+1)},\qq
\vep_i(v_1(x)):=
\sum_{l=\aa}^\bb\frac{1}{\TY(D,l,i)}.
\label{gamma-ep1}
\end{eqnarray}
Set $e_i^c(v_1(x))=v_1(x')$ ($x=(\TY(x,l,m))$, $x'=({\TY(x,l,m)}')$), 
then ${\TY(x,l,m)}'$ is 
given by
\begin{equation}
{\TY(x,l,m)}'=\begin{cases}
\displaystyle \TY(x,l,i)
\frac{
\displaystyle\sum_{p=\aa}^{l-1}\frac{1}{\TY(D,p,i)}+
\displaystyle\sum_{p=l}^\bb\frac{c}{\TY(D,p,i)}}
{\displaystyle
\sum_{p=\aa}^{l}\frac{1}{\TY(D,p,i)}+
\displaystyle\sum_{p=l+1}^\bb\frac{c}{\TY(D,p,i)}},&{\rm if}\,\,m=i,\\
\TY(x,l,m),&{\rm otherwise}.
\end{cases}
\label{eic1}
\end{equation}
\end{pro}
For $\cV_2$, we also obtain the similar description as follows:
Set ${\bf c}:=\max(k-i,1)$ and ${\bf d}
=\min(k,n-i)$. 
Note that ${\bf d}\geq {\bf c}$.

\begin{pro}\label{V2}
The $\ge_n$-geometric crystal structure of $\cV_2
=(\cV_2,\{\ovl e_i\},\{\ovl \gamma_i\},\{\ovl\vep_i\})_{i=0,1,\cd,n-1}$ 
is described as follows: for $v_2(y)=Y_{w_2}(\TY(y,i,j))u_2$ and 
$i=0,1,2,\cd, n-1$, set
\begin{equation}
\TY(E,l,i)(y):=
\frac{\TY(y,l,i)(\TY(y,l+1,i)\cd \TY(y,\dd-1,i)\TY(y,\dd,i))^2}
{\TY(y,l+1,i-1)\TY(y,l+2,i-1)\cd \TY(y,\dd,i-1)\TY(y,\dd+1,i-1)\cdot
\TY(y,l,i+1)\TY(y,l+1,i+1)\cd \TY(y,\dd-1,i+1)\TY(y,\dd,i+1)}.
\label{Eli}
\end{equation}
\begin{eqnarray}
&&\ovl\gamma_i(v_2(y)):=
\TY(E,\cc,i)\frac{\TY(y,\cc,i)}
{\TY(y,\cc,i-1)\TY(y,\cc-1,i+1)},\qq
\ovl\vep_i(v_2(y)):=
\sum_{l=\cc}^{\dd}\frac{1}{\TY(E,l,i)}.
\label{gamma-ep2}
\end{eqnarray}
Set $\ovl e_i^c(v_2(y))=v_2(y')$ ($y=(\TY(y,l,m))$, 
$y'=({\TY(y,l,m)}')$, then ${\TY(y,l,m)}'$ is 
given by
\begin{equation}
{\TY(y,l,m)}'=\begin{cases}
\displaystyle \TY(y,l,i)
\frac{\displaystyle\sum_{p=\cc}^{l-1}\frac{1}{\TY(E,p,i)}+
\displaystyle\sum_{p=l}^{\dd}\frac{c}{\TY(E,p,i)}}
{\displaystyle
\sum_{p=\cc}^{l}\frac{1}{\TY(E,p,i)}+
\displaystyle\sum_{p=l+1}^{\dd}\frac{c}{\TY(E,p,i)}},&{\rm if}\,\,m=i,\\
\TY(y,l,m),&{\rm otherwise}.
\end{cases}
\label{eic2}
\end{equation}
\end{pro}
In the following sections, we shall define a birational
 isomorphism between $\cV_1$ and $\cV_2$, and using the isomorphism
we shall patch them together to obtain an affine geometric crystal 
$\cV(\TY(A,1,n))$.

\renewcommand{\thesection}{\arabic{section}}
\section{Lattice-Path Combinatorics}
\setcounter{equation}{0}
\renewcommand{\theequation}{\thesection.\arabic{equation}}

In this section, we shall treat rectangular lattices with certain weight
on each edge. Here $m\times n$-lattice means a lattice whose vertical 
length is $m$ and horizontal length is $n$. 
Thus, for example, $2\times3$-rectangular lattice is visualized as 
$\threethree(,,,,,)$.

\subsection{Weight of Paths}\label{weight-path}
Now, we fix two positive integers $n,k$ with $k\leq n$ and 
take a $(k-1)\times(n-k)$-rectangular lattice.
We associate a coordinate $(i,j)$ 
$(1\leq i\leq k<i+j\leq n+1)$
to each lattice point 
on the $i$-th horizontal line from the bottom and $(i+j-k)$-th line
from the left and put a weight $\TY(x,i,j)$ ,
which is denoted by 
$L_1[n,k]=L_1[n,k]((\TY(x,i,j)))$ and called $(\TY(x,i,j))$
-weighted lattice.
Similarly, we define the $(\TY(y,i,j))$-weighted lattice 
$L_2[n,k]=L_2[n,k]((\TY(y,i,j)))$ by associating a coordinate $(i,j)$
and a weight $\TY(y,i,j)$
$(1\leq i\leq k\leq i+j\leq n)$ to each lattice point:

\vskip20pt
\setlength{\unitlength}{1mm}
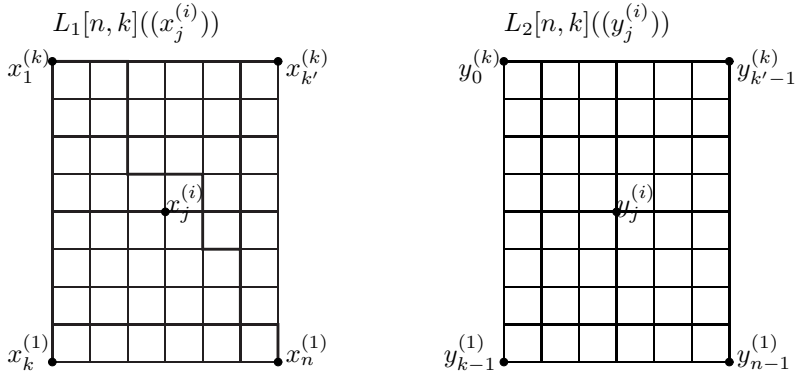
\begin{picture}(67,47)(-11,-5)
\Thicklines
{\path(0,40)(10,40)(10,25)(20,25)(20,15)(25,15)(25,5)(30,5)(30,0)}
\thinlines
\put(0,0){\grid(30,40)(5,5)}
\put(0,0){$\circle*{1}$}
\put(30,0){$\circle*{1}$}\put(0,40){$\circle*{1}$}\put(30,40){$\circle*{1}$}
\put(15,20){$\circle*{1}$}
\put(15,20){$\TY(x,i,j)$}
\put(0,44){$L_1[n,k]((\TY(x,i,j)))$}
\put(-6,0){$\TY(x,1,k)$}
\put(31,0){$\TY(x,1,n)$}
\put(-6,38){$\TY(x,k,1)$}
\put(31,38){$\TY(x,k,k')$}
\put(60,0){\grid(30,40)(5,5)}
\put(60,0){$\circle*{1}$}
\put(90,0){$\circle*{1}$}
\put(60,40){$\circle*{1}$}
\put(90,40){$\circle*{1}$}
\put(75,20){$\circle*{1}$}
\put(75,20){$\TY(y,i,j)$}
\put(60,44){$L_2[n,k]((\TY(y,i,j)))$}
\put(52,0){$\TY(y,1,k-1)$}
\put(91,0){$\TY(y,1,n-1)$}
\put(54,38){$\TY(y,k,0)$}
\put(91,38){$\TY(y,k,k'-1)$}
\end{picture}

Indeed, those lattices  are modelled after coordinates
$x=(\TY(x,i,j)), y=(\TY(y,i,j))\in(\bbC^\times)^{k k'}$ as in Sect.4. 

For the $x$-weighted lattice $L_1[n,k](x)$, suppose that 
each horizontal strip 
$s=(i,j)\frac{\qq}{\qq}(i,j+1)$
has a weight $wt(s)=\frac{\TY(x,i,j)}{\TY(x,i,j+1)}$
and each vertical strip has a weight 1.
Similarly, 
for the $y$-weighted lattice $L_2[n,k](y)$, suppose that 
each horizontal strip has a weight 1 and each vertical strip
$s=(i,j)\frac{\qq}{\qq}(i-1,j+1)$
has a weight $wt(s)=\frac{\TY(y,i-1,j+1)}{\TY(y,i,j)}$.

Set
\[
P_i[n,k]:=\{{\hbox{shortest paths 
on $L_i[n,k]$ from $(k,2-i)$ to $(1,n+1-i)$.}}\}\q(i=1,2).
\]
For a path $p\in P_1[n,k]$  we use the expression
$p=(s_1,s_2\cd,s_k)$ where each $s_i$ is consecutive horizontal strips on the 
$i$-th horizontal line from the top, that is, $s_i$ with 
length $m$ is in the form:
\begin{equation}
s_i=(k-i+1,j)\frac{\q}{\q}(k-i+1,j+1)\frac{\q}{\q}\cd\frac{\q}{\q}
\label{strips}
(k-i+1,j+m).
\end{equation}
Note that we allow that some strips have length 0.
For a path $p\in P_2[n,k]$, we also use the similar expression by 
considering vertical strips.

To each path $p\in P_1[n,k](x)$ (resp.$P_2[n,k](y)$)
 we associate the weight 
\begin{equation}
x(p)({\rm resp.}\,\,y(p)):=\prod_{s:{\rm strip\,\,in\,\,}p}wt(s).
\label{xpyp}
\end{equation}
For $x$(resp.$y$)-weighted lattice 
$L_1[n,k](x)$(resp.$L_2[n,k](y)$), 
define $\vep=\vep(x)$ (resp. $E_2=E_2(y)$)
as a total sum of weights of all shortest paths, that is, 
$\vep=\vep(x):=\sum_{p\in P_1[n,k](x)}x(p)$
(resp. $E_2=E_2(y):=\sum_{p\in P_2[n,k](y)}y(p)$).
For $i=1,2$, define the following set of partial paths in $L_i[n;k]$:
\begin{eqnarray*}
&&P_i[n,k]^{(l)}_m:=\{{\hbox{shortest(partial)
paths from $(l,m)$  to $(1,n+1-i)$}}\} \\
&&{P_i[n,k]^{*(l)}_m}:=\{{\hbox{shortest (partial)
paths from $(k,2-i)$  to $(l,m)$}}\} 
\end{eqnarray*}
We also define the weight $x(p)$ (resp. $y(p)$) 
of a path  $p\in P_i[n,k]^{(l)}_m$, 
$P_i[n,k]^{*(l)}_m$ by the similar way as in \eqref{xpyp} and 
their total summations:
\begin{eqnarray}
&&\TY(X,l,m)=\TY(X,l,m)(x):=\sum_{p\in \TY({P_1[n,k]},l,m)}x(p),\q
\TYS(X,l,m)=\TYS(X,l,m)(x):=\sum_{p\in {P_1[n,k]}^{*(l)}_m}x(p),\\
&&\TY(Y,l,m)=\TY(Y,l,m)(y):=\sum_{p\in \TY({P_2[n,k]},l,m)}y(p),\q
\TYS(Y,l,m)=\TYS(Y,l,m)(y):=\sum_{p\in {P_2[n,k]}^{*(l)}_m}y(p).
\end{eqnarray}
\begin{lem}
We have the formulae:
\begin{eqnarray}
&&\TY(X,l,m)=\TY(X,l-1,m+1)+\frac{\TY(x,l,m)}{\TY(x,l,m+1)}
\TY(X,l,m+1),\label{x=x+x}\\
&&\TY(Y,l,m)=\TY(Y,l,m+1)+\frac{\TY(y,l-1,m+1)}{\TY(y,l,m)}
\TY(Y,l-1,m+1),\\
&&\TYS(X,l,m+1)=\TYS(X,l+1,m)+\frac{\TY(x,l,m)}{\TY(x,l,m+1)}
\TYS(X,l+1,m),\label{*x=x+x}\\
&&\TYS(Y,l,m)=\TYS(Y,l,m-1)
+\frac{\TY(y,l,m)}{\TY(y,l+1,m-1)}\TYS(Y,l+1,m-1).
\label{y=y+y}
\end{eqnarray}
\end{lem}
{\sl Proof.}
We observe that 
the set of
paths $P_1[n,k]^{(l)}_m$ 
are divided into the following two sets:
one is the set of paths 
through $(l-1,m+1)$ and the other is the set of paths 
through $(l,m+1)$.
Thus, \eqref{x=x+x} is immediate consequence of the observation, and
also \eqref{*x=x+x} and \eqref{y=y+y} are obtained similarly.
\qed

Now, we define  rational maps $\Sigma:(\bbC^\times)^{kk'}\to(\bbC^\times)^{kk'}
(x\mapsto y)$
and  $\Xi:(\bbC^\times)^{kk'}\to(\bbC^\times)^{kk'}(y\mapsto x)$ by 
\begin{eqnarray}
&&\TY(y,l,m)=\TY({\Sigma(x)},l,m):=
\TY(x,l+1,m)\frac{\TY(X,l,m)(x)}{\TY(X,l+1,m)(x)},
\q(1\leq l\leq k\leq l+m\leq n),\label{sigxi}\\
&&\TY(x,l,m)=\TY({\Xi(y)},l,m)
:=\TY(y,l,m)\frac{\TYS(Y,l-1,m)(y)}{\TYS(Y,l,m)(y)}
\q(1\leq l\leq k<l+m\leq n+1),\label{xisig}
\end{eqnarray}
where 
$\TY(X,l,m)=\frac{1}{\TY(x,1,n)}$ if $l+m=k$ and 
$\TYS(Y,0,m)=\frac{1}{\TY(y,k,0)}$.
Note that $\TY(y,k,m)=\TY(X,k,m)$ and then $\TY(y,k,0)=\frac{1}{\TY(x,1,n)}$.
\begin{thm}\label{birat}
The rational maps $\Sigma$ and $\Xi$ are both birational and bi-positive.
Furthermore, they are inverse to each other, that is, 
$\Sigma=\Xi^{-1}$ and $\Xi=\Sigma^{-1}$.
\end{thm}
{\sl Proof.}
Bi-positivity of $\Sigma$ and $\Xi$ follows from their
explicit forms. 
Thus, to prove the theorem, it suffices to show
the birationality.
Since 
\begin{eqnarray}
&&(\Xi\circ\Sigma(x))^{(l)}_m=
\TY(x,l+1,m)\cdot\frac{\TY(X,l,m)(x)}{\TY(X,l+1,m)(x)}
\cdot
\frac{\TYS(Y,l-1,m)(\Sigma(x))}{\TYS(Y,l,m)(\Sigma(x))},\\
&&(\Sigma\circ\Xi(y))^{(l)}_m=
\TY(y,l+1,m)\cdot\frac{\TYS(Y,l,m)(y)}{\TYS(Y,l+1,m)(y)}
\cdot
\frac{\TY(X,l,m)(\Xi(y))}{\TY(X,l+1,m)(\Xi(y))},
\end{eqnarray}
we shall show
\begin{eqnarray}
&&\frac{\TY(x,l,m)}{\TY(x,l+1,m)}=\frac{\TY(X,l,m)(x)}{\TY(X,l+1,m)(x)}
\cdot
\frac{\TYS(Y,l-1,m)(\Sigma(x))}{\TYS(Y,l,m)(\Sigma(x))},
\label{xxy}\\
&&\frac{\TY(y,l,m)}{\TY(y,l+1,m)}
=\frac{\TYS(Y,l,m)(y)}{\TYS(Y,l+1,m)(y)}
\cdot
\frac{\TY(X,l,m)(\Xi(y))}{\TY(X,l+1,m)(\Xi(y))},
\label{yyx}
\end{eqnarray}
The formulae \eqref{xxy} and \eqref{yyx} are immediate from the
following lemma. 
\begin{lem}\label{xxyyyx}
We have the formula
\begin{eqnarray}
&&\TY(x,l,m)=\TY(X,l,m)(x)\TYS(Y,l-1,m)(\Sigma(x))\q
(1\leq l\leq k<l+m\leq n+1),
\label{xxy2}\\
&&\TY(y,l,m)=\TYS(Y,l,m)(y)\TY(X,l,m)(\Xi(y))\q
(1\leq l\leq k\leq l+m\leq n)
\label{yyx2}
\end{eqnarray}
\end{lem}
{\sl Proof of Lemma \ref{xxyyyx}.}
Since the proof of \eqref{yyx2} is similar to \eqref{xxy2}, we shall
show \eqref{xxy2} only.
Let $S$ be the set of the coordinates of the lattice points 
in $L_1[n,k](x)$, that is ,
$S:=\{(l,m)|1\leq l\leq k<l+m\leq n+1\}$. 
We define a total order on $S$ by setting:
$(l,m)<(l',m')$ if $l>l'$,  or $l=l'$ and $m<m'$.
The proof of \eqref{xxy2} is proceeded by induction on this order.

First, let us see the case $l=k$, which is shown by induction on 
$m$. The case $m=1$ is clear from that 
\[
\TYS(Y,k-1,1)(\Sigma(x))
=\frac{\TY(y,k-1,1)}{\TY(y,k,0)}(\Sigma(x))
=\TY(x,k,1)\frac{\TY(X,k-1,1)\TY(X,k,0)}{\TY(X,k,1)}
=\frac{\TY(x,k,1)}{\TY(X,k,1)}
\]
Thus, for $m>1$
we obtain from \eqref{x=x+x},
\eqref{y=y+y} and the hypothesis of the induction
that 
\begin{eqnarray*}
\TYS(Y,k-1,m)(\Sigma(x))&=&\TYS(Y,k-1,m-1)(\Sigma(x))
+\frac{\TY(y,k-1,m)}{\TY(y,k,m-1)}\TYS(Y,k,m-1)(\Sigma(x))
=\frac{\TY(x,k,m-1)}{\TY(X,k,m-1)}
+\frac{\TY(x,k,m)\TY(X,k-1,m)}{\TY(X,k,m-1)\TY(X,k,m)}\\
&=&\TY(x,k,m)\frac{\TY(X,k-1,m)+
\frac{\TY(x,k,m-1)}{\TY(x,k,m)}\TY(X,k,m)}{\TY(X,k,m-1)\TY(X,k,m)}
=\frac{\TY(x,k,m)}{\TY(X,k,m)}.
\end{eqnarray*}
Therefore, we have \eqref{xxy2} for $l=k$. 

Let us show $\TYS(Y,l-1,m)(\Sigma(x))=\frac{\TY(x,l,m)}{\TY(X,l,m)}$ for 
$l<k$. 
By the induction hypothesis, we may assume 
\begin{equation}
\TYS(Y,l-1,m-1)(\Sigma(x))=\frac{\TY(x,l,m-1)}{\TY(X,l,m-1)},\qq
\TYS(Y,l,m-1)(\Sigma(x))=\frac{\TY(x,l+1,m-1)}{\TY(X,l+1,m-1)}.
\label{ind-ass}
\end{equation}
Here, by \eqref{x=x+x}, \eqref{y=y+y} and 
the induction hypothesis \eqref{ind-ass} 
we have
\begin{eqnarray*}
&&\TYS(Y,l-1,m)(\Sigma(x))=\TYS(Y,l-1,m-1)(\Sigma(x))
+\frac{\TY(y,l-1,m)}{\TY(y,l,m-1)}
\TYS(Y,l,m-1)(\Sigma(x))=
\frac{\TY(x,l,m-1)}{\TY(X,l,m-1)}
+\frac{\TY(x,l,m)\frac{\TY(X,l-1,m)}{\TY(X,l,m)}}
{\TY(x,l+1,m-1)\frac{\TY(X,l,m-1)}{\TY(X,l+1,m-1)}}
\cdot\frac{\TY(x,l+1,m-1)}{\TY(X,l+1,m-1)}\\
&&=\frac{\TY(x,l,m-1)}{\TY(X,l,m-1)}+\frac{\TY(x,l,m)\TY(X,l-1,m)}
{\TY(X,l,m-1)\TY(X,l,m)}
=\TY(x,l,m)\frac{\TY(X,l-1,m)+\frac{\TY(x,l,m-1)}{\TY(x,l,m)}\TY(X,l,m)}
{\TY(X,l,m-1)\TY(X,l,m)}=\frac{\TY(x,l,m)}{\TY(X,l,m)}
\end{eqnarray*}
Now, we have completed the proof of Lemma \ref{xxyyyx}.\qed

The lemma shows \eqref{xxy} and \eqref{yyx}, and then \eqref{sigxi}
and \eqref{xisig}, which completes the proof of the theorem. \qed

\subsection{Triangular decomposition of paths}\label{U-V}

Here we consider the $x$-weighted lattice $L_1[n,k](x)$ 
($x=(\TY(x,i,j))$) 
and define the triangular decomposition
of the set of paths on $L_1[n,k]$. 

A path $p$ in $P_1[n,k]$ is called a {\it path above }(resp. 
{\it path below}) $(l,m)$ if 
for any point $(l,j)$ on $p$, we have $j>m$ (resp. $j<m$).
Set 
\begin{eqnarray*}
&&LP[n,k]^{(l)}_m:=\{p\in P_1[n,k]\,|\,p{\hbox{ is a path through }}(l,m)\},\\
&&LA[n,k]^{(l)}_m:=\{p\in P_1[n,k]\,|\,p{\hbox{ is a path above }}(l,m)\},\\
&&LB[n,k]^{(l)}_m:=\{p\in P_1[n,k]\,|\,p{\hbox{ is a path below }}(l,m)\},
\end{eqnarray*}
and define 
\[
\TY(R,l,m)(x):=\sum_{p\in LP[n,k]^{(l)}_m}x(p),\q
\TY(U,l,m)(x):=\sum_{p\in LA[n,k]^{(l)}_m}x(p),\q
\TY(V,l,m)(x):=\sum_{p\in LB[n,k]^{(l)}_m}x(p).
\]
It is obvious from the definitions that
\begin{equation}
\vep(x)=\TY(R,l,m)(x)+\TY(U,l,m)(x)+\TY(V,l,m)(x),\q
\TY(U,l-1,m)=\TY(U,l,m)+\TY(R,l,m),\q
\TY(V,l+1,m)=\TY(V,l,m)+\TY(R,l,m).
\label{EUV}
\end{equation}
\begin{lem}\label{path-zr}
We obtain the following formula:
\begin{eqnarray}
&&\TY(U,l,m+1)=\TY(U,l+1,m)+\TYS(X,l+1,m)\cdot \TY(X,l+1,m+1)
\frac{\TY(x,l+1,m)}{\TY(x,l+1,m+1)},\label{U=U+XX}\\
&&\TY(U,l-1,m)=\TY(U,l-1,m+1)+\TYS(X,l,m)\cdot \TY(X,l-1,m+1)
\label{U=U+XX2}\\
&&\TY(V,l+1,m+1)=\TY(V,l+1,m)+\TYS(X,l+1,m)\cdot \TY(X,l,m+1).
\label{V=V+XX}
\end{eqnarray}
\end{lem}
{\sl Proof.}
By the definition of $LA[n,k]^{(l)}_m$ we find that
$\TY({LA[n,k]},l+1,m)$ is a subset of $\TY({LA[n,k]},l,m+1)$ and 
their difference $\TY({LA[n,k]},l,m+1)\setminus \TY({LA[n,k]},l+1,m)$
coincides with 
the set of paths which go through $(l+1,m)$ and $(l+1,m+1)$,
which shows \eqref{U=U+XX}.
Considering similarly, the difference 
$\TY(LA,l-1,m)\setminus\TY(LA,l-1,m+1)$ is just the set of paths
going through $(l,m)$ and $(l-1,m+1)$, which implies
\eqref{U=U+XX2}.
Similarly, we obtain that 
$\TY({LB[n,k]},l+1,m)$ is a subset of $\TY({LB[n,k]},l+1,m+1)$ and 
their difference is same as the set of paths which 
go through $(l+1,m)$ and $(l,m+1)$. 
\qed
\renewcommand{\thesection}{\arabic{section}}
\section{Birational map $\osigma$}
\setcounter{equation}{0}
\renewcommand{\theequation}{\thesection.\arabic{equation}}

\subsection{Birational and Bi-positive isomorphism between $\cV_1$ and
$\cV_2$}

We shall define a birational and bi-positive isomorphism between $\cV_1$ and
$\cV_2$ by virtue of  the birational maps
$\Sigma$ and $\Xi$ introduced in the
previous section (Theorem \ref{birat}).
\begin{pro}\label{pro-bi}
Let $\cV_1$ and $\cV_2$ be geometric crystals introduced in \ref{v1v2}, 
and $\Sigma$ and $\Xi$ the birational maps as in the previous section. 
Define a rational map $\ovl\sigma:\cV_1\to\cV_2$ by 
$\ovl\sigma(v_1(x)):=v_2(\Sigma(x))$.
Then, the rational map $\ovl\sigma$ is a birational and bi-positive
map and its inverse is given as 
$\ovl\sigma^{-1}(v_2(y))=v_1(\Xi(y))$,
where $v_1(x)$ and $v_2(y)$ are as in \eqref{v1} and \eqref{v2}
respectively.
\end{pro}
{\sl Proof.}
Since we have the following commutative diagram:
\begin{equation}
\xymatrix{
\linethickness{20pt}
{(\bbC^\times)}^{kk'}\ar@{->}[r]_\Sigma\ar@{->}[d]_{v_1}&{(\bbC^\times)}^{kk'}
\ar@{->}[d]_{v_2}\\
\cV_1\ar@{->}[r]_{\ovl\sigma}&\cV_2
},
\end{equation}
and $v_1$, $v_2$ and $\Sigma$ are birational, it is clear that
 $\ovl\sigma$ is birational.
By the definition of positive structure of geometric crystal, 
bi-positivity  of $\Sigma$ means that of $\ovl\sigma$. \qed

\subsection{Properties of $\ovl\sigma$}

To induce an $\TY(A,1,n)$-geometric crystal structure on $\cV_1$, 
we shall see some crucial properties of $\ovl\sigma$.
Let $\cV_1=(\cV_1,\{e_i\},\{\gamma_i\},\{\vep_i\})_{i\in I_0}$
and $\cV_2=(\cV_1,\{\ovl e_i\},\{\ovl\gamma_i\},\{\ovl\vep_i\})_{i\in I_n}$
be the geometric crystals as above.
\begin{pro}\label{intertwine}
For $i\in\{1,2,\cd,n-1\}=I_0\cap I_n$, we have 
\begin{equation}
\ovl\sigma\circ e_i^c=\ovl e_i^c\circ\ovl\sigma,\q
\gamma_i=\ovl\gamma_i\circ\ovl\sigma,\q
\vep_i=\ovl\vep_i\circ\ovl\sigma
\end{equation}
\end{pro}
{\sl Proof.}
Let $y=\Sigma(x)$. Then,  we have 
\begin{equation}
\TY(y,c,m)\TY(y,c+1,m)\cd \TY(y,d,m)
=\TY(x,c+1,m)\TY(x,c+2,m)\cd \TY(x,d+1,m)
\frac{\TY(X,c,m)}{\TY(X,d+1,m)}\q
(1\leq  c\leq  d\leq k),
\label{yy=xxx}
\end{equation}
where if $(l,m)$ is out of the lattice $L_1[n,k]$, 
we understand $\TY(X,l,m)=1$.
Therefore, applying \eqref{yy=xxx} to 
\eqref{gamma-ep2} we have
\begin{eqnarray*}
\ovl\gamma_i(\ovl\sigma(v_1(x)))=
\frac{(\TY(x,\cc+1,i)\TY(x,\cc+2,i)\cd\TY(x,\dd+1,i)
\frac{\TY(X,\cc,i)}{\TY(X,\dd+1,i)})^2}{\TY(x,\cc+1,i-1)
\TY(x,\cc+2,i-1)\cd\TY(x,\dd+2,i-1)
\frac{\TY(X,\cc,i-1)}{\TY(X,\dd+2,i-1)}
\TY(x,\cc,i+1)
\TY(x,\cc+1,i+1)\cd\TY(x,\dd+1,i+1)
\frac{\TY(X,\cc-1,i+1)}{\TY(X,\dd+1,i+1)}}
\end{eqnarray*}
Here we set
\[
x_i(p,q):=\TY(x,p,i)\TY(x,p+1,i)\cd\TY(x,q,i),\q
y_i(p,q):=\TY(y,p,i)\TY(y,p+1,i)\cd\TY(y,q,i).
\]
Let us consider the following four cases (a)--(d):\\
(a) $k-i\geq 1$, $n-i\geq k$.
(b) $k-i<1$, $n-i\geq k$.
(c) $k-i< 1$, $n-i<  k$.
(d) $k-i\geq 1$, $n-i< k$.
\\

\nd (a) In this case, $\aa=k-i+1$, $\bb=k$, 
${\cc}=k-i$ and ${\dd}=k$.
Then we have 
\begin{eqnarray*}
\ovl\gamma_i(\ovl\sigma(v_1(x)))
&&=\frac{\left(x_i(k-i+1,k)\frac{\TY(X,k-i,i)}{\TY(X,k+1,i)}\right)^2}
{x_{i-1}(k-i+1,k)\frac{\TY(X,k-i+1,i-1)}{\TY(X,k+1,i-1)}
x_{i+1}(k-i,k)\frac{\TY(X,k-i-1,i+1)}{\TY(X,k+1,i+1)}}\\
&&=
\frac{\left(x_i(k-i+1,k)\TY(X,k-i,i)\right)^2}
{x_{i-1}(k-i+1,k)\TY(X,k-i+1,i-1)
x_{i+1}(k-i,k)\TY(X,k-i-1,i+1)}\\
&&=
\frac{\left(x_i(k-i+1,k)\frac{1}{\TY(x,1,n)}\right)^2}
{x_{i-1}(k-i+1,k)\frac{1}{\TY(x,1,n)}
x_{i+1}(k-i,k)\frac{1}{\TY(x,1,n)}}\\
&&=
\frac{x_i(k-i+1,k)^2}
{x_{i-1}(k-i+1,k)
x_{i+1}(k-i,k)}=\gamma_i(v_1(x)),
\end{eqnarray*}
where note that $\TY(X,k+1,m)=\TY(x,k+1,m)=1$ for any $m$ and $\TY(X,l,m)=
\frac{1}{\TY(x,1,n)}$ if $l+m=k$.

\nd
(b) In this case,  $\aa=1$, $\bb=k$, $\cc=1$ and $\dd=k$.
Then, we find that 
\begin{eqnarray*}
\ovl\gamma_i(\ovl\sigma(v_1(x)))&&=
\frac{\left(x_i(2,k)\TY(X,1,i)\right)^2}
{x_{i-1}(2,k)\TY(X,1,i-1)\cdot
x_{i+1}(2,k)\TY(X,1,i+1)}
\\
&&=
\frac{\left(x_i(2,k)\frac{\TY(x,1,i)}{\TY(x,1,n)}\right)^2}
{x_{i-1}(2,k)\frac{\TY(x,1,i-1)}{\TY(x,1,n)}x_{i+1}(2,k)\frac{\TY(x,1,i+1)}{\TY(x,1,n)}}
=\frac{x_i(1,k)^2}{x_{i-1}(1,k)x_{i+1}(1,k)}
=\gamma_i(v_1(x)),
\end{eqnarray*}
where note that $\TY(x,k+1,i-1)=1$.

\nd 
(c) In this case, $\aa=k-i+1$, $\bb=n-i+1$, 
${\cc}=k-i$ and ${\dd}=n-i$, 
Then, we find that 
\begin{eqnarray*}
\ovl\gamma_i(\ovl\sigma(v_1(x)))&&=
\frac{\left(x_i(k-i+1,n-i+1)\frac{\TY(X,k-i+1,i)}{\TY(X,n-i+1,i)}\right)^2}
{x_{i-1}(k-i+1,n-i+2)\frac{\TY(X,k-i+1,i-1)}{\TY(X,n-i+2,i-1)}
x_{i+1}(k-i,n-i+1)\frac{\TY(X,k-i-1,i+1)}{\TY(X,n-i,i+1)}}\\
&&=\frac{x_i(k-i+1,n-i+1)^2}{x_{i-1}(k-i+1,n-i+2)x_{i+1}(k-i,n-i+1)}
=\gamma_i(v_1(x)),
\end{eqnarray*}
where note that $\TY(x,n-i+1,i+1)=1$ and $\TY(X,l,m)=1$ if $l+m=n+1$.

\nd
(d) In this case,  $\aa=1$, $\bb=n-i+1$, $\cc=1$ and $\dd=n-i$.
Then, we find that 
\begin{eqnarray*}
&&\ovl\gamma_i(\ovl\sigma(v_1(x)))=
\frac{\left(x_i(2,n-i+1)\frac{\TY(X,1,i)}{\TY(X,n-i+1,i)}\right)^2}
{x_{i-1}(2,n-i+2)\frac{\TY(X,1,i-1)}{\TY(X,n-i+2,i-1)}
x_{i+1}(2,n-i)\frac{\TY(X,1,i+1)}{\TY(X,n-i,i+1)}}\\
&&=
\frac{\left(x_i(2,n-i+2)\frac{\TY(x,1,i)}{\TY(x,1,n)}\right)^2}
{x_{i-1}(2,n-i+2)\frac{\TY(x,1,i-1)}{\TY(x,1,n)}x_{i+1}(2,n-i+1)
\frac{\TY(x,1,i+1)}{\TY(x,1,n)}}
=\frac{x_i(1,n-i+1)^2}{x_{i-1}(1,n-i+2)x_{i+1}(1,n-i+1)}
=\gamma_i(v_1(x)),
\end{eqnarray*}
where note that $\TY(x,n-i+1,i+1)=1$. 
Thus, we obtain $\ovl\gamma_i(\ovl\sigma(v_1(x)))=\gamma_i(v_1(x))$ ($i=1,2,\cd, n-1$).

Next, we shall see $\ovl\vep_i(\ovl\sigma(x))=\vep_i(x)$ for $i=1,2,\cd, n-1$.
It follows from the explicit forms \eqref{Dli} and  \eqref{Eli} that 
\begin{equation}
\TY(E,l,i)(\ovl\sigma(v_1(x)))=
\TY(D,l+1,i)\frac{\TY(X,l,i)\TY(X,l+1,i)}{\TY(X,l+1,i-1)\TY(X,l,i+1)},
\label{E=DXX}
\end{equation}
where $\TY(D,l,i)=\TY(D,l,i)(v_1(x))$, $\TY(X,l,i)=\TY(X,l,i)(v_1(x))$, etc.
We get  $\TY(X,l+1,i-1)=\TY(X,l,i)+\frac{\TY(x,l+1,i-1)}{\TY(x,l+1,i)}\TY(X,l+1,i)$
by \eqref{x=x+x}. Thus, substituting this to \eqref{E=DXX}, we have
\begin{equation}
\frac{1}{\TY(E,l,i)(\ovl\sigma(v_1(x)))}=
\frac{\TY(X,l,i+1)}{\TY(D,l+1,i)\TY(X,l,i)\TY(X,l+1,i)}\left(
\TY(X,l,i)+\frac{\TY(x,l+1,i-1)}{\TY(x,l+1,i)}\TY(X,l+1,i)\right).
\label{1E1D}
\end{equation}
As above we consider the cases (a)--(d).

\nd (a) In this case, $\aa=k-i+1$, $\bb=k$, 
${\cc}=k-i$ and ${\dd}=k$.
\begin{eqnarray*}
&&\ovl\vep_i(\ovl\sigma(v_1(x)))=\sum_{l=k-i}^{k}\frac{1}{\TY(E,l,i)(\ovl\sigma(v_1(x)))}\\
&&=\frac{\TY(X,k-i,i+1)}{\TY(D,k-i+1,i)\TY(X,k-i+1,i)}
+\sum_{l=k-i+1}^{k-1}\left(\frac{\TY(X,l,i+1)}{\TY(D,l+1,i)\TY(X,l+1,i)}+
\frac{\TY(X,l,i+1)\TY(x,l+1,i-1)}{\TY(D,l+1,i)\TY(X,l,i)\TY(x,l+1,i)}\right)
+\frac{\TY(X,k,i+1)}{\TY(X,k,i)}\\
&&=\left(\frac{\TY(X,k-i,i+1)}{\TY(D,k-i+1,i)\TY(X,k-i+1,i)}+
\frac{\TY(X,k-i+1,i+1)\TY(x,k-i+2,i-1)}{\TY(D,k-i+2,i)\TY(X,k-i+1,i)\TY(x,k-i+2,i)}\right)
+\left(\frac{\TY(X,k,i+1)}{\TY(X,k,i)}+\frac{\TY(X,k-1,i+1)}{\TY(D,k,i)\TY(X,k,i)}\right)\\
&&+\sum_{l=k-i+2}^{k-1}\left(\frac{\TY(X,l-1,i+1)}{\TY(D,l,i)\TY(X,l,i)}+
\frac{\TY(X,l,i+1)\TY(x,l+1,i-1)}{\TY(D,l+1,i)\TY(X,l,i)\TY(x,l+1,i)}\right)
=\frac{1}{\TY(D,k-i+1,i)}+\frac{1}{\TY(D,k,i)}+\sum_{l=k-i+2}^{k-1}\frac{1}{\TY(D,l,i)}\\
&&=\sum_{l=k-i+1}^{k}\frac{1}{\TY(D,l,i)}=\vep_i(v_1(x)),
\end{eqnarray*}
where for the 4th equality we use \eqref{x=x+x}, $\TY(D,k,i)=\frac{\TY(x,k,i)}{\TY(x,k,i+1)}$
and $\TY(D,l+1,i)=\frac{\TY(x,l+1,i-1)\TY(x,l,i+1)}{\TY(x,l,i)\TY(x,l+1,i)}\TY(D,l,i)$.

\nd(b) $\aa=1$, $\bb=k$, $\cc=1$ and $\dd=k$. 
\begin{eqnarray*}
&&\ovl\vep_i(\ovl\sigma(v_1(x)))=\sum_{l=1}^{k}\frac{1}{\TY(E,l,i)(\ovl\sigma(v_1(x)))}
=\sum_{l=1}^{k-1}\left(\frac{\TY(X,l,i+1)}{\TY(D,l+1,i)\TY(X,l+1,i)}+
\frac{\TY(X,l,i+1)\TY(x,l+1,i-1)}{\TY(D,l+1,i)\TY(X,l,i)\TY(x,l+1,i)}\right)
+\frac{\TY(X,k,i+1)}{\TY(X,k,i)}\\
&&=\frac{\TY(x,2,i-1)\TY(X,1,i+1)}{\TY(x,2,i)\TY(D,2,i)\TY(X,1,i)}+
\left(\frac{\TY(X,k,i+1)}{\TY(X,k,i)}+\frac{\TY(X,k-1,i+1)}{\TY(D,k,i)\TY(X,k,i)}\right)
+\sum_{l=2}^{k-1}\left(\frac{\TY(X,l-1,i+1)}{\TY(D,l,i)\TY(X,l,i)}+
\frac{\TY(X,l,i+1)\TY(x,l+1,i-1)}{\TY(D,l+1,i)\TY(X,l,i)\TY(x,l+1,i)}\right)\\
&&=\frac{1}{\TY(D,1,i)}+\frac{1}{\TY(D,k,i)}+\sum_{l=2}^{k-1}\frac{1}{\TY(D,l,i)}
=\sum_{l=1}^{k}\frac{1}{\TY(D,l,i)}=\vep_i(v_1(x)),
\end{eqnarray*}
where for the 4th equality, we also use 
\eqref{x=x+x}, $\TY(D,k,i)=\frac{\TY(x,k,i)}{\TY(x,k,i+1)}$
and $\TY(D,l+1,i)=\frac{\TY(x,l+1,i-1)\TY(x,l,i+1)}{\TY(x,l,i)\TY(x,l+1,i)}\TY(D,l,i)$.

For the cases (c)  and (d), we can show $\ovl\vep_i(\ovl\sigma(v_1(x)))=\vep_i(v_1(x))$
similarly.

Finally, let us show $\osigma\circ e^\al_i=\ovl e^\al_i\circ \osigma$ $(i=1,2,\cd, n-1)$.
Set $v_2(y')=\osigma(e^\al_m(v_1(x)))$ and $v_2(y''):=\ovl e^\al_m(\osigma(v_1(x)))$ 
($m=1,2,\cd,n-1$) and show $\TY(y',i,j)=\TY(y'',i,j)$ for any $i\in \{1,2,\cd,k\}$ 
and $j=\cc,\cc+1,\cd,\dd$.

\nd
We shall see the case $j\ne m$.
In this case, since $\TY(y',l,j)=\TY(x,l+1,j)\frac{\TY(X,l,j)(x')}{\TY(X,l+1,j)(x')}$
and 
$\TY(y'',l,j)=\TY(y,l,j)=\TY(x,l+1,j)\frac{\TY(X,l,j)(x)}{\TY(X,l+1,j)(x)}$,
it suffices to show $\TY(X,l,j)(x)=\TY(X,l,j)(x')$ for $j\ne m$.
Here, suppose $j>m$. In this case, 
$\TY(X,l,j)$ does not depend on $\TY(x,l,m)$'s 
by its definition and then we find that 
$\TY(X,l,j)(x)=\TY(X,l,j)(x')$ for $j>m$.

Next, suppose $j=m-1$. We set 
\begin{equation}
\TY(A,l,i):=\displaystyle\sum_{p=\cc}^{l-1}\frac{1}{\TY(E,p,i)}+
\displaystyle\sum_{p=l}^{\dd}\frac{c}{\TY(E,p,i)},\qq\q
\TY(B,l,i):=\displaystyle\sum_{p=\aa}^{l-1}\frac{1}{\TY(D,p,i)}+
\displaystyle\sum_{p=l}^\bb\frac{c}{\TY(D,p,i)}.
\label{AB}
\end{equation}
Let us denote $\TY(X,l,j)(x')$ by $\TY(X',l,j)$. 
Applying the formula \eqref{x=x+x} , we get
\begin{eqnarray}
&&
\TY(X',l,m-1)=\TY(X',l-1,m)+\frac{\TY(x',l,m-1)}{\TY(x',l,m)}\TY(X',l,m)
=\left(\TY(X',l-2,m+1)+\frac{\TY(x',l-1,m)}{\TY(x',l-1,m+1)}\TY(X',l-1,m+1)\right)
+\frac{\TY(x',l,m-1)}{\TY(x',l,m)}
\left(\TY(X',l-1,m+1)+\frac{\TY(x',l,m)}{\TY(x',l,m+1)}\TY(X',l,m+1)\right)\nn\\
&&\qq\q=\left(\TY(X,l-2,m+1)+\frac{\TY(x,l,m-1)}{\TY(x,l,m+1)}\TY(X,l,m+1)\right)
+\TY(X,l-1,m+1)
\left(\frac{\TY(x',l-1,m)}{\TY(x,l-1,m+1)}+\frac{\TY(x,l,m-1)}{\TY(x',l,m)}\right).
\label{ww}
\end{eqnarray}
Now, we have 
\begin{eqnarray*}
&&\frac{\TY(x',l-1,m)}{\TY(x,l-1,m+1)}+\frac{\TY(x,l,m-1)}{\TY(x',l,m)}
=\frac{1}{\TY(B,l,m)}\left(\frac{\TY(x,l-1,m)}{\TY(x,l-1,m)}\TY(B,l-1,m)+
\frac{\TY(x,l,m-1)}{\TY(x,l,m)}\TY(B,l+1,m)\right)\\
&&=\frac{1}{\TY(B,l,m)}\left(\frac{\TY(x,l-1,m)}{\TY(x,l-1,m)}\left(\TY(B,l,m)
+\frac{\al}{\TY(D,l-1,m)}-\frac{1}{\TY(D,l-1,m)}\right)+
\frac{\TY(x,l,m-1)}{\TY(x,l,m)}\left(\TY(B,l,m)
-\frac{\al}{\TY(D,l,m)}+\frac{1}{\TY(D,l,m)}\right)
\right)\\
&&=\left(\frac{\TY(x,l-1,m)}{\TY(x,l-1,m+1)}+\frac{\TY(x,l,m-1)}{\TY(x,l,m)}\right)
+(\al-1)\left(\frac{\TY(x,l-1,m)}{\TY(x,l-1,m+1)\TY(D,l-1,m)}
-\frac{\TY(x,l,m-1)}{\TY(x,l,m)\TY(D,l,m)}\right)
=\frac{\TY(x,l-1,m)}{\TY(x,l-1,m+1)}+\frac{\TY(x,l,m-1)}{\TY(x,l,m)},
\end{eqnarray*}
where for the last equality we use the relation
$\TY(D,l-1,m)=\frac{\TY(x,l-1,m)\TY(x,l,m)}{\TY(x,l,m-1)\TY(x,l-1,m+1)}\TY(D,l,m)$.
Thus, substituting this to \eqref{ww}, we obtain 
\[
\TY(X',l,m-1)=
\left(\TY(X,l-2,m+1)+\frac{\TY(x,l,m-1)}{\TY(x,l,m+1)}\TY(X,l,m+1)\right)
+\TY(X,l-1,m+1)
\left(\frac{\TY(x,l-1,m)}{\TY(x,l-1,m+1)}+\frac{\TY(x,l,m-1)}{\TY(x,l,m)}\right)
=\TY(X,l,m-1).
\]
Let us show the cases $j<m-1$ using descending induction on $j$. Applying \eqref{x=x+x}
and the induction hypothesis, we have
\[
\TY(X',l,j)=\TY(X',l-1,j+1)+\frac{\TY(x',l,j)}{\TY(x',l,j+1)}\TY(X',l,j+1)
=\TY(X,l-1,j+1)+\frac{\TY(x,l,j)}{\TY(x,l,j+1)}\TY(X,l,j+1)=\TY(X,l,j),
\]
which completes the proof of  $\TY(X',l,j)=\TY(X,l,j)$ for $j\ne m$ and then 
$\TY(y',l,j)=\TY(y'',l,j)$ for $j\ne m$.

Finally, let us show the case $m=j$, that is, 
$\TY(y',l,m)=\TY(y'',l,m)$ for $l=1,2,\cd,k$.
Let $\TY(A,l,i)$ and $\TY(B,l,i)$ be as in \eqref{AB}. 
By the explicit forms of $e^\al_m$, $\ovl e^\al_m$ and $\osigma$, we have 
\begin{equation}
\TY(y',l,m)=\TY(x,l+1,m)\frac{\TY(B,l+1,m)\TY(X',l,m)}{\TY(B,l+2,m)\TY(X',l+1,m)},\qq
\TY(y'',l,m)=\TY(x,l+1,m)\frac{\TY(A,l,m)(\osigma(x))\TY(X,l,m)}{\TY(A,l+1,m)(\osigma(x))
\TY(X,l+1,m)}
\label{yyy}
\end{equation}
Thus, we may show 
\begin{equation}
\TY(B,l+1,m)\TY(X',l,m)=\TY(A,l,m)(\osigma(x))\TY(X,l,m)
\label{BXAX}
\end{equation}
Suppose $\al=1$. Since we have $\TY(B,l+1,m)(x)|_{\al=1}=\vep_m(x)$, 
$\TY(A,l,m)(\sigma(x))|_{\al=1}=\ovl\vep_m(\osigma(x))$, $\TY(X',l,m)=\TY(X,l,m)$ and 
$\vep_m(x)=\ovl\vep_m(\osigma(x))$, we know that \eqref{BXAX} holds for the case $\al=1$.

Suppose that $\al$ is generic.
Since for $j>m$ we have $\TY(X',l,j)=\TY(X,l,j)$, we get 
\begin{eqnarray}
&&\TY(B,l+1,m)\TY(X',l,m)=
\TY(B,l+1,m)\left(\frac{\TY(x',l,m)}{\TY(x,l,m+1)}\TY(X,l,m+1)+\TY(X,l-1,m+1)\right)
\nn \\ 
&&=\TY(B,l+1,m)\left(\frac{\TY(x,l,m)\TY(B,l,m)}{\TY(x,l,m+1)\TY(B,l+1,m)}
\TY(X,l,m+1)+\TY(X,l-1,m+1)\right)
=\frac{\TY(x,l,m)}{\TY(x,l,m+1)}\TY(B,l,m)\TY(X,l,m+1)+\TY(B,l+1,m)\TY(X,l-1,m+1)\nn\\
&&=\frac{\TY(x,l,m)}{\TY(x,l,m+1)}\TY(X,l,m+1)\left(
\sum_{p=\aa}^{l-1}\frac{1}{\TY(D,p,m)}+\sum_{p=l}^\bb\frac{\al}{\TY(D,p,m)}\right)
+\TY(X,l-1,m)\left(
\sum_{p=\aa}^{l}\frac{1}{\TY(D,p,m)}+\sum_{p=l+1}^\bb\frac{\al}{\TY(D,p,m)}\right),
\label{BXD}\\
&&\TY(A,l,m)(\osigma(x))\TY(X,l,m)=
\TY(X,l,m)\left(\sum_{p=\cc}^{l-1}\frac{1}{\TY(E,p,i)(\osigma(x))}+
\sum_{p=l}^{\dd}\frac{\al}{\TY(E,p,i)(\osigma(x))}\right)\nn \\
&&=\TY(X,l,m)\left(\sum_{p=\cc}^{l-1}
\frac{\TY(X,p,m+1)\TY(X,p+1,m-1)}{\TY(D,p+1,m)\TY(X,p,m)\TY(X,p+1,m)}
+\sum_{p=l}^{\dd}
\frac{\al \TY(X,p,m+1)\TY(X,p+1,m-1)}{\TY(D,p+1,m)\TY(X,p,m)\TY(X,p+1,m)}
\right).
\label{AXD}
\end{eqnarray}
Therefore, since we obtained \eqref{BXAX} for $\al=1$, 
we may show the coefficients of $\al$ in \eqref{BXD} and \eqref{AXD} 
coincide each other, that is, 
\begin{equation}
\frac{\TY(x,l,m)}{\TY(x,l,m+1)}\TY(X,l,m+1)
\sum_{p=l}^\bb\frac{1}{\TY(D,p,m)}
+\TY(X,l-1,m)\sum_{p=l+1}^\bb\frac{1}{\TY(D,p,m)}
=\TY(X,l,m)\sum_{p=l}^{\dd}
\frac{\TY(X,p,m+1)\TY(X,p+1,m-1)}{\TY(D,p+1,m)\TY(X,p,m)\TY(X,p+1,m)}.
\label{al-co}
\end{equation}
We shall show \eqref{al-co} by descending induction on $l$.
Suppose $l=\dd$. If $k\leq n-m$, we have $\l=\dd=\min(k,n-m)=k$ and 
$\bb=\min(k,n-m+1)=k=l$, and then 
\begin{eqnarray*}
&&{\hbox{L.H.S. of \eqref{al-co}}}=\frac{\TY(x,k,m)\TY(X,k,m)}{\TY(x,k,m+1)\TY(D,k,m)}=\TY(X,k,m+1),
\\
&&{\hbox{R.H.S. of \eqref{al-co}}}=\TY(X,k,m)\frac{\TY(X,k,m+1)\TY(X,k+1,m-1)}{\TY(D,k+1,m)
\TY(X,k,m)\TY(x,k+1,m)}=\TY(X,k,m+1).
\end{eqnarray*}
In the case $k>n-m$, we have $l=\dd=n-m$ and  $\bb=n-m+1$,  and then
\begin{eqnarray*}
&&{\hbox{L.H.S. of \eqref{al-co}}}=
\frac{\TY(x,n-m,m)\TY(X,n-m,m+1)}{\TY(x,n-m,m+1)}
\left(\frac{1}{\TY(D,n-m,m)}+\frac{1}{\TY(D,n-m+1,m)}\right)
+\frac{\TY(X,n-m-1,m+1)}{\TY(D,n-m+1,m)}\\
&&=\frac{\TY(X,n-m,m)}{\TY(D,n-m+1,m)}
+\frac{\TY(x,n-m,m)\TY(X,n-m,m+1)}{\TY(x,n-m,m)\TY(D,n-m,m)}
=\frac{1}{\TY(D,n-m+1,m)}\left(\TY(X,n-m,m)
+\frac{\TY(x,n-m+1,m-1)\TY(X,n-m+1,m)}{\TY(x,n-m+1,m)}\right)
=\frac{\TY(X,n-m+1,m-1)}{\TY(D,n-m+1,m)},\\
&&{\hbox{R.H.S. of \eqref{al-co}}}=
\TY(X,n-m,m)\frac{\TY(X,n-m,m+1)\TY(X,n-m+1,m-1)}{\TY(D,n-m+1,m)\TY(X,n-m,m)\TY(X,n-m+1,m)}
=\frac{\TY(X,n-m+1,m-1)}{\TY(D,n-m+1,m)},
\end{eqnarray*}
where note that $\TY(X,n-m+1,m)=\TY(X,n-m,m+1)=1$.
Thus, we obtain \eqref{al-co} for $l=\dd$.
Here, suppose $l<\dd$. By the induction hypothesis, we have 
\begin{eqnarray*}
&&{\hbox{R.H.S. of \eqref{al-co}}}=
\TY(X,l,m)\frac{\TY(X,l,m+1)\TY(X,l+1,m-1)}{\TY(D,l+1,m)\TY(X,l,m)\TY(X,l+1,m)}
+\frac{\TY(X,l,m)}{\TY(X,l+1,m)}\times{\hbox{L.H.S. of \eqref{al-co} for $l+1$}}\\
&&=\frac{\TY(X,l,m+1)\TY(X,l+1,m-1)}{\TY(D,l+1,m)\TY(X,l+1,m)}+\TY(X,l,m)
\left(\frac{\TY(x,l+1,m)\TY(X,l+1,m)}{\TY(x,l+1,m)\TY(D,l+1,m)\TY(X,l+1,m)}+
\sum_{p=l+2}^\bb\frac{1}{\TY(D,p,m)}\right)\\
&&=\frac{\TY(x,l,m)\TY(X,l,m+1)}{\TY(x,l,m+1)\TY(D,l,m)}+\TY(X,l,m)\sum_{p=l+1}^\bb
\frac{1}{\TY(D,p,m)}=
{\hbox{L.H.S. of \eqref{al-co}}},
\end{eqnarray*}
which shows \eqref{al-co} and then we obtain $\TY(y',l,m)=\TY(y'',l,m)$, that is, 
$\osigma\circ e^\al_i=\ovl e^\al_i\circ\osigma$ for $i=1,2,\cd,n-1$.
Now, we have completed the proof of Proposition \ref{intertwine}.\qed

\renewcommand{\thesection}{\arabic{section}}
\section{Affine geometric crystal $\cV(\TY(A,1,n))$}
\setcounter{equation}{0}
\renewcommand{\theequation}{\thesection.\arabic{equation}}

\subsection{0-structures on  $\cV_1$}\label{ss0}
In this subsection, we induce the $0$-structures on $\ge_0$-geometric crystal 
$\cV_1$. Since $\cV_2$ is a $\ge_n$-geometric crystal, it is equipped with 
the $0$-structure, on the other hand, $\cV_1$ does not hold 
$0$-structure. We shall, however, define $0$-structure on $\cV_1$ through 
the birational isomorphism $\osigma$ by:
\begin{equation}
\gamma_0:=\ovl\gamma_0\circ\osigma,\q
\vep_0:=\ovl\vep_0\circ\osigma,\q
e_0^c:=\osigma^{-1}\circ\ovl e_i^c\circ\osigma
\label{0st}
\end{equation}
The following is one of the main theorems 
of this paper, which will be shown in the later
subsections.
\begin{thm}\label{0-str}
$\ge_0$-geometric crystal $\cV_1$ becomes an 
$\TY(A,1,n)$-geometric crystal 
if it is equipped with $0$-structures as in \eqref{0st}, which will be denoted 
by $\cV(\TY(A,1,n))$. Furthermore, 
$v_1:(\bbC^\times)^{kk'}\to \cV(\TY(A,1,n))(=\cV_1)$ 
gives a positive structure.

\end{thm}
Before showing this theorem, let us see the explicit forms of 
$\gamma_0, \vep_0$ and $e_0^c$.
\begin{pro}\label{exp-0}
Let us set $v_1(x')=e^c_0(v_1(x))$ $(c\in(\bbC^\times)$. 
Then we have the following:
\begin{eqnarray}
&&\gamma_0(v_1(x))=\frac{1}{\TY(x,1,n)\TY(x,k,1)},\qq
\vep_0(v_1(x))=\TY(x,1,n)\vep(x),\label{g0ep0}\\
&&\TY(x',l,m)=\begin{cases}
\displaystyle 
\TY(x,l,m)\frac{\TY(\al,l,m)}{\TY(\al,l+1,m)}&{\rm if }\,\,
(l,m)\ne(1,n),\\
\displaystyle 
\frac{\TY(x,1,n)}{c}&{\rm if}\,\,(l,m)=(1,n).
\end{cases}
\label{e0}
\end{eqnarray}
where $\TY(\al,l,m)=\TY(\al,l,m)(c)=\TY(U,l-1,m)+c\TY(V,l,m)$, and 
$\TY(U,l,m)$ and 
$\TY(V,l,m)$ are as in \ref{U-V}.

\end{pro}
{\sl Proof.}
First, let us show $\gamma_0(v_1(x))=\ovl\gamma_0(\osigma(v_1(x)))$.
Set $v_2(y)=v_2((\TY(y,i,j)))=\osigma(v_1(x))$.  
By the explicit form of $\osigma$ we know that
$\TY(y,l,m)=\TY(x,l+1,m)\frac{\TY(X,l,m)}{\TY(X,l+1,m)}$ and then have
\begin{eqnarray*}
\ovl\gamma_0(\osigma(v_1(x)))=\ovl\gamma_0(v_2(y))=
\frac{{\TY(y,k,0)}^2}{\TY(y,k-1,1)\TY(y,k,1)}
=\frac{\frac{1}{{\TY(x,1,n)}^2}}{\left(\frac{\TY(x,k,1)}{\TY(x,1,n)\TY(X,k,1)}\right)
\TY(X,k,1)}
=\frac{1}{\TY(x,1,n)\TY(x,k,1)}=\gamma_0(v_1(x)).
\end{eqnarray*}
Next, let us show 
$\vep_0(v_1(x))=\ovl\vep_0(\osigma(v_1(x)))$.
By \eqref{gamma-ep2}, we get $\ovl\vep_0(v_2(y))=\frac{\TY(y,k,1)}{\TY(y,k,0)}$
and then 
\[
\ovl\vep_0(v_2(y))=\frac{\TY(y,k,1)}{\TY(y,k,0)}
=\frac{\TY(X,k,1)}{\frac{1}{\TY(x,1,n)}}=
\TY(x,1,n)\TY(X,k,1)=\TY(x,1,n)\vep(x)=\vep_0(v_1(x)).
\]
Finally, it is sufficient to show $\ovl e^c_0\circ\osigma=\osigma\circ e_0^c$.
Set $v_2(y')=\ovl e^c_0\circ\osigma(x)(=\ovl e_0^c(y))$ and 
$v_2(y'')=\osigma\circ e_0^c(x)(=\osigma(x'))$. 
We have $\TY(y',k,0)=c\TY(y,k,0)=\frac{c}{\TY(x,1,n)}$ 
and $\TY(y'',k,0)=\frac{1}{\TY(x',1,n)}=\frac{c}{\TY(x,1,n)}$, which
shows the case $(l,m)=(k,0)$.

Now, we define the total order $\prec$ 
on the set of coordinates $F:=\{(l,m)|
1\leq l\leq k\leq l+m\leq n\}\setminus\{(k,0)\}$ of the lattice 
$L_2[n,k]$ by setting $(l,m)\prec(l',m')\Leftrightarrow l<l'$,  and
$m>m'$ if $l=l'$. Thus, $(1,n-1)$ is the minimum and $(k,1)$ is the maximum.
We shall show $\TY(y',l,m)=\TY(y'',l,m)$ by the induction according to $\prec$.
Let us see the case $(l,m)=(1,n-1)$.
\begin{eqnarray*}
&&\TY(y',1,n-1)=\TY(y,1,n-1)=\TY(x,2,n-1)\frac{\TY(X,1,n-1)}{\TY(X,2,n-1)}
=\frac{\TY(x,1,n-1)\TY(x,2,n-1)}{\TY(x,1,n)},\\
&&\TY(y'',1,n-1)=\frac{\TY(x',1,n-1)\TY(x',2,n-2)}{\TY(x',1,n)}
=\frac{\TY(x,1,n-1)\frac{\TY(\al,1,n-1)}{\TY(\al,2,n-1)}\cdot
\TY(x,2,n-1)\frac{\TY(\al,2,n-1)}{\TY(\al,3,n-2)}}{\frac{\TY(x,1,n)}{c}}
=\frac{\TY(x,1,n-1)\vep(x)\cdot\frac{\TY(x,2,n-1)}{c \vep(x)}}
{\frac{\TY(x,1,n)}{c}}
=\frac{\TY(x,1,n-1)\TY(x,2,n-1)}{\TY(x,1,n)},
\end{eqnarray*}
which shows $\TY(y',1,n-1)=\TY(y'',1,n-1)$.

Let us see $(l,m)\succ (1,n-1)$.  We assume 
$\TY(y',q,r)=\TY(y'',q,r)$ for $q<l$ and any $r$.
Thus, multiplying both sides for $q=1,2,\cd,l-1$  we get 
$\TY(x',2,r)\TY(x',3,r)\cd\TY(x',l,r)
\frac{\TY(X',1,r)}{\TY(X',l,r)}
=\TY(x,2,r)\TY(x,3,r)\cd \TY(x,l,r)
\frac{\TY(X,1,r)}{\TY(X,l,r)}$
and then for $r=m,m+1$
\begin{eqnarray}
&&\frac{\TY(\al,2,m)\TY(X',1,m)}{\TY(\al,l+1,m)\TY(X',l,m)}
=\frac{\TY(X,1,m)}{\TY(X,l,m)},\qq
\frac{\TY(\al,2,m+1)\TY(X',1,m+1)}{\TY(\al,l+1,m+1)\TY(X',l,m+1)}
=\frac{\TY(X,1,m+1)}{\TY(X,l,m+1)}.
\label{xxx=xxx}
\end{eqnarray}
Substituting $\TY(X',1,m+1)=\frac{\TY(x',1,m+1)}{\TY(x',1,n)}
=\frac{c\TY(x,1,m+1)\TY(\al,1,m+1)}{\TY(x,1,n)\TY(\al,2,m+1)}$ to 
\eqref{xxx=xxx}, we obtain 
\begin{equation}
\TY(X',l,m+1)=\frac{c\TY(\al,1,m+1)\TY(X,l,m+1)}{\TY(\al,l+1,m+1)}.
\label{xcal}
\end{equation}
Since we also assume that $\TY(y',l,p)=\TY(y'',l,p)$ for $p>m$, 
we have $\TY(x',l+1,m+1)\frac{\TY(X',l,m+1)}{\TY(X',l+1,m+1)}
=\TY(x,l+1,m+1)\frac{\TY(X,l,m+1)}{\TY(X,l+1,m+1)}$. Substituting this to 
$\TY(X',l+1,m)=\TY(X',l,m+1)
+\frac{\TY(x',l+1,m)}{\TY(x',l+1,m+1)}\TY(X',l+1,m+1)$, we find that 
\begin{equation}
\TY(X',l+1,m)=\TY(X',l,m+1)\left(1+\frac{\TY(x',l+1,m)\TY(X,l+1,m+1)}
{\TY(x,l+1,m+1)\TY(X,l,m+1)}\right).
\label{44}
\end{equation}
Now, applying \eqref{xcal} to \eqref{44}, we obtain
\begin{equation}
\TY(X',l+1,m)=\frac{c\TY(\al,1,m+1)\TY(X,l,m+1)}{\TY(\al,l+1,m+1)}
\left(1+\frac{\TY(x',l+1,m)\TY(X,l+1,m+1)}
{\TY(x,l+1,m+1)\TY(X,l,m+1)}\right).
\label{99}
\end{equation}
Using \eqref{xxx=xxx} and \eqref{99}, we get 
\begin{equation}
\TY(y'',l,m)=
\TY(x',l+1,m)\frac{\TY(X',l,m)}{\TY(X',l+1,m)}
=\TY(x,l+1,m)\frac{\TY(\al,l+1,m)\TY(X,l,m)}{\TY(\al,l+2,m)\TY(X,l,m+1)+
\frac{\TY(x,l+1,m)}{\TY(x,l+1,m+1)}\TY(\al,l+1,m)\TY(X,l+1,m+1)},
\label{yyxal}
\end{equation}
where for the last equality we use $\TY(\al,1,m)=\TY(\al,1,m+1)=\vep$.
Here, 
\begin{eqnarray}
&&\TY(\al,l+2,m)\TY(X,l,m+1)+
\frac{\TY(x,l+1,m)}{\TY(x,l+1,m+1)}\TY(\al,l+1,m)\TY(X,l+1,m+1)\nn\\
&&=(\TY(U,l+1,m)+c\TY(R,l+1,m)+c\TY(V,l+1,m))\TY(X,l+1,m)+
\frac{\TY(x,l+1,m)}{\TY(x,l+1,m+1)}
(\TY(U,l+1,m)+\TY(R,l+1,m)+c\TY(V,l+1,m))\TY(X,l+1,m+1)\nn\\
&&=(\TY(U,l+1,m)+c\TY(R,l+1,m))\left(\TY(X,l,m+1)+
\frac{\TY(x,l+1,m)}{\TY(x,l+1,m+1)}\TY(X,l+1,m+1)\right)+
\TY(R,l+1,m)(c\TY(X,l,m+1)+\frac{\TY(x,l+1,m)}{\TY(x,l+1,m+1)}\TY(X,l+1,m+1))
\nn\\
&&=(\TY(U,l+1,m)+c\TY(R,l+1,m))\TY(X,l+1,m)+
\TYS(X,l+1,m)\TY(X,l+1,m)(c\TY(X,l,m+1)
+\frac{\TY(x,l+1,m)}{\TY(x,l+1,m+1)}\TY(X,l+1,m+1))\nn\\
&&=\TY(X,l+1,m)\left(\TY(U,l+1,m)
+\frac{\TY(x,l+1,m)}{\TY(x,l+1,m+1)}\TYS(X,l+1,m)\TY(X,l+1,m+1)
+c(\TY(R,l+1,m)+\TYS(X,l+1,m)\TY(X,l,m+1))\right)\nn\\
&&=\TY(X,l+1,m)(\TY(U,l,m+1)+c\TY(V,l+1,m+1))=
\TY(X,l+1,m)\TY(\al,l+1,m+1),
\label{denomi}
\end{eqnarray}
where for the third equality we use the formula \eqref{x=x+x} and 
formula $\TY(R,l+1,m)=\TYS(X,l+1,m)\TY(X,l+1,m)$ and, 
for the last equality we use the formula \eqref{U=U+XX} and 
\eqref{V=V+XX} in Lemma \ref{EUV}. 
Applying this to \eqref{yyxal}, we obtain 
\[
\TY(y'',l,m)=\TY(x,l+1,m)
\frac{\TY(X,l,m)\TY(\al,l+1,m+1)}{\TY(X,l+1,m)\TY(\al,l+1,m+1)}
=\TY(x,l+1,m)\frac{\TY(X,l,m)}{\TY(X,l+1,m)}=\TY(y,l,m)=\TY(y',l,m),
\]
which shows $\osigma\circ e_0^c=\ovl e_0^c\circ\osigma$.\qed
\subsection{Proof of Theorem \ref{0-str}}

In the rest of the section, for simplicity, let us denote 
$v_1(x)$ by $x$ and $v_2(y)$ by $y$, e.g., $\vep_i(x)$ for $\vep_i(v_1(x))$.

The positivity of $v_1$ on $\cV(\TY(A,1,n))$ 
is trivial from the explicit forms of 
$\gamma_i$, $\vep_i$ and $\eit^c$ $(i\in I$) as in Proposition \ref{V1}
and Proposition \ref{exp-0}.

To prove Theorem \ref{0-str}, it suffices to show the relations 
of the geometric crystals related to $0$-structures, that is, 
\begin{eqnarray}
&&\gamma_0(e_i^c(x))=c^{a_{i0}}\gamma_0(x),\q
\gamma_i(e_0^c(x))=c^{a_{0i}}\gamma_i(x),\q(i=0,1,\cd,n),
\label{0gamma0}\\
&&\vep_0(e^c_0(x))=c^{-1}\vep_0(x),\label{ep0}\\
&&e_0^ce_i^d=e_i^de_0^c,\q(i=2,3,\cd, n-1),\label{0i}\\
&&e_0^ce_1^{cd}e_0^d=e_1^de_0^{cd}e_1^c,\label{010}\\
&&e_0^ce_n^{cd}e_0^d=e_n^de_0^{cd}e_n^c.\label{0n0}
\end{eqnarray}
As for \eqref{0gamma0}, by Proposition \ref{intertwine} and \eqref{0st},
for $i=0,1,2,\cd,n-1$ we have
\[
\gamma_0(e_i^c(x))=\ovl\gamma_0(\osigma(e_i^c(x)))=
\ovl\gamma_0(\ovl e_i^c\osigma(x))=c^{a_{i0}}\ovl\gamma_0(\osigma(x))
=c^{a_{i0}}\gamma_0(x).
\]
Similarly, for $i=1,\cd,n-1$ we obtain
\[
\gamma_i(e_0^c(x))=\gamma_i(\osigma^{-1}\circ\ovl e_i^c\circ\osigma(x))
=\ovl\gamma_i(\ovl e_i^c\circ\osigma(x))=c^{a_{0i}}\ovl\gamma_i(\osigma(x))
=c^{a_{0i}}\gamma_i(x).
\]
For $i=n$ case, by the explicit form of $\gamma_0$ we have 
\[
\gamma_0(e_n^c(x))=\frac{1}{c\TY(x,1,n)\TY(x,k,1)}=c^{-1}\gamma_0(x).
\]
Similarly, by the explicit form of $\gamma_n$ we get
\[
\gamma_n(e_0^c(x))=\frac{\left(\frac{\TY(x,1,n)}{c}\right)^2}
{\frac{\TY(\al,1,n-1)}{\TY(\al,2,n-1)}\TY(x,1,n-1)
\frac{\TY(\al,2,n-1)}{\TY(\al,3,n-1)}\TY(x,2,n-1)}
=\frac{{\TY(x,1,n)}^2}{c\TY(x,1,n-1)\TY(x,2,n-1)}=c^{-1}\gamma_n(x),
\]
where note that $\TY(\al,1,n-1)=\vep$ and $\TY(\al,3,n-1)=c\vep$.

As for \eqref{ep0}, due to \eqref{0st} and the fact $\ovl\vep_0(\ovl e_0^c(v_2(y)))
=c^{-1}\ovl\vep_0(v_2(y))$ we have
\[
\vep_0(e^c_0(x))=\ovl\vep_0(\osigma(e_0^c((x))))=
\ovl\vep_0(\ovl e_0^c(\osigma(x)))
=c^{-1}\ovl\vep_0(\osigma(x))=c^{-1}\vep_0(x).
\]
Since we have
$\ovl e_0^c\ovl e_i^d=\ovl e_i^d\ovl e_0^c$ ($i=2,\cd, n-1$) and 
$\ovl e_0^c\ovl e_1^{cd}\ovl e_0^d =\ovl e_1^d\ovl e_0^{cd}\ovl e_1^c$ 
on $\cV_2$, by using Proposition \ref{intertwine} and 
\eqref{0st} we get \eqref{0i} and \eqref{010}.
Thus, the remaining case to show 
is \eqref{0n0}, which is the most difficult one since
we have to show it by direct calculations.

\subsection{Formula for Functions on Paths }
In this subsection, for a function or a map $f$ on $\cV_1$, 
we shall denote $f(e_0^d(x))$ $(d\in \bbC^\times$) 
by $\ovl f(x)$ for simplicity (except for $\ovl e_i,\ovl\gamma_i,\ovl\vep_i$).
Then, e.g., $\TY(\ovl X,l,m)$ implies $\TY(X,l,m)(e_0^d(x))$.

We shall give several formula for the functions related to the actions of $e_0^c$.
\begin{lem}
We have the following formula:
\begin{eqnarray}
&&\TYO(R,1,n-1)=\frac{\vep d\TY(R,1,n-1)}{\TY(R,2,n-1)+d\TY(R,1,n-1)},\q
\TYO(R,2,n-1)=\frac{\vep \TY(R,2,n-1)}{\TY(R,2,n-1)+d\TY(R,1,n-1)},
\label{R1}\\
&&\TYO(U,1,m)=\frac{\vep \TY(U,1,m)}{\TY(U,1,m)+d\TY(R,1,m)},\q
\TYO(R,1,m)=\frac{\vep d\TY(R,1,m)}{\TY(U,1,m)+d\TY(R,1,m)}\,\,(m<n),
\label{U1}\\
&&\TYO(U,l,m)=\frac{\vep\TY(U,l,m)}{\TY(\al,l+1,m)(d)}, \q
\TYO(V,l,m)=\frac{d \vep\TY(V,l,m)}{\TY(\al,l,m)(d)}, \q
\TYO(R,l,m)=\frac{d \vep^2\TY(R,l,m)}{\TY(\al,l,m)(d)\TY(\al,l+1,m)(d)}\,\,(l\geq 1),
\label{lm}\\
&&\TYO(X,l,m)=\frac{d \vep\TY(X,l,m)}{\TY(\al,l+1,m)(d)},\qq
{\TYO(X,l,m)}^*=\frac{\vep\TYS(X,l,m)}{\TY(\al,l,m)(d)}\q(l,m)\ne(1,n).
\label{XX*}
\end{eqnarray}
\end{lem}
{\sl Proof.}
Since we have $\ovl e_0^c \ovl e_0^d(y)=\ovl e_0^{cd}(y)$, we get 
$e_0^c e_0^d(x)=e_0^{cd}(x)$. Thus, we have $\TY((e_0^{cd}(x)),1,n-1)=
\TY((e_0^ce_0^d(x)),1,n-1)$ and then 
\begin{equation}
\TYO(R,2,n-1)+c\TYO(R,1,n-1)=\ovl\vep\frac{\TY(R,2,n-1)+cd\TY(R,1,n-1)}
{\TY(R,2,n-1)+d\TY(R,1,n-1)}.
\label{R2}
\end{equation}
By the formula $\vep_0(e_0^c(x))=c^{-1}\vep_0(x)$ in \eqref{ep0} and 
$\vep_0(x)=\TY(x,1,n)\vep(x)$ in \eqref{g0ep0}, 
one has $\vep(e_0^d(x))=\vep(x)$. 
Applying this to \eqref{R2}, one gets 
\begin{equation*}
\TYO(R,2,n-1)+c\TYO(R,1,n-1)=\vep\frac{\TY(R,2,n-1)+cd\TY(R,1,n-1)}
{\TY(R,2,n-1)+d\TY(R,1,n-1)}
\end{equation*}
and comparing the terms with or without $c$, one gets \eqref{R1}.

The formulae \eqref{U1} are the special case of \eqref{lm}. 
Hence, it is sufficient to show \eqref{lm}.
As considered above $\TY((e_0^{cd}(x)),l,m)=\TY((e_0^ce_0^d(x)),l,m)$,
we have 
\begin{equation}
\frac{\TY(U,l-1,m)+cd\TY(V,l,m)}{\TY(U,l,m)+cd\TY(V,l+1,m)}
=\frac{\TY(U,l-1,m)+d\TY(V,l,m)}{\TY(U,l,m)+d\TY(V,l+1,m)}
\frac{\TYO(U,l-1,m)+c\TYO(V,l,m)}{\TYO(U,l,m)+c\TYO(V,l+1,m)}
\label{UUVV}
\end{equation}
Multiplying the both sides of \eqref{UUVV}
for $l=1,2,\cd,l$, then one has
\[
\frac{\vep}{\TY(U,l,m)+cd\TY(V,l+1,m)}
=\frac{\vep}{\TY(U,l,m)+d\TY(V,l+1,m)}
\frac{\ovl\vep}{\TYO(U,l,m)+c\TYO(V,l+1,m)}.
\]
Since we have seen above $\vep(x)=\ovl\vep(x)$, one gets
\begin{equation}
{\TYO(U,l,m)+c\TYO(V,l+1,m)}=\vep\cdot
\frac{\TY(U,l,m)+cd\TY(V,l+1,m)}{\TY(U,l,m)+d\TY(V,l+1,m)}.
\end{equation}
Comparing the terms with or without $c$, 
we have first two formulae in \eqref{lm}.
The second formula in \eqref{lm} is obtained from
$\TYO(R,l,m)=\TYO(V,l+1,m)-\TYO(V,l,m)$. 
Here note that $\TY(\al,l,m)(d)=\TY(U,l-1,m)+d\TY(V,l,m)$.

Due to \eqref{U=U+XX2} and $\TY(R,l,m)=\TYS(X,l,m)\TY(X,l,m)$
one has 
$\frac{\TY(X,l-1,m+1)}{\TY(X,l,m)}
=\frac{\TY(U,l-1,m)-\TY(U,l-1,m+1)}{\TY(R,l,m)}$, and then 
sending $x$ to $e_0^d(x)$ and using the formula \eqref{lm} one has
\begin{equation}
\frac{\TYO(X,l-1,m+1)}{\TYO(X,l,m)}
=\frac{\TYO(U,l-1,m)-\TYO(U,l-1,m+1)}
{\TYO(R,l,m)}=\frac{\TY(\al,l+1,m)(d)(\TY(U,l-1,m)-\TY(U,l-1,m+1))}
{\TY(\al,l,m+1)(d)\TY(R,l,m)}
=\frac{\TY(\al,l+1,m)(d)\TY(X,l-1,m+1)}
{\TY(\al,l,m+1)(d)\TY(X,l,m)},
\end{equation}
for the second equality we use 
\begin{equation*}
\TYO(U,l-1,m)-\TYO(U,l-1,m+1)=
\frac{\TY(U,l-1,m)}{\TY(\al,l,m)}-\frac{\TY(U,l-1,m+1)}{\TY(\al,l,m+1)}=
\frac{d\ep(\TY(U,l-1,m)-\TY(U,l-1,m+1))}{\TY(\al,l,m)\TY(\al,l,m+1)},
\end{equation*}
which is obtained by using $\TY(V,l,m)=\vep-\TY(U,l-1,m)$.
Thus, one can find
\begin{equation}
\frac{\TYO(X,l,m)}{\TY(X,l,m)}
=\frac{\TY(\al,l,m+1)}{\TY(\al,l+1,m)}
\frac{\TYO(X,l-1,m+1)}{\TY(X,l-1,m+1)}
=\cd=
\frac{\TY(\al,l,m+1)}{\TY(\al,l+1,m)}
\frac{\TY(\al,l-1,m+2)}{\TY(\al,l,m+1)}
\cd
\frac{\TY(\al,2,l+m-1)}{\TY(\al,3,l+m-2)}
\frac{\TYO(X,1,l+m-1)}{\TY(X,1,l+m-1)}
\label{Xalalal}
\end{equation}
Here applying $\TY(X,1,l+m-1)=\frac{\TY(x,1,l+m-1)}{\TY(x,1,n)}$, 
$\TYO(X,1,l+m-1)=\frac{\TYO(x,1,l+m-1)}{\TYO(x,1,n)}$, 
$\TYO(x,1,n)=\frac{\TY(x,1,n)}{c}$ and $\TYO(x,1,l+m-1)
=\frac{\vep\TY(x,1,l+m-1)}{\TY(\al,2,l+m-1)}$ to \eqref{Xalalal}, we obtain
the first one of \eqref{XX*}. The second one of \eqref{XX*}
is obtained by the similar method using the formula:
\[
\frac{\TYS(X,l,m)}{\TYS(X,l-1,m+1)}
=\frac{\TY(U,l-1,m)-\TY(U,l-1,m+1)}{\TY(R,l-1,m+1)},
\]
which follows from \eqref{U=U+XX2} and 
$\TY(R,l-1,m+1)=\TYS(X,l-1,m+1)\TY(X,l-1,m+1)$ also.\qed

\medskip

Here, we define the subset of $P_1[n,k]$:
${}^1P_1[n,k]:=\{p\in P_1[n,k]|p{\hbox{ is through }}(1,n-1)\}$ and 
${}^2P_1[n,k]:=\{p\in P_1[n,k]|p{\hbox{ is through }}(2,n-1)\}$.
Note that since all paths must go through $(1,n-1)$ or $(2,n-1)$, 
we have $P_1[n,k]={}^1P_1[n,k]\sqcup {}^2P_1[n,k]$.
For $j=1,2$, set 
\[
{}^j\TY({LW},l,m):=\TY({LW},l,m)\cap {}^j P_1[n,k],\q
{}^j\TY({{P_1[,k]}},l,m):=\TY({{P_1[n,k]}},l,m)\cap {}^j P_1[n,k],
\]
where $W=P,A,B$.
For example, an element in ${}^1\TY(LA,l,m)$ is a path
through $(1,n-1)$ and above $(l,m)$.
Here for $j=1,2$ we define:
\begin{eqnarray*}
&&\jTY(X,l,m):=\sum_{p\in \jTY({{P_1[n,k]}},l,m)}x(p),\qq
\jTY(R,l,m):=\sum_{p\in \jTY(LP,l,m)}x(p),\\
&&\jTY(U,l,m):=\sum_{p\in \jTY(LA,l,m)}x(p),\qq
\jTY(V,l,m):=\sum_{p\in \jTY(LB,l,m)}x(p),
\end{eqnarray*}
where note that $\TY(W,l,m)=\aTY(W,l,m)+\bTY(W,l,m)$
$(W=X,R,U,V)$.
Let us see the formulae related to these functions.
\begin{lem}\label{formula-1}
For a function $f(x)$ on $\cV_1$ let us denote $\ovl f(x):=f(e^d_0(x))$ as before.
We get the following formula:
\begin{eqnarray}
&&\aTYO(U,l,m)=\frac{d\vep^2\cdot\aTY(U,l,m)}
{(\TY(R,2,n-1)+d\TY(R,1,n-1))\TY(\al,l+1,m)(d)},\label{1U}\\
&&\bTYO(V,l,m)=\frac{d\vep^2\cdot \bTY(V,l,m)}
{(\TY(R,2,n-1)+d\TY(R,1,n-1))\TY(\al,l,m)(d)},\label{2V}\\
&&\aTYO(X,l,m)=\frac{d\vep \TY(\al,l,m)(d)}
{(\TY(R,2,n-1)+d\TY(R,1,n-1))\TYS(X,l,m)(d)}
\left(\frac{\aTY(U,l-1,m)}{\TY(\al,l,m)(d)}
-\frac{\aTY(U,l,m)}{\TY(\al,l+1,m)(d)}
\right),\label{XUU1}\\
&&\aTYO(X,l,m)=\frac{d\vep \TY(\al,l+1,m-1)(d)}
{(\TY(R,2,n-1)+d\TY(R,1,n-1))\TYS(X,l+1,m-1)(d)}
\left(\frac{\aTY(U,l,m-1)}{\TY(\al,l+1,m-1)(d)}
-\frac{\aTY(U,l,m)}{\TY(\al,l+1,m)(d)}
\right),\label{XUU2}
\end{eqnarray}
where note that \eqref{XUU2} is not immediate from \eqref{XUU1}.
\end{lem}
{\sl Proof.}
First, let us show \eqref{XUU1} and  \eqref{XUU2} 
by descending induction on $m$.
One can easily know that the formula \eqref{x=x+x} is valid for 
$\aTY(X,l,m)$ and $\bTY(X,l,m)$, namely,
\begin{equation}
\jTY(X,l,m)=\jTY(X,l-1,m+1)+\frac{\TY(x,l,m)}{\TY(x,l,m+1)}\,
\jTY(X,l,m+1)\q(j=1,2).
\label{x=x+x12}
\end{equation}
Then, we also have 
\begin{equation}
\aTYO(X,l,m)=\aTYO(X,l-1,m+1)+\frac{\TYO(x,l,m)}{\TYO(x,l,m+1)}\,
\aTYO(X,l,m+1).
\label{1barx=x+x}
\end{equation}
Applying the induction hypothesis on the right hand side of 
\eqref{1barx=x+x},
more precisely, substituting \eqref{XUU2} for $\aTYO(X,l-1,m+1)$ and 
\eqref{XUU1} for  $\aTYO(X,l,m+1)$ to \eqref{x=x+x12}, 
we obtain 
\begin{equation}
\aTYO(X,l,m)=\frac{\vep d\TY(\al,l,m)}{\TYS(X,l,m)\TYS(X,l,m+1)\xi}
(A-B+C-D),
\label{ABCD}
\end{equation}
where 
\begin{eqnarray*}
&&A:=\frac{\TY(x,l,m)\TY(\al,l+1,m+1)\TYS(X,l,m)\aTY(U,l-1,m+1)}
{\TY(x,l,m+1)\TY(\al,l,m+1)\TY(\al,l+1,m)},\q
B:=\frac{\TY(x,l,m)\TYS(X,l,m)\aTY(U,l,m+1)}{\TY(x,l,m+1)\TY(\al,l+1,m)},\\
&&C:=\frac{\TYS(X,l,m+1)\aTY(U,l-1,m)}{\TY(\al,l,m)},\q
D:=\frac{\TYS(X,l,m+1)\aTY(U,l-1,m+1)}{\TY(\al,l,m+1)},\q
\xi:=\TY(R,2,n-1)+d\TY(R,1,n-1).
\end{eqnarray*}
By using \eqref{*x=x+x}, we get 
\begin{equation}
A-D=\frac{\aTY(U,l-1,m+1)}{\TY(\al,l,m+1)\TY(\al,l+1,m)}
\left(\TYS(X,l,m+1)(\TY(\al,l+1,m+1)-\TY(\al,l+1,m))-\TYS(X,l+1,m)
\TY(\al,l+1,m+1)\right).\label{ad}
\end{equation}
The following  formulae are immediate from the definitions:
\begin{eqnarray}
&&\TY(\al,l+1,m+1)-\TY(\al,l+1,m)=(d-1)\TYS(X,l+1,m)\TY(X,l,m+1),\\
&&(d-1)\TY(R,l,m+1)-\TY(\al,l+1,m+1)=\TY(\al,l,m+1).
-\TY(U,l-1,m+1)-d\TY(V,l,m+1).
\end{eqnarray}
Applying these formulae to \eqref{ad}, we get 
\begin{eqnarray}
A-D&=&\frac{\aTY(U,l-1,m+1)\TYS(X,l+1,m)}{\TY(\al,l,m+1)\TY(\al,l+1,m)}
\left((d-1)\TY(R,l,m+1)-\TY(\al,l+1,m+1)\right)=
=\frac{\aTY(U,l-1,m+1)\TYS(X,l+1,m)}{\TY(\al,l,m+1)\TY(\al,l+1,m)}
(-\TY(U,l-1,m+1)-d\TY(V,l,m+1))\nn \\
&=&-\frac{\aTY(U,l-1,m+1)\TYS(X,l+1,m)}{\TY(\al,l,m+1)\TY(\al,l+1,m)}
\TY(\al,l,m+1)=
-\frac{\aTY(U,l-1,m+1)\TYS(X,l+1,m)}{\TY(\al,l+1,m)}.\label{ad2}
\end{eqnarray}
Here applying \eqref{*x=x+x} and 
$\aTY(R,l,m+1)=\aTY(X,l,m+1)\TYS(X,l,m+1)$ to \eqref{ad2}, we obtain
\begin{equation}
(A-D)-B=-\frac{1}{\TY(\al,l+1,m)}\left(\frac{\TY(x,l,m)}{\TY(x,l,m+1)}
\TYS(X,l,m)\aTY(U,l,m+1)+\TYS(X,l+1,m)\aTY(U,l-1,m+1)\right)
=-\frac{\TYS(X,l,m+1)}{\TY(\al,l+1,m)}(\TY(X,l,m+1)\TYS(X,l+1,m)+
\aTY(U,l,m+1))=-\frac{\TYS(X,l,m+1)}{\TY(\al,l+1,m)}\aTY(U,l,m).
\label{abd}
\end{equation}
Thus, finally applying this to \eqref{ABCD} we get
\begin{equation}
\aTYO(X,l,m)=\frac{\vep d\TY(\al,l,m)}{\TYS(X,l,m)\TYS(X,l,m+1)\xi}
(A-B+C-D)
=\frac{\vep d\TY(\al,l,m)}{\TYS(X,l,m)\xi}
\left(\frac{\aTY(U,l-1,m)}{\TY(\al,l,m)}-\frac{\aTY(U,l,m)}{\TY(\al,l+1,m)}
\right), \label{ABCD2}
\end{equation}
which means \eqref{XUU1}. 
So, let us show \eqref{XUU2} from \eqref{XUU1}. Let $P$ (resp. $Q$)
be the right-hand side of \eqref{XUU1} (resp. \eqref{XUU2}).
Then, one has
\begin{equation}
\frac{\xi}{d\vep}(P-Q)=
\frac{\aTY(U,l,m+1)}{\TYS(X,l+1,m-1)}
-\frac{\aTY(U,l-1,m)}{\TYS(X,l,m)}
-\frac{\TY(\al,l+1,m-1)\aTY(U,l,m)}{\TYS(X,l+1,m-1)\TY(\al,l+1,m)}
+\frac{\TY(\al,l,m)\aTY(U,l,m)}{\TYS(X,l,m)\TY(\al,l+1,m)}.
\label{PQ}
\end{equation}
\begin{lem}
We get the following formula:
\begin{eqnarray}
&&\frac{\aTY(U,l,m+1)}{\TYS(X,l+1,m-1)}
-\frac{\aTY(U,l-1,m)}{\TYS(X,l,m)}
=\aTY(U,l,m)\left(\frac{1}{\TYS(X,l+1,m-1)}-\frac{1}{\TYS(X,l,m)}\right),
\label{PC}\\
&&\TY(\al,l,m)\TYS(X,l+1,m-1)-\TY(\al,l+1,m-1)\TYS(X,l,m)
=\TY(\al,l+1,m)(\TYS(X,l+1,m-1)-\TYS(X,l,m)).
\label{QD}
\end{eqnarray}
\end{lem}
{\sl Proof.}
By considering similarly to \eqref{U=U+XX2}, we get the formula
\begin{equation}
\aTY(U,l,m)=\aTY(U,l,m+1)+\TYS(X,l+1,m)\cdot \aTY(X,l,m+1).
\label{1U1UXX}
\end{equation}
It follows from \eqref{1U1UXX} that we get \eqref{PC}. 

Using the formula \eqref{U=U+XX2}, $\TY(U,l-1,m)=\TY(U,l,m)+\TY(R,l,m)$
and $\vep=\TY(U,l-1,m)+\TY(V,l,m)$ we have
\begin{eqnarray*}
{\hbox{ L.H.S of }}\eqref{QD}&&=
\TY(U,l-1,m)\TYS(X,l+1,m-1)-\TY(U,l,m-1)\TYS(X,l,m)
+d(\TY(V,l,m)\TYS(X,l+1,m-1)-\TY(V,l+1,m-1)\TYS(X,l,m))\\
&&=(1-d)(\TY(U,l-1,nm)\TYS(X,l+1,m-1)-\TY(U,l,m-1)\TYS(X,l,m))
+d\vep(\TYS(X,l+1,m-1)-\TYS(X,l,m))\\
&&=((1-d)\TY(U,l,m)+d\vep)(\TYS(X,l+1,m-1)-\TYS(X,l,m))
=\TY(\al,l+1,m)(\TYS(X,l+1,m-1)-\TYS(X,l,m)).\qed
\end{eqnarray*}
Thus, applying \eqref{PC} and \eqref{QD} to \eqref{PQ} we obtain 
$P=Q$, which implies \eqref{XUU2}. Now, the induction proceeds and then
we complete the proof of \eqref{XUU1} and \eqref{XUU2}.

Next, let us show \eqref{1U}.
By \eqref{XX*} and \eqref{XUU1} one gets
\begin{equation}
\aTY(\ovl R,l,m)=\TYOS(X,l,m)\cdot\aTYO(X,l,m)
=\frac{d\vep^2}{\TY(R,2,n-1)+d\TY(R,1,n-1)}
\left(\frac{\aTY(U,l-1,m)}{\TY(\al,l,m)(d)}
-\frac{\aTY(U,l,m)}{\TY(\al,l+1,m)(d)}\right).\label{RUU}
\end{equation}
Applying \eqref{RUU} to 
$\aTYO(U,l,m)=\sum_{j=\min(n+1-m,k)}^{l+1}\aTYO(R,j,m)$, one gets \eqref{1U}.

Finally, let us show \eqref{2V}.
One has
\begin{equation}
\bTYO(V,l,m)=\TYO(V,l,m)-\aTYO(V,l,m)=\frac{d\vep}{\xi\TY(\al,l,m)(d)}
\left(\xi\TY(V,l,m)+\vep\aTY(U,l-1,m)-\TY(\al,l,m)(d)\TY(R,1,n-1)\right)
=:\frac{d\vep}{\xi\TY(\al,l,m)(d)}\cdot Z,
\label{VVUR}
\end{equation}
where $\xi=\TY(R,2,n-1)+d\TY(R,1,n-1)$. Then,
\begin{align*}
Z&=\TY(V,l,m)(\xi-d\TY(R,1,n-1))+\vep\aTY(U,l-1,m)-\TY(R,1,n-1)\TY(U,l-1,m)\\
&=\TY(V,l,m)\TY(R,2,n-1)+(\vep-\TY(R,1,n-1)\aTY(U,l-1,m)-
\TY(R,1,n-1)\bTY(U,l-1,m)\\
&=\TY(V,l,m)\TY(R,2,n-1)+\TY(R,2,n-1)\aTY(U,l-1,m)-\TY(R,1,n-1)\bTY(U,l-1,m)
=\TY(R,2,n-1)(\TY(V,l,m)+\aTY(U,l-1,m))-\TY(R,1,n-1)\bTY(U,l-1,m)\\
&=\TY(R,2,n-1)(\bTY(V,l,m)+\TY(R,1,n-1))-\TY(R,1,n-1)\bTY(U,l-1,m)=\vep\cdot \bTY(V,l,m),
\end{align*}
which shows \eqref{2V}.\qed
\subsection{Proof of \eqref{0n0}}

Now, we shall show \eqref{0n0} by virtue 
of the formulae in the previous subsection.
Set $\til x=e_0^ce_n^{cd}e_0^d(x)$ and $x''=e_n^de_0^{cd}e_n^c(x)$ 
where note that 
$x,\til x, x''$ mean $v_1(x),v_1(\til x),v_1(x'')\in \cV_1$.
It suffices to show that $\TY(\til x,l,m)=\TY(x'',l,m)$ for any $(l,m)$,
which will be done by using the induction on $(l,m)\in L_1[n,k]$ 
according to the total
order $\prec$:$(l,m)\prec(l',m')\Leftrightarrow $ 
$l<l'$, or if $l=l'$, $m>m'$.

Since we have $\TY(e_n^c(x),1,n)=c\TY(x,1,n)$ and 
$\TY(e_0^d(x),1,n)=\frac{\TY(x,1,n)}{d}$, we obtain
$\TY(\til x,1,n)=\TY(x'',1,n)=\TY(x,1,n)$, 
which shows the case $(l,m)=(1,n)$.
Next, let us see the case $(l,m)=(1,n-1)$.
Let $\TY(\al,l,m)=\TY(\al,l,m)(c)=\TY(U,l-1,m)+c\TY(V,l,m)$ be as above.
Thus, we have $\TY(\al,1,n-1)=\TY(R,1,n-1)+\TY(R,2,n-1)$ and 
$\TY(\al,2,n-1)=c\TY(R,1,n-1)+\TY(R,2,n-1)$, where note that 
$\TY(R,1,n-1)$ depends on $\TY(x,1,n)$ but not does $\TY(R,2,n-1)$ so.
Hence, we find that
\begin{equation}
\TY(\ovl \al,1,n-1)(c)=\frac{1}{cd}\TY(\ovl R,1,n-1)+\TY(\ovl R,2,n-1),\q
\TY(\ovl \al,2,n-1)(c)=\frac{1}{d}\TY(\ovl R,1,n-1)+\TY(\ovl R,2,n-1).
\label{1n-1}
\end{equation}
Here note that 
\begin{equation}
\TY(\til x,l,m)=\TY(\ovl x,l,m)\left.\frac{\TY(\ovl\al,l,m)(c)}
{\TY(\ovl\al,l+1,m)(c)}\right|_{\TY(\ovl x,1,n)\to cd\TY(\ovl x,1,n)}
=\TY(x,l,m)\frac{\TY(\al,l,m)(d)}{\TY(\al,l+1,m)(c)}\cdot
\left.\frac{\TY(\ovl\al,l,m)(c)}
{\TY(\ovl\al,l+1,m)(d)}\right|_{
\TY(\ovl x,1,n)\to cd\TY(\ovl x,1,n)}
\label{tilx-al}
\end{equation}
Thus, it follows from \eqref{R1}, \eqref{1n-1} and \eqref{tilx-al} that 
\begin{eqnarray}
\TY(\til x, 1,n-1)&=&
\TY(x,1,n-1)\frac{\TY(\al,1,n-1)(d)}{\TY(\al,2,n-1)(d)}
\cdot \frac{\frac{1}{cd}\TY(\ovl R,1,n-1)+\TY(\ovl R,2,n-1)}
{\frac{1}{d}\TY(\ovl R,1,n-1)+\TY(\ovl R,2,n-1)}\nn
 \\
&=&\TY(x,1,n-1)\frac{\TY(R,1,n-1)+\TY(R,2,n-2)}{d\TY(R,1,n-1)+\TY(R,2,n-1)}
\cdot \frac{\frac{\vep d}{cd}\TY(R,1,n-1)+\vep \TY(R,2,n-1)}
{\frac{\vep d}{d}\TY(R,1,n-1)+\vep \TY(R,2,n-1)}
=\TY(x,1,n-1)\frac{\frac{1}{c}\TY(R,1,n-1)+\TY(R,2,n-1)}
{d\TY(R,1,n-1)+\TY(R,2,n-1)}.
\label{tilx1n-1}
\end{eqnarray}
On the other hand, we have 
\begin{eqnarray*}
\TY(x'',1,n-1)=\TY((e_n^de_0^{cd}e_n^c(x)),1,n-1)
=\TY((e_0^{cd}e_n^c(x)),1,n-1)
=\TY(x,1,n-1)\left.\frac{\TY(R,1,n-1)+\TY(R,2,n-1)}
{cd\TY(R,1,n-1)+\TY(R,2,n-1)}\right|_{\TY(x,1,n)\to c \TY(x,1,n)}
=\TY(x,1,n-1)\frac{\frac{1}{c}\TY(R,1,n-1)+\TY(R,2,n-1)}
{d\TY(R,1,n-1)+\TY(R,2,n-1)}.
\end{eqnarray*}
Then, by this and \eqref{tilx1n-1} we find that 
$\TY(\til x,1,n-1)=\TY(x'',1,n-1)$.

Let us see the case $(1,m)$ with $m<n-1$. 
Here for a function $f(x)$ 
we shall denote $f(e_n^c x)$ (resp. $f(e_n^{cd}e_0^d x)$)
by $\what f(x)$ (resp. $\ovl{\ovl f}(x)$).
The following are immediate from the explicit action of $e_n^{cd}$:
\begin{equation}
\TY(\ovl{\ovl R},1,m)=\frac{1}{cd}\TY(\ovl R,1,m),\qq
\ovl{\ovl \vep}=\frac{1}{cd}\TY(\ovl R,1,n-1)+\TY(\ovl R,2,n-1).\label{R-ep}
\end{equation}
It follows from \eqref{U1}, \eqref{tilx-al} and \eqref{R-ep} that
\begin{eqnarray}
\TY(\til x,1,m)&=&
\TY(\ovl x,1,m)\frac{\TY(\ovl{\ovl \al},1,m)(c)}
{\TY(\ovl{\ovl \al},2,m)(c)}=\TY(x,1,m)\frac{\TY(\al,1,m)(d)}{\TY(\al,2,m)(d)}
\left(\frac{\TY(\ovl{\ovl U},1,m)+\TY(\ovl{\ovl R},1,m)}
{\TY(\ovl{\ovl U},1,m)+c \TY(\ovl{\ovl R},1,m)}\right)
\nn \\
&=&\TY(x,1,m)\frac{\vep}{\TY(U,1,m)+d\TY(R,1,m)}\cdot
\frac{\vep \frac{\frac{1}{c}\TY(R,1,n-1)+\TY(R,2,n-2)}
{d\TY(R,1,n-1)+\TY(R,2,n-1)}}
{\ovl{\ovl \vep}-\TY(\ovl{\ovl R},1,m)+c\TY(\ovl{\ovl R},1,m)}
=\TY(x,1,m)\frac{\frac{1}{c}\TY(R,1,n-1)+\TY(R,2,n-1)}
{(d-\frac{1}{c})\TY(R,1,m)+\frac{1}{c}\TY(R,1,n-1)+\TY(R,2,n-1)}.\label{tilx1m}
\end{eqnarray}
We also have the following:
\begin{equation}
\TY(\what R,1,m)=\frac{1}{c}\TY(R,1,m),\qq
\what \vep=\frac{1}{c}\TY(R,1,n-1)+\TY(R,2,n-1),\qq
\TY(\what U,1,m)=\what\vep-\TY(\what R,1,m).\label{hats}
\end{equation}
Using these formulae, we get
\begin{equation*}
\TY(x'',1,m)=\TY(x,1,m)\frac{\TY(\what \al,1,m)(cd)}{\TY(\what\al,2,m)(cd)}
=\TY(x,1,m)\frac{\what\vep}{\TY(\what U,1,m)+cd\TY(\what R,1,m)}
=\TY(x,1,m)\frac{\frac{1}{c}\TY(R,1,n-1)+\TY(R,2,n-1)}{(d-\frac{1}{c})\TY(R,1,m)
+\frac{1}{c}\TY(R,1,n-1)+\TY(R,2,n-1)},
\end{equation*}
which coincides with \eqref{tilx1m} and then we have
$\TY(\til x, 1,m)=\TY(x'',1,m)$ for $m<n-1$.

Now, let us see the general case $(l,m)$ with $l>1$.
We have 
\[
\TY(\til x,l,m)=\TY(\ovl{\ovl x},l,m)\frac{\TY(\ovl{\ovl \al},l,m)(c)}
{\TY(\ovl{\ovl \al},l+1,m)(c)}=\TY(x,l,m)
\frac{\TY(\al,l,m)(d)}{\TY(\al,l+1,m)(d)}
\cdot\frac{\TY(\ovl{\ovl \al},l,m)(c)}{\TY(\ovl{\ovl \al},l,m)(c)},\qq
\TY(\what x,l,m)=
\TY(x,l,m)\frac{\TY(\what{\al},l,m)(cd)}{\TY(\what{\al},l+1,m)(cd)}.
\]
Therefore, to obtain $\TY(\til x,l,m)=\TY(x'',l,m)$ we shall show 
\begin{equation}
\frac{\TY(\al,l,m)(d)}{\TY(\al, l+1,m)(d)}\cdot 
\frac{\TY(\ovl{\ovl \al},l,m)(c)}{\TY(\ovl{\ovl \al},l+1,m)(c)}=
\frac{\TY(\what \al,l,m)(cd)}{\TY(\what \al,l+1,m)(cd)}.
\label{alalal}
\end{equation}
To show \eqref{alalal}, it suffices to show that there exists a function 
$\Theta$ which does not depend on $(l,m)$ and possibly depends on 
$x,c,d$ such that 
\begin{equation}
\TY(\al,l,m)(d)\cdot 
\TY(\ovl{\ovl \al},l,m)(c)=
\Theta\TY(\what \al,l,m)(cd).
\label{Theta}
\end{equation}
\begin{lem}
We have the following formulae:
\begin{eqnarray}
\TY(\what U,l,m)&=&\frac{1}{c}\aTY(\what U,l,m)+\bTY(\what U,l,m),
\label{hat-U}\\
\TY(\what V,l,m)&=&\frac{1}{c}\aTY(\what V,l,m)+\bTY(\what V,l,m),
\label{hat-V}\\
\bTY(\ovl U,l,m)&=&\frac{\vep}{\TY(\al,l+1,m)(d)}
\left(\TY(U,l,m)-\frac{d\vep}{\xi}\aTY(U,l,m)\right),\label{2U}\\
\aTY(V,l,m)&=&\frac{d\vep^2}{\xi}\left(\frac{\TY(R,1,n-1)}{\vep}-\frac{\aTY(U,l-1,m)}{\TY(\al,l,m)(d)}\right).\label{1V}
\end{eqnarray}
where $\xi:=\TY(\al,2,n-1)(d)=\TY(R,2,n-1)+d\TY(R,1,n-1)$.
\end{lem}
{\sl Proof.}
The formula \eqref{2U}
is an immediate consequence of \eqref{U1}, \eqref{1U}
and the fact $\TY(U,l,m)=\aTY(U,l,m)+\bTY(U,l,m)$. 
The formula \eqref{1V} is obtained similarly.

We know that $\aTY(U,l,m)\cdot \TY(x,1,n)$ 
and $\bTY(U,l,m)$ do not depend on $\TY(x,1,n)$, which induces \eqref{hat-U}.
By considering similarly, we also obtain \eqref{hat-V}.\qed

Due to this lemma and Lemma \ref{formula-1}, we obtain 
\begin{eqnarray*}
\TY(\ovl{\ovl \al},l,m)(c)&=&
\TY(\ovl{\ovl U},l-1,m)+c\TY(\ovl{\ovl V},l,m)=
(\aTY(\ovl{\ovl U},l-1,m)+\bTY(\ovl{\ovl U},l-1,m))+c
(\aTY(\ovl{\ovl V},l,m)+\bTY(\ovl{\ovl V},l,m))\\
&=&\frac{1}{cd}\aTY({\ovl U},l-1,m)+\bTY({\ovl U},l-1,m)
+c(\frac{1}{cd}\aTY({\ovl V},l,m)+\bTY({\ovl V},l,m))\\
&=&\frac{\vep^2}{\xi\TY(\al,l,m)(d)}
\left(\frac{1}{c}\aTY(U,l-1,m)+\bTY(U,l-1,m)+d\aTY(V,l,m)
+cd\cdot \bTY(V,l,m)\right)\\
&=&\frac{\vep^2}{\xi\TY(\al,l,m)(d)}\left(\TY(\what U,l-1,m)+cd
\TY(\what V,l,m)\right)=\frac{\vep^2}{\xi\TY(\al,l,m)(d)}\TY(\what\al,l,m)(cd),
\end{eqnarray*}
which shows \eqref{Theta} and $\Theta=\frac{\vep^2}{\xi}$ 
then \eqref{alalal}. Hence, we completed to prove
$\TY(\til x,l,m)=\TY(x'',l,m)$ and then the relation
\eqref{0n0}. \qed

\renewcommand{\thesection}{\arabic{section}}
\section{Ultra-discretization of $\cV(\TY(A,1,n))$}
\setcounter{equation}{0}
\renewcommand{\theequation}{\thesection.\arabic{equation}}

For basic notions of crystals, coherent family of 
perfect crystals and their limit we refer the reader to \cite{KKM} 
(See also \cite{KMN1,KMN2}).

\subsection{Crystal $B^{k,\ify}$}
We review the $\TY(A,1,n)$-crystal $B^{k,\ify}$ following \cite{OSS}.
\label{kify}
We define the crystal $B^{k,\ify}$ by setting
\begin{equation}
B^{k,\ify}:=\left\{
(b_{ji})_{1\leq j\leq k,\\\syl j\leq i\leq j+k'}\left|
b_{ji}\in\bbZ, \sum_{i=j}^{j+k'}b_{ji}=0{\hbox{ for any }}j\right.\right\},
\end{equation}
where $k'=n+1-k$. We also define $b_\ify$ in $B^{k,\ify}$ to be an element 
whose all entries are 0.
For $b=(b_{ji})\in B^{k,\ify}$, $\vep_i(b), \vp_i(b)$, $\eit(b)$ and 
$\fit(b)$  for $i= 1,2,\cd, n$ are defined as follows:\\
For $i\in I\setminus\{0\}$, let us set $\beta=\max(0,i-k')$ and 
$\gamma=\min(k,i)$. Fix $b=(b_{ji})\in \kify$.
For $c\in \bbZ$ with $\beta<c\leq \gamma$ and $i\in I$, 
let us set 
\begin{equation}
\Gamma_i(c):=\sum_{\beta<j<c}(b_{ji}-b_{j+1,i+1}), \q
\Gamma_{i,\min}:=\min\{\Gamma_i(c)|\beta<c\leq \gamma\}.
\label{gammac}
\end{equation}
Define
\begin{equation}
c_0:=\min\{c|\beta<c\leq \gamma,\,\,\Gamma_i(c)=\Gamma_{i,\min}\},\q
c_1:=\max\{c|\beta<c\leq \gamma,\,\,\Gamma_i(c)=\Gamma_{i,\min}\}.
\label{c0c1}
\end{equation}
For $i\in I_0$ the functions $\vep_i(b)$ and $\vp_i(b)$ on $\kify$
are defined by
\begin{equation}
\vep_i(b):=\sum_{\beta\leq j<c_0}(b_{j+1,i+1}-b_{ji}),\q
\vp_i(b):=\sum_{c_1 \leq j\leq \gamma}(b_{ji}-b_{j+1,i+1}).
\label{vepvp}
\end{equation}
The Kashiwara operators $\eit$ and $\fit$ are given as follows: 
set $b'=(b'_{p,q})=\eit(b)$ and 
$b''=(b''_{p,q})=\fit(b)$ $(i\in I_0)$ 
and then
\begin{equation}
b'_{p,q}=b_{p,q}-\del_{p,c_0}\del_{q,i+1}+\del_{p,c_0}\del_{q,i},\q
b''_{p,q}=b_{p,q}-\del_{p,c_1}\del_{q,i}+\del_{p,c_1}\del_{q,i+1}.
\label{eitfit}
\end{equation}
Note that $\fit^{-1}=\eit$.
Next, let us define $\vep_0$, $\vp_0$, $\til e_0$ and $\til f_0$:
Set 
\begin{equation}
C:=\{(c_0,c_1,\cd,c_k)\in\bbN^{k+1}|1=c_0<c_1<\cd c_k\leq n+1\}
\label{CCC}
\end{equation}
and define a partial order $c \preceq c'$ on $C$ by $c_j\leq c'_j$ for any $j$.
For fixed $b=(b_{ji})\in \kify$, let $\Del_b:C\to \bbZ_{\geq 0}$ by 
\[
\Del_b(c)=\sum_{j=1}^k \sum_{c_{j-1}<i<c_j}b_{ji},
\]
and set 
$C^b_{\min}:=\{c\in C|\Del_b(c){\hbox{ is minimal.}}\}$.
By \cite[5.3.]{OSS}, we obtain the following lemma:
\begin{lem}
Fix $b=(b_{ji})\in \kify$.
\begin{enumerate}
\item 
There exists unique $\TY(c,e,)\in C$ such that 
\begin{eqnarray*}
&\Del_b(\TY(c,e,))\leq \Del_b(c)&{\rm if}\,\, \TY(c,e,)\preceq c,\\
&\Del_b(\TY(c,e,))< \Del_b(c)&{\rm if}\,\, \TY(c,e,)\not\preceq c,
\end{eqnarray*}
for any $c\in C$.
\item
There exists unique $\TY(c,f,)\in C$ such that 
\begin{eqnarray*}
&\Del_b(\TY(c,f,))\leq \Del_b(c)&{\rm if}\,\, c\preceq \TY(c,f,),\\
&\Del_b(\TY(c,f,))< \Del_b(c)&{\rm if}\,\, c\not\preceq \TY(c,f,),
\end{eqnarray*}
for any $c\in C$.
\end{enumerate}
\end{lem}
Note that $c^{(e)}, c^{(f)}\in \TY(C,b,\min)$ and then $\Del_b(c^{(e)})=\Del_b(c^{(f)})$.
Here, using these $\TY(c,e,)$ and $\TY(c,f,)$, for $b=(b_{ji})\in\kify$ 
we define
\begin{eqnarray}
&&\vep_0(b):=-b_{k,n+1}-\Del_b(\TY(c,e,)),\q
\vp_0(b):=-b_{1,1}-\Del_b(\TY(c,f,)),
\label{vep0}\\
&&b'_{ji}:=b_{ji}-\del_{i,\TY(c,e,j-1)}+\del_{i,\TY(c,e,j)},\label{e00}\\
&&b''_{ji}:=b_{ji}-\del_{i,\TY(c,f,j)}+\del_{i,\TY(c,f,j-1)},\label{f00}
\end{eqnarray}
where $b'=\til e_0(b)$ and $b''=\til f_0(b)$. Note that 
$\til f_0^{-1}=\til e_0$.
Then, by these definitions the crystal $\kify$ becomes 
a limit of coherent family of perfect crystals
$\{B^{k,l}\}_{l\geq l}$ (\cite{KKM,OSS}).
Note that $\wt_0(b)=\vp_0(b)-\vep_0(b)=-b_{1,1}+b_{k,n+1}$ since 
$\Del_b(c^{(e)})=\Del_b(c^{(f)})$. 
\subsection{Crystal structure on ${\mathcal UD}(\cV(\TY(A,1,n)))$}

Let us see the crystal structure on $(\UD(\cV(\TY(A,1,n))),
\UD(e_{i,v_1}),\UD(\gamma_i\circ v_1),\UD(\vep_i\circ v_1))_{i\in I}$, where 
$v_1:(\bbC^\times)^{kk'}\to \cV(\TY(A,1,n))(=\cV_1)$ 
is the positive structure on $\cV(\TY(A,1,n))$.
For a horizontal strip $s=(l,m)\frac{\qq}{\qq}(l,m+1)$ on $L_1[n,k]$, 
let $uw(s):=\TY(x,l,m)-\TY(x,l,m+1)$ and for a path $p\in P_1[n,k]$
define $ux(p):=\sum_{\hbox{\small s:strip in }p}uw(s)$.
\begin{pro}\label{ex-gcry}
The explicit forms of $(\UD(\cV(\TY(A,1,n))),
\UD(e_{i,v_1}),\UD(\gamma_i\circ v_1),\UD(\vep_i\circ v_1))$
are as follows:
$\UD(\cV(\TY(A,1,n)))=\bbZ^{kk'}$ where 
we denote the coordinate $\UD(\TY(x,i,j))$ also by $\TY(x,i,j)$.
\begin{eqnarray}
&&\UD(\gamma_0\circ v_1)=-\TY(x,1,n)-\TY(x,k,1),\\
&&\UD(\gamma_i\circ v_1)=
2\sum_{j=\aa}^\bb\TY(x,j,i)-\sum_{j=\aa}^{\bb+1}\TY(x,j,i)
-\sum_{j=\aa-1}^\bb\TY(x,j,i),\label{udgi}\\
&&\UD(\vep_0\circ v_1)=\TY(x,1,n)+\max_{p\in P_1[n,k]}
(ux(p)),\label{udv0}\\
&&\UD(\vep_i\circ v_1)=\max_{\aa\leq l\leq \bb}(\TYO(D,l,i)),\label{udvi}\\
&&\TY(\UD(e_{0,v_1})^d(x),l,m)=\TY(x,l,m)+\max(\max_{p\in \TY({LA[n,k]},l-1,m)}
(ux(p)),d+\max_{p\in \TY({LB[n,k]},l,m)}(ux(p)))\\
&&\qq\qq -\max(\max_{p\in \TY({LA[n,k]},l,m)}(ux(p)),
d+\max_{p\in \TY({LB[n,k]},l+1,m)}(ux(p)))\q(d\in\bbZ),\nn \\
&&\TY(\UD(e_{i,v_1})^d(x),l,m)=\begin{cases}
\displaystyle \TY(x,l,m)+\max(\max_{\aa\leq p\leq l-1}(\TYO(D,p,i)),\max_{l\leq p\leq \bb}
(d+\TYO(D,p,i)))&\\
\displaystyle
-\max(\max_{\aa\leq p\leq l}(\TYO(D,p,i)),\max_{l+1\leq p\leq \bb}
(d+\TYO(D,p,i))))&{\rm if}\,\, m=i,\\
\TY(x,l,m),&{\rm otherwise,}\label{udex}
\end{cases}
\end{eqnarray}
where $d\in\bbZ$, $1\leq i\leq n$, $(l,m)\in L_1[n,k]$, 
$\aa$ and $\bb$ are as in \ref{v1v2} and 
\begin{equation}
\TYO(D,l,i)=-\UD(\TY(D,l,i))
=-\TY(x,l,i)-2\sum_{j=l+1}^\bb
\TY(x,j,i)+\sum_{j=l+1}^{\bb+1}\TY(x,j,i-1)+\sum_{j=l}^\bb\TY(x,j,i+1).
\label{bdli}
\end{equation}
\end{pro}
\nd
{\sl Proof.}
All the formula are obtained trivially applying  the formulae
\[
\UD(X\cdot Y^{\pm})=\UD(X)\pm \UD(Y),\qq
\UD(X+Y)=\max(X,Y),
\]
to the explicit forms of 
$e_i^d,\gamma_i,\vep_i$ ($i\in I$) as in
\eqref{Dli}, \eqref{gamma-ep1}, \eqref{eic1}, \eqref{g0ep0} and 
\eqref{e00}.
\qed

In the sequel, we shall
denote $\UD(e_{i,v_1})$, $\UD(\gamma_i\circ v_1)$ and 
$\UD(\vep_i\circ v_1)$ by
$\eit'$, $\wt'_i$ and $\vep'_i$ respectively. 
\subsection{Isomorphism between ${\mathcal UD}(\cV)$ and $\kify$}

We shall show that the crystal $\kify$ and $\bbZ^{kk'}=
\UD(\TY(A,1,n))$ $(k'=n-k+1)$ are 
isomorphic to each other. 
Define the map $\Omega:\bbZ^{kk'}=\UD(\TY(A,1,n))\to\kify$ as follows:
Set $(b_{ji})=\Omega(x)$:
\begin{equation}
b_{ji}:=\TY(x,k-j+1,i)-\TY(x,k-j+1,i-1),
\label{omegaji}
\end{equation}
where we understand that 
$\TY(x,l,m)$ is 0 if $(l,m)$ is out of the lattice $L_1[n,k]$.
The following theorem gives an affirmative answer to the conjecture 
for $\TY(A,1,n)$ and generic Dynkin index $k$.
\begin{thm}\label{ud-thm}
The map $\Omega$ gives an isomorphism of crystals between 
$\kify$ and $\bbZ^{kk'}=\UD(\TY(A,1,n))$.
\end{thm}
\nd
{\sl Proof.}
First, let us see the well-definedness of the map $\Omega$, that is,
we may show that 
\begin{equation}
\sum_{i=j}^{j+k'}\Omega(x)_{ji}=0\q{\hbox{ for any }}j.
\label{omega0}
\end{equation}
The left hand-side of \eqref{omega0} is 
\[
(\TY(x,j'',j)-\TY(x,j'',j-1))-(\TY(x,j'',j+1)-\TY(x,j'',j))+\cd+
(\TY(x,j'',j+k')-\TY(x,j'',j+k'-1))=
\TY(x,j'',j+k')-\TY(x,j'',j-1)=0\,\,(j''=k-j+1),
\]
since $(j'',j+k')$ and $(j'',j-1)$ are out of $L_1[n,k]$. 
Indeed, 
$j''+(j+k')=n+2>n+1$ and $j''+j-1=k$. 

Next, we shall see the bijectivity of $\Omega$. It can be done easily by 
giving the inverse of $\Omega$:
\[
\TY(\Omega^{-1}(b),j,i)=b_{j'',j''}+b_{j'',j''+1}+\cd+b_{j'',i-1}+b_{j'',i}.
\]
It is clear to show that this is an inverse of $\Omega$. Indeed, 
denoting the right hand-side of the above formula by $\TY(\Theta,j,i)$, 
we obtain
\[
\TY(\Theta,j'',i)-\TY(\Theta,j'',i-1)=b_{j,i}.
\]
Thus, the map $\Omega$ is bijective.

Finally, let us show for any $x\in\bbZ^{kk'}(=\UD(\cV(\TY(A,1,n))))$ and 
$i\in I$,
\begin{enumerate}
\item
$\wt'_i(x)=\wt_i\circ\Omega(x)$.
\item
$\vep'_i(x)=\vep_i\circ\Omega(x)$.
\item
$\eit'(x)=\Omega^{-1}\circ \eit\circ \Omega(x)$.
\end{enumerate}
(i)
First note that $\Del_b(c^{(e)})=\Del_b(c^{(f)})=\min\{\Del_b(c)|c\in C\}$
for $b\in\kify$. 
Thus, we get $\wt_0(b)=\vp_0(b)-\vep_0(b)=-b_{1,1}+b_{k,n+1}$ and then 
\[
\wt_0\circ\Omega(x)=-\TY(x,k,1)-\TY(x,1,n)=\wt'_0(x).
\]
Suppose that $i\ne0$.
By the definition of $c_0$ and $c_1$ in \eqref{c0c1}, we have 
$\sum_{j=c_0}^{c_1-1}(b_{j,i}-b_{j+1,i+1})=0$. 
By using this, we obtain
\begin{eqnarray}
\wt_i(b)=\vp_i(b)-\vep_i(b)=
\sum_{\beta\leq j<c_0}B_{ji}+\sum_{c_1\leq j\leq \gamma}B_{ji}=
\sum_{\beta\leq j\leq \gamma}B_{ji},
\end{eqnarray}
where $B_{ji}:=b_{j,i}-b_{j+1,i+1}$.
We have
\begin{eqnarray*}
&&\wt_i\circ\Omega(x)=\sum_{\beta\leq j\leq\gamma}\TY(x,k-j+1,i)-\TY(x,k-j+1,i-1)
-\TY(x,k-j,i+1)+\TY(x,k-j,i).
\end{eqnarray*}
Note that 
\begin{eqnarray*}
&&\beta=\max(0,i-k')=\max(0,i-n+k-1)=k+\max(-k,i-n-1)=k-\min(k,n-i+1)=k-\bb,\\
&&\gamma=\min(k,i)=k+1+\min(-1,i-k-1)=k+1-\max(1,k-i+1)=k+1-\aa.
\end{eqnarray*}
Hence, we get
\[
\sum_{\beta\leq j\leq \gamma}\TY(x,k-j+1,i)+\TY(x,k-j,i)
=\sum_{j=k-\gamma+1}^{k-\beta+1}\TY(x,j,i)+\sum_{j=k-\gamma}^{k-\beta}\TY(x,j,i)
=2\sum_{j=\aa}^\bb\TY(x,j,i),
\]
since $\TY(x,k-\gamma,i)=\TY(x,\aa-1,i)=0$ 
and $\TY(x,k-\beta+1,)=\TY(x,\bb+1,i)=0$. Therefore, by \eqref{udgi} 
\begin{eqnarray*}
\wt_i\circ\Omega(x)&=&
\sum_{\beta\leq j\leq \gamma}(\TY(x,k-j+1,i)-\TY(x,k-j+1,i-1)
-\TY(x,k-j,i+1)+\TY(x,k-j,i))
\\&=&2\sum_{j=\aa}^\bb\TY(x,j,i)-\sum_{j=\aa}^{\bb+1}\TY(x,j,i-1)-
\sum_{j=\aa-1}^\bb\TY(x,j,i+1)=\wt'_i(x).
\end{eqnarray*}

\nd (ii)
Let us see the relation between the set of paths $P_1[n,k]$ and 
the set $C$ as in \eqref{CCC}.
\begin{lem}\label{path-c}
There exists one-to-one correspondence between $P_1[n,k]$ and 
$C$ defined as follows:
For $c=(1,c_1,\cd,c_{k-1},n+1)\in C$, let us define a path 
$p=(s_1,s_2\cd,s_k)$ by 
\[
s_j=(k-j+1,c_{j-1})\frac{\qq}{\qq}(k-j+1,c_j-1),
\]
where $s_j$ is consecutive horizontal strips as in \eqref{strips} with 
length $c_j-c_{j-1}-1$. 
\end{lem}
The proof is immediate. Let us denote this correspondence 
by $\pi:C\to P_1[n,k]$.

Using this correspondence, one has for $(b_{ji})=\Omega(x)$ and $c\in C$
\begin{eqnarray}
\qq \Del_{\Omega(x)}(c)&=&\sum_{j=1}^k\sum_{c_{j-1}<i<c_j}b_{ji}=
\sum_{j=1}^k\sum_{c_{j-1}<i<c_j}(\TY(x,k-j+1,i)-\TY(x,k-j+1,i))
\nn \\
&=&\sum_{j=1}^k(\TY(x,k-j+1,c_j -1)-\TY(x,k-j+1,c_{j-1}))
=-\sum_{j=1}^k(\TY(x,k-j+1,c_{j-1} )-\TY(x,k-j+1,c_{j}-1))
=-ux(\pi(c)).\label{xp-Del}
\end{eqnarray}
Then, for $(b_{ji})=\Omega(x)$ by \eqref{udv0} we get
\begin{eqnarray*}
&&\vep_0\circ\Omega(x)=
-b_{k,n+1}-\Del_{\Omega(x)}(c^{(e)})=\TY(x,1,n)-\min_{c\in C}\{\Del_{\Omega(x)}(c)\}
\\
&&=\TY(x,1,n)-\min_{c\in C}\{-ux(\pi(c)\}
=\TY(x,1,n)+\max_{p\in P_1[n,k]}\{\-ux(p)\}=\vep'_0(x).
\end{eqnarray*}
Next, let us see the case $i\ne0$. By the definition of $c_0$ 
as in \eqref{c0c1}, for $b=(b_{ji})\in\kify$ one has
\begin{eqnarray*}
\vep_i(b)&=&\sum_{\beta\leq j<c_0}(b_{j+1,i+1}-b_{ji})
=b_{\beta,+1,i+1}-b_{\beta,i}-\sum_{\beta<j<c_0}(b_{j,i}-b_{j+1,i+1})\\
&=&b_{\beta,+1,i+1}-b_{\beta,i}-\Gamma_i(c_0)=
b_{\beta,+1,i+1}-b_{\beta,i}-\min_{\beta<c\leq \gamma}
\{\sum_{\beta< j<c}(b_{ji}-b_{j+1,i+1})\}\\
&=&b_{\beta+1,i+1}-b_{\beta,i}+\max_{\beta<c\leq \gamma}
\{\sum_{j=\beta+1}^{c-1}(b_{j+1,i+1}-b_{ji})\}.
\end{eqnarray*}
Due to this formula, we get
\begin{eqnarray*}
\vep_i\circ\Omega(x)&=&
\TY(x,k-\beta,i+1)-\TY(x,k-\beta,i)-\TY(x,k-\beta+1,i)
+\TY(x,k-\beta+1,i-1)\\
&&+\max_{\beta<c\leq \gamma}
\{\sum_{j=\beta+1}^{c-1}\TY(x,k-j,i+1)-\TY(x,k-j,i)-\TY(x,k-j+1,i)
+\TY(x,k-j+1,i-1)\}.
\end{eqnarray*}
Here one gets
\begin{eqnarray*}
&&\max_{\beta<c\leq \gamma}
\{\sum_{j=\beta+1}^{c-1}\TY(x,k-j,i+1)-\TY(x,k-j,i)-\TY(x,k-j+1,i)
+\TY(x,k-j+1,i-1)\}
=\max_{\beta<c\leq \gamma}
\{\sum_{j=k-c+1}^{\bb-1}\TY(x,j,i+1)-\TY(x,j,i)-\TY(x,j+1,i)
+\TY(x,j+1,i-1)\}\\
&&=\max_{\aa\leq l\leq \bb}
\{\sum_{j=l}^{\bb-1}\TY(x,j,i+1)-\TY(x,j,i)-\TY(x,j+1,i)
+\TY(x,j+1,i-1)\},
\end{eqnarray*}
where the last equality follows from 
that setting  $l=k-c+1$ we have  $l=k-\gamma+1=\aa$ (resp. $l=k-\beta=\bb$)
if $c=\gamma$ (resp. $c=\beta+1$).
Therefore, by \eqref{udvi} one has 
\begin{eqnarray*}
\vep_i\circ\Omega(x)&=&
\TY(x,\bb,i+1)-\TY(x,\bb,i)-\TY(x,\bb+1,i)
+\TY(x,\bb+1,i-1)
+\max_{\aa\leq l\leq \bb}
\{\sum_{j=l}^{\bb-1}\TY(x,j,i+1)-\TY(x,j,i)-\TY(x,j+1,i)
+\TY(x,j+1,i-1)\}\\
&=&
\max_{\aa\leq l\leq \bb}\{
-\TY(x,l,i)-2\sum_{j=l+1}^\bb
\TY(x,j,i)+\sum_{j=l+1}^{\bb+1}\TY(x,j,i-1)+\sum_{j=l}^\bb\TY(x,j,i+1)\}
=\max_{\aa\leq l\leq \bb}\{\TYO(D,l,i)\}=\vep'_i(x).
\end{eqnarray*}
(iii)
Suppose $i=0$.
Let us show $\til e_0=\Omega\circ \til e'_0\circ\Omega^{-1}$.
The following lemma is immediate from the definitions of paths.
\begin{lem}
For $c=(c_0,c_1,\cd, c_{k+1})\in C$, set $p=\pi(c)\in P_1[n,k]$. 
\begin{enumerate}
\item
Path $p$ is above $(l,m)$ iff $c_{k-l}>m$.
\item
Path $p$ is below $(l,m)$ iff $c_{k-l+1}\leq m$.
\item
Path $p$ is through $(l,m)$ iff $c_{k-l}\leq m<c_{k-l+1}$.
\end{enumerate}
\end{lem}
Since $\TY(U,l-1,m)$ (resp. $\TY(V,l,m)$) is 
the total sum of all weight of paths above $(l-1,m)$ (resp. below $(l,m)$),
by using \eqref{xp-Del} and the above lemma, we obtain 
\begin{eqnarray*}
\UD(\TY(U,l-1,m)(\Omega^{-1}(b))=
\max_{c:c_{k-l+1}>m}\{-\Del_b(c)\},\q
\UD(\TY(V,l,m)(\Omega^{-1}(b))=
\max_{c:c_{k-l+1}\leq m}\{-\Del_b(c)\}.
\end{eqnarray*}
Then we obtain
\begin{equation}
\UD(\TY(\al,l,m)(d))(\Omega^{-1}(b))=\max\left(
\max_{c:c_{k-l+1}>m}\{-\Del_b(c)\}, d+
\max_{c:c_{k-l+1}\leq m}\{-\Del_b(c)\}\right)=:
\TY(\bbA,l,m)(d).
\label{ud-al}
\end{equation}
Thus, 
\begin{equation}
\UD(e_0)^d(\Omega^{-1}(b))^{(l)}_m=
\TY(x,l,m)(\Omega^{-1}(b))+\TY(\bbA,l,m)(d)-\TY(\bbA,l+1,m)(d),
\label{ud-e0}
\end{equation}
Therefore, we get
\begin{eqnarray}
&&(\Omega\circ\UD(e_0)^d(\Omega^{-1}(b)))_{ji}=
\UD(e_0)^d(\Omega^{-1}(b))^{(k-j+1)}_i-
\UD(e_0)^d(\Omega^{-1}(b))^{(k-j+1)}_{i-1}\nn
\\
&&=\TY(x,k-j+1,i)(\Omega^{-1}(b))+\TY(\bbA,k-j+1,i)(d)-\TY(\bbA,k-j+2,i)(d)
-\TY(x,k-j+1,i-1)(\Omega^{-1}(b))-\TY(\bbA,k-j+1,i-1)(d)+\TY(\bbA,k-j+2,i-1)(d)\nn \\
&&=b_{ji}+\TY(\bbA,k-j+1,i)(d)-\TY(\bbA,k-j+2,i)(d)
-\TY(\bbA,k-j+1,i-1)(d)+\TY(\bbA,k-j+2,i-1)(d),
\label{4A}
\end{eqnarray}
Let us denote the last formula of \eqref{4A} by $\bbP(d)$.
Now we know that for $b=(b_{ji})\in\kify$ it suffices to show 
\begin{equation}
\til e_0(b)_{ji}=b_{ji}-\del_{i,\ovl c_{j-1}}+\del_{i,\ovl c_j}
=\bbP(1),
\label{fin-e0}
\end{equation}
where set $c^{(e)}=(\ovl c_0,\ovl c_1,\cd,\ovl c_{k+1})$.
One has 
\begin{equation}
\TY(\bbA,k-j+1,i)(1)=-\min(\min_{c:c_j>i}\{\Del_b(c)\},
-1+\min_{c:c_j\leq i}\{\Del_b(c)\}).
\label{bbP1}
\end{equation}
We shall see the cases :
(I) $\ovl c_j=i$. 
(II) $\ovl c_j>i$.
(III) $\ovl c_j<i$.
\begin{enumerate}
\item[(I)]
Suppose $\ovl c_j=i$. 
One has $\min_{c:c_j\leq i}\{\Del_b(c)\}=\Del_b(c^{(e)})$.
Let $c'\in C$ be an element such that 
$\min_{c:c_j> i}\{\Del_b(c)\}=\Del_b(c')$ and $c'_j>i$. 
By the definition of $c^{(e)}$, we find that 
if $c'\succ c^{(e)}$, $\Del_b(c')\geq \Del_b(c^{(e)})$ and 
if $c'\not\succ c^{(e)}$, $\Del_b(c')> \Del_b(c^{(e)})$.
Then, in any case we get 
$\TY(\bbA,k-j+1,i)(1)=1-\Del_b(c^{(e)})$ and 
$\TY(\bbA,k-j+1,i-1)(1)=-\Del_b(c^{(e)})$, and then
\[
\TY(\bbA,k-j+1,i)(1)-\TY(\bbA,k-j+1,i-1)(1)=1.
\]
\item[(II)] Suppose $\ovl c_j>i$.
One has $\min_{c:c_j> i}\{\Del_b(c)\}=\Del_b(c^{(e)})$. 
Let $c'\in C$ be an element such that 
$\min_{c:c_j\leq i}\{\Del_b(c)\}=\Del_b(c')$ and $c'_j\leq i$, which means 
$c'\not\succ c^{(e)}$ and then $\Del_b(c')> \Del_b(c^{(e)})$.
Therefore, we obtain $\TY(\bbA,k-j+1,i)(1)=-\Del_b(c^{(e)})$. 
Next, we also have $\min_{c:c_j> i-1}\{\Del_b(c)\}=\Del_b(c^{(e)})$
since $\ovl c_j>i>i-1$. 
Let $c''\in C$ be an element such that 
$\min_{c:c_j\leq i-1}\{\Del_b(c)\}=\Del_b(c'')$ and $c''_j\leq i-1$, which means 
$c''\not\succ c^{(e)}$ and then $\Del_b(c'')> \Del_b(c^{(e)})$.
Thus, we have $\TY(\bbA,k-j+1,i-1)(1)=-\Del_b(c^{(e)})$. 
Therefore, we obtain $\TY(\bbA,k-j+1,i)(1)-\TY(\bbA,k-j+1,i-1)(1)=0$.
\item[(III)]
Suppose $\ovl c_j<i$.
One gets $\min_{c:c_j\leq i}\{\Del_b(c)\}=\Del_b(c^{(e)})$. 
Let $c'\in C$ be an element such that 
$\min_{c:c_j> i}\{\Del_b(c)\}=\Del_b(c')$ and $c'_j>i$.
If $c'\succ c^{(e)}$, then $\Del_b(c')\geq \Del_b(c^{(e)})$, and 
if $c'\not\succ c^{(e)}$, then $\Del_b(c')> \Del_b(c^{(e)})$. In both cases, we have
$\TY(\bbA,k-j,i)(1)=1-\Del_b(c^{(e)})$. 
Next, we also have $\min_{c:c_j\leq i-1}\{\Del_b(c)\}=\Del_b(c^{(e)})$
since $\ovl c_j<i$. 
It is easy to get $\min_{c:c_j> i-1}\{\Del_b(c)\}\geq \Del_b(c^{(e)})$.
Thus, we have $\TY(\bbA,k-j,i-1)(1)=1-\Del_b(c^{(e)})$ and 
then we obtain $\TY(\bbA,k-j,i)(1)-\TY(\bbA,k-j,i-1)(1)=0$.
By the result from  (I),(II) and (III), we find that 
\begin{equation}
\TY(\bbA,k-j,i)(1)-\TY(\bbA,k-j,i-1)(1)=\del_{i,\ovl c_j},\qq
\TY(\bbA,k-j+1,i)(1)-\TY(\bbA,k-j+1,i-1)(1)=\del_{i,\ovl c_{j-1}}.
\label{AAdel}
\end{equation}
Substituting \eqref{AAdel} in \eqref{4A} and setting $d=1$, we obtain
\[
(\Omega\circ e'_0\circ\Omega^{-1}(b))_{ji}=
(\Omega\circ \UD(e_0)\circ\Omega^{-1}(b))_{ji}=b_{ji}
+\del_{i,\ovl c_j}-\del_{i,\ovl c_{j-1}}
=\til e_0(b)_{ji}.
\]
\end{enumerate}
Here note that if we set $d=-1$, we obtain the action of $\til f_0$.

Finally, we consider the case $i\in\{1,2,\cd,n\}$.
The action of Kashiwara operators $\eit$ and $\fit$ on $\kify$
are  given in \eqref{eitfit}.
By the explicit form of $\TY(D,l,m)$ in \eqref{Dli}, for 
$x\in\bbZ^{kk'}$ one gets
\[
\UD(\TY(D,p,i))(x)
=\TY(x,p,i)-\TY(x,p,i+1)-\TY(x,\bb+1,i-1)+
\sum_{r=p+1}^\bb 2\TY(x,r,i)-\TY(x,r,i-1)-\TY(x,r,i+1),
\]
and then using $\TY(\Omega^{-1}(b),j,i)=\sum_{s=k-j+1}^i b_{k-j+1,s}$, 
for $b=(b_{ji})\in\kify$ one has
\begin{equation}
\UD(\TY(D,p,i))(\Omega^{-1}(b))
=\sum_{r=p+1}^{\bb+1}b_{k-r+1,i}-\sum_{r=p}^\bb b_{k-r+1,i+1}=:\cD_{p,i}
\end{equation}
Thus, it follows from \eqref{udex} that for $d\in\bbZ$, 
\begin{eqnarray*}
\UD(e_i)^d(\Omega^{-1}(b))^{(l)}_m
&=&\sum_{s=1}^m b_{k-l+1,s}+\del_{i,m}
\left(\max\left(\max_{\aa\leq p\leq l-1}(-\cD_{p,i}),d+
\max_{l\leq p\leq \bb}(-\cD_{p,i})\right)\right.\\
&&-\left.\max\left(\max_{\aa\leq p\leq l}(-\cD_{p,i}),d+
\max_{l+1\leq p\leq \bb}(-\cD_{p,i})\right)\right).
\end{eqnarray*}
And then, 
\begin{eqnarray} 
&&\Omega\circ\UD(e_i)^d(\Omega^{-1}(b))_{l,m}
=\UD(e_i)^d(\Omega^{-1}(b))^{(k-l+1)}_m
-\UD(e_i)^d(\Omega^{-1}(b))^{(k-l+1)}_{m-1}
\label{oeoi}\\
\hspace{-30pt}&&=b_{l,m}+\del_{i,m}
\left(\max\left(\max_{\aa\leq  p\leq l'-1}(-\cD_{p,i}),d+
\max_{l'\leq  p\leq \bb}(-\cD_{p,i})\right)
-\max\left(\max_{\aa\leq  p\leq l'}(-\cD_{p,i}),d+
\max_{l'+1\leq  p\leq\bb}(-\cD_{p,i})\right)\right)\nn \\
\hspace{-30pt}&&-\del_{i,m-1}
\left(\max\left(\max_{\aa\leq  p\leq l'-1}(-\cD_{p,i}),d+
\max_{l'\leq  p\leq \bb}(-\cD_{p,i})\right)
-\max\left(\max_{\aa\leq  p\leq l'}(-\cD_{p,i}),d+
\max_{l'+1\leq  p\leq \bb}(-\cD_{p,i})\right)\right),\nn
\end{eqnarray}
where $l'=k-l+1$.
Here, recalling $\beta=k-\bb$ and $\gamma=k+1-\aa$ as above, we get
\begin{eqnarray*}
&&\max_{\aa\leq  p\leq l'-1}(-\cD_{p,i})
=-\min_{l+1 \leq  p\leq \gamma}(\cD_{p,i})
=-\min_{l+1 \leq  p\leq \gamma}(\sum_{r=\beta}^{p-1} (b_{r,i}-b_{r+1,i+1})
\\
&&\max_{l'\leq  p\leq \bb}(-\cD_{p,i})
=-\min_{\beta+1 \leq  p\leq l}(\cD_{p,i})
=-\min_{\beta+1 \leq  p\leq l}(\sum_{r=\beta}^{p-1} (b_{r,i}-b_{r+1,i+1}).
\end{eqnarray*}
Since $\sum_{r=\beta}^p (b_{r,i}-b_{r+1,i+1}
=\Gamma_i(p)+(b_{\beta,i}-b_{\beta+1,i+1})$ and 
in general, $\min(A+x,B+x)-\min(C+x,D+x)=\min(A,B)-\min(C,D)$, 
by setting $d=1$ we obtain
\begin{eqnarray*}
\hspace{-30pt}&&\Omega\circ\UD(e_i)(\Omega^{-1}(b))_{l,m}\\
\hspace{-30pt}&&=b_{l,m}+(\del_{i,m-1}-\del_{i,m})\times
\left(\min\left(\min_{l+1 \leq  p\leq \gamma}(\Gamma_i(p)),-1+
\min_{\beta+1 \leq  p\leq l}(\Gamma_i(p))\right)\right.\\
\hspace{-30pt}&&\left.
-\min\left(\min_{l \leq  p\leq \gamma}(\Gamma_i(p)),-1+
\min_{\beta+1 \leq  p\leq l-1}(\Gamma_i(p))\right)\right)
=:b_{i,m}+(\del_{i,m-1}-\del_{i,m})\times G.
\end{eqnarray*}
 Now, to finalize the proof, let us show 
\[
G=-\del_{l,c_0},
\]
where $c_0$ is as in \eqref{c0c1}.
By the definition of $c_0$ in\eqref{c0c1}, we find that 
\begin{enumerate}
\item If $l>c_0$, $G=(-1+\Gamma_i(c_0))-(-1+\Gamma_i(c_0))=0$.
\item If $l<c_0$, $G=\Gamma_i(c_0)-\Gamma_i(c_0)=0$.
\item If $l=c_0$, $G=(-1+\Gamma_i(c_0))-\Gamma_i(c_0)=-1$.
\end{enumerate}
Thus, we have $G=-\del_{l,c_0}$. Therefore, 
for $i\in\{1,2,\cd,n\}$ we obtain
\[
\Omega\circ\UD(e_i)(\Omega^{-1}(b))_{l,m}
=b_{l,m}+\del_{m,i}\del_{l,c_0}-\del_{m,i+1}\del_{l,c_0}=
\eit(b)_{l,m}.
\]
Here, we have completed the proof of Theorem \ref{ud-thm}.\qed

\renewcommand{\thesection}{\arabic{section}}
\section{Affine Weyl Group Action}
\setcounter{equation}{0}
\renewcommand{\theequation}{\thesection.\arabic{equation}}

In this section, as an application of the affine geometric crystal structure
on $\cV(\TY(A,1,n))$ and its ultra-discretization, we shall describe 
a birational action of affine Weyl group and a piecewise-linear 
action of affine Weyl group explicitly.

\subsection{Birational action of affine Weyl group $W(\TY(A,1,n))$}
Let $(X,\{e_i\}_{i\in I},\{\gamma_i,\}_{i\in I},
\{\vep_i\}_{i\in I})$ be a geometric crystal of $G$(or $\ge$).
\begin{thm}[\cite{BK}]\label{weyl}
Define $s_i(x):=e_i^{{\gamma_i(x)}^{-1}}(x)$ for $x\in X$.
Then, $\lan s_i|i\in I\ran$ gives a birational actions on $X$
of Weyl group $W$ associated with $G$ (or $\ge$).
\end{thm}
Indeed, we can easily find that $s_i^2={\rm id}_X$, 
$s_is_j=s_js_i$ if $a_{ij}=a_{ji}=0$,  and 
the braid relation $s_is_js_i=s_js_is_j$ holds if $a_{ij}=a_{ji}=-1$.
Since we have obtained the $\TY(A,1,n)$-geometric crystal structure on 
$\cV(\TY(A,1,n))$, we get the birational action of the 
affine Weyl group $W(\TY(A,1,n))$ on $\cV(\TY(A,1,n))$. 
Its explicit forms are as follows:
\begin{thm}\label{aff-weyl}
Let us set
\begin{equation}
\TY(F,p,m):=
\frac{\TY(x,p,m)(\TY(x,\aa,m)\cd\TY(x,p-1,m))^2}
{\TY(x,\aa+1,m-1)\cd\TY(x,p,m-1)\TY(x,\aa,m+1)\cd\TY(x,p-1,m+1)}.
\label{Fpm}
\end{equation}
Then, we obtain the action of $s_i$:
\begin{eqnarray}
&&\TY(s_i(x),l,m)=\begin{cases}
\displaystyle \TY(x,l,i)
\frac{
\displaystyle\sum_{p=a}^{l-1}\TY(F,p,i)+
\displaystyle\sum_{p=l}^b\frac{1}{\TY(D,p,i)}}
{\displaystyle
\sum_{p=a}^{l}\TY(F,p,i)+
\displaystyle\sum_{p=l+1}^b\frac{1}{\TY(D,p,i)}},&{\rm if}\,\,m=i,\\
\TY(x,l,m),&{\rm otherwise},
\end{cases}(i=1,\cd,n),
\label{si}\\
&&\TY(s_0(x),l,m)=\begin{cases}
\displaystyle \TY(x,l,m)\frac{\TY(U,l-1,m)+\TY(x,1,n)\TY(x,k,1)
\TY(V,l,m)}{\TY(U,l,m)+\TY(x,1,n)\TY(x,k,1)
\TY(V,l+1,m)},&{\rm if}\,\,(l,m)\ne(1,n),\\
\displaystyle \frac{1}{\TY(x,k,1)}&{\rm if}\,\,(l,m)=(1,n).
\end{cases}
\label{s0}
\end{eqnarray}
\end{thm}
{\sl Proof.}
We know that $\gamma_i(x)=\frac{\TY(D,\aa-1,i)}{\TY(x,\aa-1,i)}$ 
for $i=1,\cd,n$ by \eqref{gamma-ep1},
and $\gamma_0(x)=\frac{1}{\TY(x,1,n)\TY(x,k,1)}$ for $i=0$
by \eqref{g0ep0}. Then,
due to Theorem \ref{weyl} substituting these formulae to \eqref{eic1}
and \eqref{e0}, 
we get \eqref{si} and \eqref{s0}.\qed

\subsection{Piecewise-linear action of the affine Weyl group 
$W(\TY(A,1,n))$}
We shall describe the action of $W(\TY(A,1,n))$ 
on the crystal $B^{k,\ify}$.

As an application of Theorem \ref{weyl} and ultra-discretization of
positive geometric crystals, we obtain the following easily.
\begin{thm}\label{ud-weyl}
Let $(X,\{e_i\}_{i\in I},\{\gamma_i\}_{i\in I},
\{\vep_i\}_{i\in I})$ be a positive geometric crystal of $G$(or $\ge$)
and $(B,\{\til e_i\}_{i\in I},\{\wt_i\}_{i\in I},\{\vep_i\}_{i\in I})$ its
ultra-discretized crystal.
Define $\til s_i(x):=\eit^{-\wt_i(b)}(b)$ for $b\in B$.
Then, $\lan \til s_i|i\in I\ran$ gives a piecewise-linear 
actions on $B$ of Weyl group $W$ associated with $G$ (or $\ge$).
\end{thm}
\begin{thm}\label{ud-aff-weyl}
Let $\beta$, $\gamma$ and $\Del_b(c)$ be as in \ref{kify}.
For $p\in\{\beta+1,\cd,\gamma\}$ and $b=(b_{j,i})\in \kify$ set 
\begin{equation}
\til\Gamma_i(p):=\sum_{p\leq j\leq \gamma}(b_{j,i}-b_{j+1,i+1}),\q
\TY(\Del,1,b)(c):=\Del_b(c)+b_{11},\q
\TY(\Del,k,b)(c):=\Del_b(c)+b_{k,n+1}.
\end{equation}
Then the action of $\til s_i$ is as follows:
\begin{eqnarray}
&&\til s_i(b)_{l,m}=
b_{l,m}+(\del_{i,m-1}-\del_{i,m})\times\label{ud-si}\\
&&\times \left(
\min\left(\min_{l+1\leq p\leq \gamma}\til\Gamma_i(p),
\min_{\beta+1\leq p\leq l}\til\Gamma_i(p)\right)
-\min\left(\min_{l\leq p\leq \gamma}\til\Gamma_i(p),
\min_{\beta+1\leq p\leq l-1}\til\Gamma_i(p)\right)
\right)\q(i=1,2,\cd,n),\nn \\
&&\til s_0(b)_{l,m}=b_{l,m}+\min\left(\min_{c_{l-1}>m}\TY(\Del,1,b)(c),
\min_{c_{l-1}\leq m}\TY(\Del,k,b)(c)\right)
-\min\left(\min_{c_{l}>m}\TY(\Del,1,b)(c),
\min_{c_{l}\leq m}\TY(\Del,k,b)(c)\right)\label{ud-s0}\\
&&+\min\left(\min_{c_{l}\geq m}\TY(\Del,1,b)(c),
\min_{c_{l}< m}\TY(\Del,k,b)(c)\right)
-\min\left(\min_{c_{l-1}\geq m}\TY(\Del,1,b)(c),
\min_{c_{l-1}<m}\TY(\Del,k,b)(c)\right).\nn
\end{eqnarray}
\end{thm}
{\sl Proof.}
As for the cases $i\ne0$, substituting $d=-\wt_i(b)=
\sum_{\beta+1\leq j\leq \gamma}(b_{ji}-b_{j+1,i+1})$ 
for \eqref{oeoi} we obtain the formula \eqref{ud-si}.

As for the $i=0$-case, substituting $d=-\wt_0(b)=
b_{11}-b_{k,n+1}$ for \eqref{4A}, we get the formula \eqref{ud-s0}.
\qed

\bibliographystyle{amsalpha}

\end{document}